\documentclass[12pt]{amsart}
\font\emailfont=cmtt10

\usepackage{amsmath,amsthm,amsfonts,amscd,flafter,epsf}

\newcommand{\Tors}{\mathrm{Tors}}

\newcommand{\HF}{HF}

\newtheorem{theorem}{Theorem}[section]
\newtheorem{prop}[theorem]{Proposition}

\newtheorem{lemma}[theorem]{Lemma}

\newtheorem{defn}[theorem]{Definition}

\def\endproof{\relax\ifmmode\expandafter\endproofmath\else
  \unskip\nobreak\hfil\penalty50\hskip.75em\hbox{}\nobreak\hfil\bull
  {\parfillskip=0pt \finalhyphendemerits=0 \bigbreak}\fi}
\def\endproofmath$${\eqno\bull$$\bigbreak}
\def\bull{\vbox{\hrule\hbox{\vrule\kern3pt\vbox{\kern6pt}\kern3pt\vrule}\hrule}}

\newcounter{bean}
\newcommand{\Q}{\mathbb{Q}}
\newcommand{\R}{\mathbb{R}}

\newcommand{\Z}{\mathbb{Z}}

\newcommand{\OneHalf}{\frac{1}{2}}

\newcommand{\Zmod}[1]{\Z/{#1}\Z}

\newcommand{\Coker}{\mathrm{Coker}}

\newcommand{\cm}{\cdot}

\newcommand{\Nbd}[1]{{\mathrm{nd}}(#1)}
\newcommand{\nbd}[1]{\Nbd{#1}}
\newcommand{\CDisk}{\mathbb D}

\newcommand{\ModSWfour}{\mathcal{M}}
\newcommand{\ModFlow}{\ModSWfour}

\newcommand{\SpinC}{{\mathrm{Spin}}^c}

\newcommand{\goesto}{\mapsto}

\newcommand{\DBar}{\overline{\partial}}

\newcommand\Wedge{\Lambda}

\newcommand\Hom{\mathrm{Hom}}

\newcommand\abuts\Rightarrow
\newcommand\Sym{\mathrm{Sym}}

\newcommand{\ModDeg}{\mathcal N}
\newcommand{\UnparModDeg}{\widehat \ModDeg}

\newcommand\RelSpinC{\underline{\SpinC}}
\newcommand\relspinc{\underline{\spinc}}

\newcommand\Filt{\mathcal F}

\newcommand\x{\mathbf x}
\newcommand\w{\mathbf w}
\newcommand\z{\mathbf z}

\newcommand\q{\mathbf q}
\newcommand\y{\mathbf y}

\newcommand\ModSphere{\ModFlow\left({\mathbb S}\longrightarrow 
\Sym^{g-1}(\Sigma_{1})\times \Sym^2(\Sigma_{2})\right)}
\newcommand\ModSpheres\ModSphere
\newcommand\CF{CF}

\newcommand\CFa{\widehat{CF}}

\newcommand\CFm{\CF^-}

\newcommand\HFm{\HF^-}

\newcommand\HFa{\widehat{HF}}

\newcommand\gr{\mathrm{gr}}
\newcommand\Mas{\mu}
\newcommand\UnparModSp{\widehat \ModSp}
\newcommand\UnparModFlow\UnparModSp
\newcommand\Mod\ModSp

\newcommand{\cald}{{\mathcal D}}

\newcommand\PD{\mathrm{PD}}

\newcommand{\spinc}{\mathfrak s}

\newcommand\ModMaps{\mathcal M}
\newcommand\ModSp\ModMaps

\newcommand\Ta{{\mathbb T}_{\alpha}}
\newcommand\Tb{{\mathbb T}_{\beta}}
\newcommand\Tc{{\mathbb T}_{\gamma}}

\newcommand\dbar{\overline\partial}

\newcommand\alphas{\mbox{\boldmath$\alpha$}}

\newcommand\xis{\mbox{\boldmath$\xi$}}
\newcommand\etas{\mbox{\boldmath$\eta$}}
\newcommand\betas{\mbox{\boldmath$\beta$}}
\newcommand\gammas{\mbox{\boldmath$\gamma$}}
\newcommand\deltas{\mbox{\boldmath$\delta$}}

\newcommand\PerDom{\mathcal P}

\newcommand\Field{\mathbb F}

\newcommand\Dual{\mathcal D}
\newcommand\Duality\Dual

\newcommand\iH{\mathbb H}
\newcommand\lk{\mathrm{lk}}

\newcommand\rhoL{\widetilde{\rho}}
\newcommand\nL{\widetilde{n}}
\newcommand\phiL{\widetilde{\phi}}
\newcommand\ModDegL{\widetilde{\ModDeg}}
\newcommand\UnparModDegL{\widehat{\widetilde{\ModDeg}}}

\newcommand\ModFlowL{\widetilde{\ModFlow}}
\newcommand\uL{\widetilde{u}}

\newcommand\Torsion{t}

\newcommand\RelDiff{\underline\epsilon}
\newcommand\Fill[1]{G_{#1}}
\newcommand\PerDoms{\mathcal P}
\newcommand\orL{\vec{L}}
\newcommand\orH{\vec{H}}
\newcommand\CFaa{\widehat{\widehat{\mathrm CF}}}

\newcommand\CFLa{\widehat{\mathrm CFL}}
\newcommand\CFLm{\mathrm{CFL}^-}
\newcommand\HFLm{\mathrm{HFL}^-}
\newcommand\HFLa{\widehat {\mathrm{HFL}}}
\newcommand\Kfrac{\mathbb K}
\newcommand\us{\mathbf u}
\newcommand\vs{\mathbf v}
\newcommand\ws{\mathbf w}
\newcommand\zs{\mathbf z}

\newcommand\spincrel\relspinc

\newcommand\CFK{\mathrm{CFK}}
\newcommand\HFK{\mathrm{HFK}}

\newcommand\CFKa{\widehat\CFK}

\newcommand\HFKa{\widehat\HFK}

\title[{Holomorphic disks and link invariants}] 
{Holomorphic disks, link invariants, and the multi-variable
Alexander polynomial}

\author[Peter Ozsv{\'a}th]{Peter Ozsv\'ath}
\address{Department of
Mathematics, Columbia University, 
New York, NY 10027 \newline
\indent{\emailfont{petero@math.columbia.edu}}}
\thanks{PSO was supported by NSF grant number DMS 0234311}

\author[Zolt{\'a}n Szab{\'o}]{Zolt{\'a}n Szab{\'o}} 
\address{Department of
Mathematics, Princeton University, New Jersey 08544 \newline
\indent{\emailfont{szabo@math.princeton.edu}}}
\thanks{ZSz was supported by NSF grant number DMS 0406155}

\begin{document}

\begin{abstract}  
  We define a Floer-homology invariant for links in $S^3$, and study
  its properties.
\end{abstract}

\maketitle

\section{Introduction}

The knot Floer homology of \cite{Knots} and~\cite{RasmussenThesis} is
an invariant for knots in $S^3$ whose Euler characteristic is the
Alexander polynomial of the knot.  Our aim here is to give a suitable
generalization of this invariant to links in $S^3$, giving rise to an
invariant whose Euler characteristic is the multi-variable Alexander
polynomial. 

Specifically, let $L\subset S^3$ be a link with $\ell$ components. Let
$H=H_1(S^3-L;\Z)$. Let $\Z[H]$ denote the group-ring of $H$, written
as sums $$\sum_{h\in H} a_h\cm e^h$$ where $a_h\in\Z$ is zero for all
but finitely many $h\in H$.  Note that $H\cong \Z^{\ell}$ is generated
by the meridians $\{\mu_i\}_{i=1}^\ell$ for the components $K_i$ of
$L$.  Thus, an orientation for $L$, denoted by $\orL$, induces an
identification between $\Z[H]$ and the ring of Laurent polynomials in
$\ell$ variables (corresponding to the components of $L$). Consider
now the affine lattice $\iH=\iH(L)$ over $H$, given by elements
$$\sum_{i=1}^{\ell} a_i\cm [\mu_i],$$
where $a_i\in \Q$ satisfies the
property that $2a_i + \lk(K_i,L-K_i)$ is an even integer, where here
$\lk$ denotes linking number.

We define a link invariant $\HFLa(\orL)$, which has the structure of a
vector space over the field $\Field=\Zmod{2}$ equipped with a
splitting into direct summands indexed by pairs consisting of an
integer (the ``homological grading'') and an element of 
$\iH$
$$\HFLa(\orL) \cong
\bigoplus_{d\in \Z, h \in \iH}\HFLa_d(\orL,h).$$

A few remarks are in order about these gradings.  First of all, the
ranks of the groups $\HFLa(\orL,h)$ are independent of the orientation
of $L$, but their homological gradings depend on this data.  Its
graded Euler characteristic is the Alexander polynomial in the
following sense.
 
Recall that an $\ell$-component link $L$ has a
symmetric multi-variable Alexander polynomial
$\Delta_L$. The link invariant is related to the multi-variable 
Alexander
polynomial by the relation
\begin{equation}
\label{eq:EulerHFLa}
\sum_{h\in \iH} \chi (\HFLa_*(\orL,h)) \cm e^h {=} 
\left\{\begin{array}{ll}
\left(\prod_{i=1}^\ell (T_i^{\OneHalf}-T_i^{-\OneHalf})\right)\cm \Delta_L  &{\text{if $\ell>1$}} \\
\Delta_L  &{\text{if $\ell=1$}} \\
\end{array}\right.
\end{equation}

The symmetry of the Alexander polynomial has the following manifestation
in link Floer homology: there is an identification
$$\HFLa_*(\orL,h)\cong \HFLa_{*-2\delta(h)}(\orL,-h),$$
where here
$$\delta(\sum_{i=1}^{\ell} a_i\cm [\mu_i])=\sum_{i=1}^\ell a_i.$$

This invariant should be compared with the link invariant described
in~\cite{Knots}. Specifically, that paper gives an invariant for
oriented links $\orL\subset S^3$
$$\HFKa(\orL)\cong\bigoplus_{d,s\in\Z}\HFKa_d(\orL,s),$$
using the
observation that an oriented $\ell$-component link in $S^3$ naturally
gives rise to a null-homologous knot in $\#^{\ell-1}(S^2\times S^1)$,
to which one can in turn apply the knot Floer homology functor,
obtaining a bigraded theory associated to this link. One of the
factors of the bigrading comes from the $\Z$-grading by $s\in\Z$ as
above, and the second comes from the internal homological grading of
the Heegaard Floer homology of $\#^{\ell-1}(S^2\times S^1)$ (whose
grading takes values in $\Z+\frac{\ell-1}{2}$).  Taking the Euler
characteristic of this knot Floer homology in a suitable sense
gives a normalized version of the Alexander-Conway polynomial of the
oriented link. An orientation $\orL$ of a link $L$ gives
rise to a homomorphism $o\colon H\longrightarrow \Z$, which extends
to a map of $\iH$ to the integers. 
Under the map $o$, the
multi-variable Alexander polynomial is carried to the Alexander-Conway
polynomial. This fact admits the following generalization on the
homological level:

\begin{theorem} 
  \label{thm:IdentifyWithLinkHomology}
  Let $L\subset S^3$ be a link endowed with an orientation, denoted
  $\orL$. Then, we have an identification
  $$\HFKa_{*+(\frac{\ell-1}{2})}(\orL,s)\cong \sum_{\{h\in \iH\big| o(h)=s\}} \HFLa_*(\orL,h),$$
  where  $o\colon \iH\longrightarrow \Z$ is the 
  natural homomorphism induced by the orientation.
\end{theorem}

Knot Floer homology can be viewed as the homology of the graded object
associated to a filtered chain complex whose total homology is $\Z$;
and indeed, the filtered chain homotopy type of this complex is a knot
invariant, cf.~\cite{4BallGenus}, \cite{RasmussenThesis}
(see also~\cite{EunSooLee}, \cite{RasmussenSlice}, and~\cite{Gornik}
for corresponding results in Khovanov's homology~\cite{Khovanov},
\cite{KhovanovRozansky}). Thus, there is a spectral sequence
starting with knot Floer homology, and converging to an $E^{\infty}$
term which has rank one.
We have the following
generalization of this fact to link Floer homology:

\begin{theorem}
  \label{thm:FilterHFLa}
  There is a spectral sequence whose $E_1$ term is $\HFLa(\orL)$,
  and whose  $E^\infty$ term is isomorphic to the exterior
  algebra $\Wedge^* V$, where $V$ is a vector space of rank $\ell-1$.
  In fact, the spectral sequence is an invariant of the
  link $L$.
\end{theorem}

The proof is given in Section~\ref{sec:IdentifyHeegaardFloers}, but we
pause here for a few remarks on the construction of this spectral sequence
and its meaning.

An orientation for $L$ gives a basis $(m_1,...,m_\ell)$ for $H$ (given
by the oriented meridians of the various components of $L$), and hence
a partial ordering on $H$, defined by $h_1\geq h_2$ if
$$h_1-h_2=\sum_{i} t_i \cm m_i,$$
where all the $t_i$ are non-negative integers. 
The above
theorem is proved by constructing a chain complex
$\CFLa(\orL)$
(cf. Definition~\ref{def:CFLa}) which admits both a $\Z$-grading and an 
$\iH$-filtration, i.e.
the group underlying this chain complex splits as a group
$$\CFLa(\orL)\cong \bigoplus_{d\in\Z,h\in \iH} \CFLa_d(\orL,h),$$
and its differential ${\widehat \partial}\colon \CFLa(\orL) \longrightarrow
\CFLa(\orL)$ carries 
$\CFLa_{d}(\orL,h)$ into
$$\bigoplus_{h'\leq h}\CFLa_{d-1}(\orL,h').$$
Moreover, there is an induced
differential 
$${\widehat{\widehat\partial}}\colon \CFLa_{d}(\orL,h) 
\longrightarrow \CFLa_{d-1}(\orL,h),$$
obtained by post-composing the differential ${\widehat\partial}$
with the projection map
$$\bigoplus_{h'\leq h}\CFLa(\orL,h')\longrightarrow \CFLa(\orL,h).$$
The homology of this graded object is identified with
$\HFLa_d(\orL,h)$. Moreover, the filtered
chain homotopy type of $\CFLa(\orL)$ 
is a link invariant, and its the total homology 
is identified with $\Wedge^* V$.

Succinctly, this equips $\HFLa(\orL)$ with an additional differential
mapping $\HFLa_d(\orL,h)$ into $\bigoplus_{h'<h}\HFLa_{d-1}(\orL,h')$,
in such a way that the total homology is $\Wedge^* V$. Moreover,
the absolute gradings on $\HFLa(\orL)$ are fixed so that 
the top-dimensional class in $\Wedge^* V$ is supported in degree zero.

In~\cite{AltKnots}, it was shown that the knot Floer homology of an
alternating knot is determined explicitly in terms of the Alexander
polynomial and signature of the knot. Using this result, together with
Equation~\eqref{eq:EulerHFLa} and
Theorem~\ref{thm:IdentifyWithLinkHomology}, we obtain the following:

\begin{theorem}
  \label{thm:AltLink}
  Let $L$ be an $\ell>1$-component
  oriented link with connected, alternating projection.
  Letting $\Delta_L(T_1,...,T_\ell)$ denote the multi-variable
  Alexander polynomial of $L$, write
  $$\left(\prod_{i=1}^{\ell}(T_i^{\OneHalf}-T_i^{-\OneHalf})\right)\cm \Delta_L(T_1,...,T_\ell)=\sum_{h\in\iH} a_h \cm e^h. $$
  Then, 
  $$\HFLa(\orL,h) = \Field^{|a_h|}_{(o(h)+\frac{\sigma-\ell+1}{2})},$$
  where here $\sigma$ denotes the signature of the oriented link $L$,
  and $\Field^n_{(d)}$ denotes the $n$-dimensional graded $\Field$-vector
  space supported entirely in grading $d$.
\end{theorem}

Link Floer homology satisfies a K{\"u}nneth principle for connected sums.
Specifically, 
let $L_1$ and $L_2$ be a pair of oriented links, and distinguish
components $K_1\in L_1$ and $K_2\in L_2$. Let $L_1\#L_2$ denote the
link obtained from the disjoint union of $L_1$ and $L_2$, via a
connected sum joining $K_1$ and $K_2$.
There is a natural map $$\iH(L_1)\oplus \iH(L_2)
\longrightarrow \iH(L_1\#L_2),$$ written $h_1, h_2\mapsto h_1\#
h_2$. This is the map which sends  both the meridian for $K_1$ 
and the meridian for $K_2$ (in $L_1$ and $L_2$) 
to the meridian for 
the connected sum $K_1\# K_2$.

\begin{theorem}
\label{thm:KunnethHFLa}
There is an isomorphism
$$\HFLa_*(\orL_1\#\orL_2,h)\cong \sum_{\{h_1\in \iH(L_1), h_2\in
\iH(L_2) \big| h_1\# h_2=h\}}\HFLa_*(\orL_1,h_1)\otimes \HFLa_*(\orL_2,h_2).$$
\end{theorem}

The above theorem is proved in Section~\ref{sec:Kunneth}, along with 
some of its natural generalizations.

\subsection{A further variant}

There are several variations of link Floer homology.  In
Section~\ref{sec:DefHFL}, we construct a chain complex
$\CFLm(\orL,h)$, which is a $\Z$-graded and $H$-filtered chain complex
of free modules over the ring $\Field[U_1,...,U_\ell]$.
Multiplication by $U_i$ lowers homological degree by two and it lowers
the filtration level by the basis element $m_i\in H$.  As shown in
Theorem~\ref{thm:InvarianceHFm}, the filtered chain homotopy type of
this complex is an invariant of the link.  The relationship between
this and the earlier construction is encoded in the fact that
$\CFLa(\orL)$ is gotten from $\CFLm(\orL)$ by setting each $U_i=0$.

The homology of the associated graded object is an oriented
link invariant $\HFLm(\orL)$
which is a module over the ring $\Field[U_1,...,U_\ell]$, endowed
with a $\Z$-grading (inherited from a $\Z$-grading on $\CFLm(\orL)$)
and an additional grading by elements of $\iH$ (induced from the
filtration). We denote this multi-grading
$$\HFLm(\orL) \cong \bigoplus_{d\in\Z, h\in \iH} \HFLm_d(\orL,h).$$ 
Calculating $\HFLm(\orL,h)$ is more challenging than calculating
$\HFLa(\orL,h)$, as its differential counts more holomorphic disks.  

The Euler characteristic in this case is given by the formula:

\begin{equation}
\label{eq:EulerHFLm}
\sum_{h\in \iH} \chi (\HFLm_*(\orL,h)) \cm e^h \dot{=} 
\left\{\begin{array}{ll}
\Delta_L  &{\text{if $\ell>1$}} \\
\frac{\Delta_L}{(1-T)}  &{\text{if $\ell=1$}}, \\
\end{array}\right.
\end{equation}
where here $f\dot{=}g$ means that two polynomials differ by multiplication
by units.
More succinctly, Equation~\eqref{eq:EulerHFLm} says that
the Euler characteristic of $\HFLm(\orL,h)$ is the Milnor
torsion of $L\subset S^3$~\cite{MilnorTorsion}.

The fact that $\CFLa(\orL)=\CFLm(\orL)/\{U_i=0\}_{i=1}^\ell$ as chain complexes
has the following manifestation on the level of homology, proved
in Section~\ref{sec:LinkInvariants}:

\begin{theorem}
  \label{thm:HFLaHFLm}
  For each fixed $d\in \Z$ and $h\in \iH$, the $\Field$-module
  $\HFLm_{d}(\orL,h)$ is the homology of a filtered chain complex whose
  $E_1$ term in dimension $d$ is  given by 
  $$\bigoplus_{(a_1,...,a_\ell)\geq 0} U_1^{a_1}\cm ...\cm
  U_\ell^{a_\ell} \cm \HFLa_{d+2a_1+...2a_\ell}(\orL,h+a_1\cm
  m_1+...+a_\ell\cm m_\ell).$$
\end{theorem}

For $\HFLm$, we have the following analogue of Theorem~\ref{thm:FilterHFLa}, 
proved in Section~\ref{sec:Relate}:

\begin{theorem}
  \label{thm:FilterHFLm}
  There is a spectral sequence, which is a link invariant,
  whose $E_1$ term is $\HFLm(\orL)$, and
  whose $E^{\infty}$ term is isomorphic to the
  $\Field[U_1,...,U_\ell]$-module $\Field[U]$, where each $U_i$ acts
  as multiplication by $U$. 
\end{theorem}

\subsection{About the construction}

Link Floer homology is constructed using suitable multiply-pointed
Heegaard diagrams for links. More precisely, if $L$ is a link with
$\ell$ components, we consider a Heegaard decomposition of $S^3$ as
$U_\alpha\cup U_\beta$, with the property that $L\cap U_\alpha$ and
$L\cap U_\beta$ consists of $\ell$ unknotted arcs. The link $L$ can
now be encoded in a genus $g$ Heegaard diagram for $S^3$, with
$g+\ell-1$ attaching circles $\alpha_1,...,\alpha_{g+\ell-1}$ for the
index one handles, and $g+\ell-1$ attaching circles
$\beta_1,...,\beta_{g+\ell-1}$ for the index two attaching circles,
and also $2\ell$ points $w_1,z_1,...,w_\ell,z_\ell$ where the link
crosses the mid-level. $\HFLa$ is a variant of Lagrangian Floer
homology in the $g+\ell-1$-fold symmetric product of $\Sigma$
punctured in the basepoints $w_i$ and $z_i$.
These topological considerations lead naturally to the notion of {\em balanced}
Heegaard diagrams, which are Heegaard diagrams for a (closed, oriented)
three-manifold with $\ell$ zero-
and $\ell$ three-handles.  Theorems~\ref{thm:FilterHFLa} and
~\ref{thm:FilterHFLm} are obtained from extending Heegaard Floer
homology to the case of such balanced Heegaard diagrams.

This paper is organized as follows.  In Section~\ref{sec:Algebra}, we
review some of the algebraic terms used throughout this paper. In
Section~\ref{sec:HeegaardDiagrams}, we discuss balanced Heegaard
diagrams associated to links, and also the topological data which can
be extracted from them. We also address admissibility issues which
will be required to define the Heegaard Floer complexes.  In
Section~\ref{sec:DefHFL}, we describe the Heegaard Floer homology
complexes associated to balanced Heegaard diagrams.  In this case, the
proof that $\partial^2=0$ is slightly more subtle than the usual case
considered in~\cite{HolDisk}: specifically, it is now no longer true
that the total count of {\em boundary degenerations}, disks with
boundary lying entirely in $\Ta$ (or $\Tb$), is zero. These issues are
addressed in Section~\ref{sec:Analysis}, where the analytical
preliminaries are set up. For certain technical reasons, we find it
also convenient to adopt the ``cylindrical'' approach to Heegaard
Floer homology developed by Lipshitz~\cite{LipshitzCyl}, where one
considers pseudo-holomorphic multi-section of the trivial
$\Sigma$-bundle over a disk, rather than disks in the symmetric
product of $\Sigma$.

With the help of this, in Section~\ref{sec:IdentifyHeegaardFloers},
it is established that $\partial^2=0$, and indeed, the
Heegaard Floer homology for balanced Heegaard diagrams is identified
with the usual Heegaard Floer homology.

With this background in place, the invariants for links are easy to construct,
and their invariance properties are readily verified in
Section~\ref{sec:LinkInvariants}. In particular,
Theorems~\ref{thm:FilterHFLa}, \ref{thm:HFLaHFLm},
and~\ref{thm:FilterHFLm} are quick consequences of the constructions.

With the link invariants in hand, we turn to some of their basic properties.
In Section~\ref{sec:Symmetry} we establish certain symmetry properties,
which parallel the usual symmetry of the Alexander polynomial.

In Section~\ref{sec:Euler}, we turn to the Euler characteristic
statements, verifying Equations~\eqref{eq:EulerHFLa}
and~\eqref{eq:EulerHFLm}.

In Section~\ref{sec:Relate} we relate the present form of link
homology with the earlier form derived from knot Floer homology
in~\cite{Knots}, establishing
Theorem~\ref{thm:IdentifyWithLinkHomology}.

The K{\"u}nneth principle for connected sums
(Theorem~\ref{thm:KunnethHFLa} above, and also some more general
statements for $\CFLm$) is established in Section~\ref{sec:Kunneth}.

In Section~\ref{sec:Examples}, we establish Theorem~\ref{thm:AltLink}.
We give also some principles which help computing the spectral
sequence from Theorem~\ref{thm:FilterHFLa}. These principles allow
one to determine the spectral sequence for all two-bridge
links from the signature and the multi-variable Alexander polynomial.

We illustrate these principles in some particular examples, giving
also some calculations for the two non-alternating, seven-crossing
links, as well.

\vskip.2cm

\noindent{\bf{Further remarks.}}
We have set up here link Floer homology $\HFLa$ as the homology of a
graded object associated to a filtration of a chain complex for
$\#^{\ell-1}(S^2\times S^1)$. As such, it gets the extra differentials
promised in Theorem~\ref{thm:FilterHFLa}.  If one is not interested in
this extra structure, but only $\HFLa$ as a graded group, then its
construction is somewhat more elementary than the constructions
described here. Properties of this invariant, and further
computations, are given in~\cite{LinkTwo}.

To some degree, link Floer homology can be viewed as a
categorification of the multi-variable Alexander polynomial. It is
interesting to compare this with the recent categorification of the
HOMFLY polynomial~\cite{KhovanovRozanskyII}, see also~\cite{Khovanov},
\cite{KhovanovRozansky}, \cite{DunfieldGukovRasmussen}.

\vskip.2cm

\noindent{\bf{Acknowledgements.}} We would like to thank
Robert Lipshitz, Jacob Rasmussen, and Andr{\'a}s Stipsicz for useful
conversations during the course of this work. We would also like to
thank Jiajun Wang for his input on an early version of this manuscript.

\section{Algebraic preliminaries}
\label{sec:Algebra}

We begin by fixing some terminology from homological algebra which
will be used throughout this paper.  We give the set $\Z^{\ell}$ the
partial ordering $a=(a_1,...,a_\ell)\leq b=(b_1,...,b_\ell)$ if each
$a_i\leq b_i$.  Let ${\mathfrak S}$ be an affine space for $\Z^\ell$.
Natural examples include the affine space $\iH$ introduced in the
introduction, where we think of $\Z^\ell$ as identified with the first
homology of a link complement (where the identification induced by
orientations on the link). Another is the set of relative $\SpinC$
structures over the link complement, c.f.
Section~\ref{subsec:RelSpinC} below.  The ordering on $\Z^\ell$
induces an ordering on ${\mathfrak S}$.  An ${\mathfrak S}$-filtered
module over $\Field$ is an $R$-module $M$ equipped with an exhausting
family of sub-modules $\Filt(M,a)\subseteq M$ indexed by $a\in
{\mathfrak S}$, with the containment relation $\Filt(M,a)\subseteq
\Filt(M,b)$ if $a\leq b$.  A module homomorphism $\phi\colon
M\longrightarrow M'$ with the property that for all $a\in\Z^\ell$,
$\phi(\Filt(M,a))\subseteq \Filt(M',a)$ is called a {\em morphism of
  ${\mathfrak S}$-filtered modules}.  A {\em filtered ${\mathfrak
    S}$-complex} is a chain complex for which the differential
$\partial$ is a morphism of ${\mathfrak S}$-filtered modules. The
homology of an ${\mathfrak S}$-filtered chain complex inherits a
natural $\Z^\ell$ filtration.

Two ${\mathfrak S}$-filtered chain maps $\phi_1,\phi_2\colon A
\longrightarrow B$ between ${\mathfrak S}$-filtered chain complexes
$A$ and $B$ are said to be {\em filtered chain homotopic} if there is
a morphism of ${\mathfrak S}$-filtered modules $H\colon A
\longrightarrow B$ with $\partial\circ H-H\circ\partial =
\phi_1-\phi_2$. Two ${\mathfrak S}$-filtered chain complexes $A$ and
$B$ are {\em ${\mathfrak S}$-filtered chain homotopy equivalent} if
there are filtered chain maps $f\colon A \longrightarrow B$ and
$g\colon B\longrightarrow A$ with the property that both $f\circ g$
and $g\circ f$ are filtered chain homotopic to the the corresponding
identity maps. If $C$ and $C'$ are filtered chain homotopy
equivalent, we write $C\simeq C'$.  Clearly, this forms an equivalence
relation on the set of ${\mathfrak S}$ filtered chain complexes, and
the induced equivalence class of a given ${\mathfrak S}$-filtered
chain complex is called its {\em ${\mathfrak S}$-filtered chain
  homotopy type}.

Given a ${\mathfrak S}$-filtered complex, we can form the associated graded
object
$$\gr(C)=\bigoplus_{a\in {\mathfrak S}}\gr(C,a),$$
where 
$$\gr(C,a)=
\Coker({\bigoplus_{b< a}\Filt(C,b)}\longrightarrow\Filt(C,a)).$$
Clearly the homology of the associated graded object 
of a $\Z^\ell$-filtered chain complex $C$ depends on only
the filtered chain homotopy type of $C$.

Given any set ${\mathfrak T}\subset {\mathfrak S}$ with the property that for all $a\in {\mathfrak T}$
if $b\leq a$, then $b\in {\mathfrak T}$, we can form the subcomplex
$C({\mathfrak T})\subset C$.

We will consider modules over the ring
$R=\Field[U_1,...,U_\ell]$. An ${\mathfrak S}$ filtered $R$-module
is an $R$-module whose underlying $\Field$-module (gotten by forgetting
the action of $U_i$) is ${\mathfrak S}$-filtered, and has the additional property
that
$$U^{a_1}\cm ... \cm U^{a_\ell}\cm \Filt(M,b)
\subseteq \Filt(M,b-a),$$
where $a=(a_1,...,a_\ell)$. The notions of morphisms, chain complexes,
homotopies, and homotopy type extend in a straightforward manner:
we consider maps which are simultaneously ${\mathfrak S}$ filtered and
which are also $R$-modules. If $C$ is a chain 
${\mathfrak S}$-filtered
chain complex of $R$-modules, 
the chain complex ${\widehat C}=C\otimes_{\Field[U_1,...,U_{\ell}]} \Field$
gotten by setting each $U_i=0$ is also a $\Z^{\ell}$-filtered chain complex
(whose filtered chain homotopy type depends on $C$ only up 
to its filtered chain homotopy type).

A {\em free} $\Z^n$-filtered chain complex of $R$-modules is one which admits a
homogeneous generating which freely generates the underlying complex
over $\Field[U_1,...,U_\ell]$. In particular, as a $\Field$-module,
$C$ splits as a direct sum
$$C=\bigoplus_{a\in{\mathfrak S}} C\{a\}.$$

\subsection{Operations on ${\mathfrak S}$-filtered chain complexes.}
\label{subsec:HomologicalProjection}

If $C$ is a $\Z^n$-filtered
chain complex, and $a\in\Z^n$,
then we can form the $\Z^n$ filtered complex
$C[a]$ whose underlying chain complex agrees with $C$, but
whose $\Z^n$-filtration is shifted by $a$; i.e.
the filtration $\Filt(C[a],b)=\Filt(C[a+b])$.

If $C$ is free $\Z^n$-filtered complex, and $i\in \{1,...,n\}$, we can
split the differential into components $D^b_a \colon C\{a\}
\longrightarrow C\{b\}$ with $b\leq a$; we can form a $\Z^{n-1}$
filtered chain complex $C^{(1)}$ by ``taking homology in the first
component''.  Specifically, write $D=D^1 + D'$, where $D^1$ consists
of all the components $D^b_a$ where $a$ and $b$ agree on all but the
first place. It is easy to see that $D^1$ is a differential, and hence
we can form the $\Z^{n-1}$-filtered chain complex
$$\Filt(C^{(1)},b)=
H_*\left(\bigcup_{\{a\big| (a_2,...,a_\ell)\leq b\}} \Filt(C,a),D^1\right),$$
endowed with the differential induced from $D'$.
Note that there remains an extra action of $U_1$ on this chain complex (which
does not change the filtration level).

Of course, this notion admits a straightforward adaptation
to taking the homology in the $i^{th}$ components for any $i\in\{1,...,\ell\}$.

\section{Heegaard diagrams}
\label{sec:HeegaardDiagrams}

We discuss here basic topological aspects of multiply-pointed Heegaard
diagrams, which are relevant for the study of links in
three-manifolds. The material here is mostly a straightforward
generalization of the singly-pointed case, which was studied
in~\cite{HolDisk}, and the doubly-pointed case from~\cite{Knots}.

\subsection{Heegaard diagrams for three-manifolds}

\begin{defn}
A {\em balanced $\ell$-pointed Heegaard diagram}
is a quadruple of data
$$(\Sigma,\alphas=\{\alpha_1,...\alpha_{g+\ell-1}\},
\betas=\{\beta_1,...,\beta_{g+\ell-1}\}, \ws=\{w_1,...,w_{\ell}\}),$$
where $\ell$ is some positive integer, $\Sigma$ is an oriented surface
of genus $g$, $\alphas=\{\alpha_1,...,\alpha_{g+\ell-1}\}$ is a
$g+\ell-1$-tuple of disjoint, simple closed curves which span a
$g$-dimensional sublattice of $H_1(\Sigma;\Z)$ (and hence they specify
a handlebody $U_\alpha$ which is bounded by $\Sigma$),
$\betas=\{\beta_1,...,\beta_{g+\ell-1}\}$ is a $g+\ell-1$-tuple of
disjoint, simple closed curves which span another $g$-dimensional
sublattice of $H_1(\Sigma;\Z)$ (specifying another handlebody
$U_\beta$), and $\ws$ is a collection of points in $\Sigma$ chosen as
follows.  Let $\{A_i\}_{i=1}^\ell$ denote the connected components of
$\Sigma-\alpha_1-...-\alpha_{g+\ell-1}$; let $\{B_i\}_{i=1}^\ell$
denote the connected components of
$\Sigma-\beta_1-...-\beta_{g+\ell-1}$. The points $w_i\in \Sigma$ are
constrained so that $w_i\in A_i\cap B_i$.
\end{defn}

A balanced $\ell$-pointed Heegaard diagram specifies a closed, oriented three-manifold
$Y=U_\alpha\cup_{\Sigma} U_{\beta}$,  endowing it with a cellular
decomposition whose zero-cells correspond to the components
$\{A_i\}_{i=1}^\ell$, its one-cells correspond to the circles
$\alpha_1,...,\alpha_{g+\ell-1}$, its two-cells correspond to
$\beta_1,...,\beta_{g+\ell-1}$, and its three-cells correspond to the
$\{B_i\}_{i=1}^\ell$.

\begin{defn}
  A balanced $\ell$-pointed Heegaard diagram is called {\em generic}
  if the circles $\alpha_i$ and $\beta_j$ meet transversally, for all
  $i,j\in \{1,...,g+\ell-1\}$. 
\end{defn}

Fix a connected, oriented three-manifold $Y$, and a generic
self-indexing Morse function on $Y$ which has the same number $\ell$
of index zero and three critical points. Fix also generic metric $g$,
together with a choice of $\ell$ gradient flowlines connecting each of
the index zero and three critical points. Then, there is an associated
generic balanced $\ell$-pointed Heegaard diagram for $Y$ whose
surface $\Sigma$ is the mid-level of the Morse function; $\alpha_i$ is
the locus of points on $\Sigma$ where the gradient flowlines leaving
the $i^{th}$ index one critical point meets $\Sigma$; similarly,
$\beta_i$ is the locus of points on $\Sigma$ which flow into the
$i^{th}$ index two critical point. Finally, for $i=1,...,\ell$, $w_i$
is the point on $\Sigma$ which lies on the distinguished gradient
flow-line connecting the $i^{th}$ index zero and index three critical
point.  If $(\Sigma,\alphas,\betas,\ws)$ is obtained in this manner
from a Morse function $f$, we call $f$ a Morse function {\em
  compatible} with the balanced Heegaard diagram for $Y$.

Given a generic $\ell$-pointed balanced Heegaard diagram for $Y$, it is easy
to construct a compatible Morse function $f$.

\begin{prop}
\label{prop:HeegaardMoves}
Any two generic balanced $\ell$-pointed Heegaard diagrams for $Y$ can be
connected by a sequence of the following moves:
\begin{list}
       {(\arabic{bean})}{\usecounter{bean}\setlength{\rightmargin}{\leftmargin}}
\item
\label{Moves:IsHSA} isotopies and handleslides of the $\alphas$ supported in the
  complement of $\ws$
\item 
 \label{Moves:IsHSB}isotopies and handleslides of the $\betas$ supported in the
  complement of $\ws$
\item 
\label{Moves:StabOneTwo}
index one/two stabilizations (and their inverses): forming the
  connected sum of $(\Sigma,\alphas,\betas,\ws)$ with a torus equipped
  with a new pair of curves $\alpha_g$ and $\beta_g$ which meet
  transversally in a single point
\item
  \label{Moves:ZeroThree}
index zero/three stabilizations (and their inverses):
  introducing a new pair of homotopic curves 
  $\alpha_{g+\ell}$ (disjoint from the $\alpha_i$ for $1\leq i\leq
  g+\ell-1$) and $\beta_{g+\ell}$ (disjoint from the other $\beta_i$)
  and a new basepoint $w_{\ell+1}$ in such  a manner that 
  $\alpha_{g+\ell}$ and $\beta_{g+\ell}$ are 
  homotopic in $\Sigma-w_1-...-w_{\ell}-w_{\ell+1}$, and
  each component of $\Sigma-\alpha_1-...-\alpha_{g+\ell}$ and
  $\Sigma-\beta_1-...-\beta_{g+\ell}$ contains some $w_i$.
\end{list}
\end{prop}

\begin{proof}
  It follows from standard Morse theory that any two (unpointed)
  Heegaard diagrams can be connected by moves of
  Types~\eqref{Moves:IsHSA}-\eqref{Moves:StabOneTwo}. The fact that
  this can be done in the complement of a single basepoint (the case
  $\ell=1$) can be established by trading an isotopy across the
  basepoint for a sequence of handleslides in the opposite direction,
  cf.~\cite[Proposition~\ref{HolDisk:prop:PointedHeegaardMoves}]{HolDisk}.
  
  The proof in general is established by showing that we can use the
  above moves to reduce the number of basepoints, and hence reducing
  to the case of a singly-pointed Heegaard diagram. This is done as
  follows.  Let $A$ be the two-chain with boundary a combination of
  $\alpha_1,...,\alpha_{g+\ell}$ which has $n_{w_{\ell+1}}(A)=1$ and
  $n_{w_i}(A)=0$ for all $i\leq \ell$. After a sequence of
  handleslides among the $\alphas$, we can arrange for $A$ to have
  only one boundary component, which we label $\alpha_{g+\ell}$.
  Indeed, another sequence of handleslides can be done to arrange
  furthermore for the genus of $A$ to be zero.  Let $B$ be
  the corresponding two-chain with boundary amongst the
  $\beta_1,...,\beta_{g+\ell}$, with $n_{w_{\ell+1}}(B)=1$
  and $n_{w_i}(B)=0$ for all $i\leq \ell$. Performing a
  sequence of handleslides among the $\betas$, we can reduce to the
  case where $B$ is a disk bounded by $\beta_{g+\ell}$. We
  can now form an index zero/three de-stabilization to delete
  $\alpha_{g+\ell}$, $\beta_{g+\ell}$ and $w_{g+\ell}$. The proof then
  follows by induction.
\end{proof}

\subsection{Relative $\SpinC$ structures}
\label{subsec:RelSpinC}

We pause our discussion on Heegaard diagrams to recall Turaev's
interpretation of $\SpinC$ structures on three-manifolds, see~\cite{Turaev},
compare also~\cite{KMcontact}.

Let $Y$ be a closed, oriented three-manifold. We say that two nowhere
vanishing vector fields $v$ and $v'$ are {\em homologous} if there is
a ball $B\subset Y$ with the property that $v$ and $v'$ are homotopic
(through nowhere vanishing vector fields) on the complement of $V$.
The set of equivalence classes of such vector fields can be 
naturally identified with the space $\SpinC(Y)$  of $\SpinC$ structures
over $Y$. In particular it is an affine space for $H^2(Y;\Z)$.

This notion has a straightforward generalization to the case of
three-manifolds with toroidal bounday. Specifically, let $(M,\partial
M)$ be a three-manifold with boundary consisting of a disjoint union
of tori $T_1\cup...\cup T_\ell$. The tangent bundle to the two-torus
has a canonical nowhere vanishing vector field, which is unique up to
homotopy (through nowhere vanishing vector fields).  Consider now
nowhere vector fields $v$ on $Y$ whose restriction to $\partial M$ are
identified with the canonical nowhere vanishing vector field on the
boundary tori (in particular, the vector field has no normal component
at the boundary). Two such vector fields $v$ and $v'$ are declared
homologous if there is a ball $B\subset M-\partial M$ with the
property that the restrictions of $v$ and $v'$ to $M-B$ are homotopic.
The set of homology classes of such vector fields is called the set of
{\em relative $\SpinC$ structures}, and it is an affine space for
$H^2(M,\partial M;\Z)$.  We denote this set by $\RelSpinC(M,\partial
M)$.

Multiplying vector fields by  $-1$ induces an involution on the
space of relative $\SpinC$ structures
$$J\colon \RelSpinC(M,\partial M) \longrightarrow \RelSpinC(M,\partial M).$$

If $\vec{v}$ is a vector field on a three-manifold $M$ with toroidal boundaries,
whose restriction to each bounding torus gives the canonical trivialization of the 
torus' tangent bundle, then we can consider the oriented two-plane field 
$\vec{v}^\perp$ of
vectors orthogonal to $\vec{v}$. Along $\partial M$, this two-plane field has a canonical 
trivialization, by outward pointing vectors. Hence, there is a well-defined
notion of a relative Chern class of this line field relative to its trivialization,
thought of as an element of $H^2(M,\partial M;\Z)$.
This descends to a well-defined assignment
$$c_1\colon \RelSpinC(M,\partial M)\longrightarrow H^2(M,\partial M;\Z).$$
Clearly,
$$c_1(J\cm \relspinc)=-c_1(\relspinc).$$

The reader familiar with~\cite{Knots} should be warned that in that
case, we were considering null-homologous knots, rather than links,
and hence there is a well-defined notion of a zero-surgery.
In~\cite{Knots}, relative $\SpinC$ structures were thought of as
absolute $\SpinC$ structures on this zero-surgery. This is a slightly
different point of view than the one taken here (where we no longer
have the luxury of referring to a zero-surgery).

\subsection{Intersection points and $\SpinC$ structures}
\label{subsec:IntPtsSpinCStruct}

Fix a generic $\ell$-pointed balanced Heegaard diagram for $Y$, and
let $f$ be a compatible Morse function.

Consider the $g+\ell-1$-fold symmetric product of $\Sigma$,
$\Sym^{g+\ell-1}(\Sigma)$, and let
\begin{eqnarray*}
\Ta=\alpha_1\times...\times \alpha_{g+\ell-1}
&{\text{and}}&
\Tb=\beta_1\times...\times \beta_{g+\ell-1}.
\end{eqnarray*}
Clearly, an intersection point $\x\in\Ta\cap\Tb$ corresponds to a
$g+\ell-1$-tuple of gradient flow-lines which connect all the index one
and two critical points. 

Let $\gamma_\x$ be the union of gradient flowlines passing through
each $x_i\in \x$, and $\gamma_\w$ be the union of gradient flowlines
passing through each $w_i\in\w$. The closure of
$\gamma_\x\cup\gamma_\w$ is a collection of arcs whose boundaries
consist of all the critical points of $f$. Moreover, each component
contains a pair of critical points whose indices have opposite
parities. Thus, we can modify the gradient vector field in an
arbitrarily small neighborhood of $\gamma_\x\cup\gamma_\w$ to obtain a
new vector field which vanishes nowhere in $Y$. Taking the homology
class of this vector field in the sense of Turaev~\cite{Turaev}, we
obtain a map from intersection points to $\SpinC$ structures over $Y$
\[\spinc_\w\colon \Ta\cap\Tb\longrightarrow \SpinC(Y).
\]
It is easy to see that this is well-defined, i.e. independent of the choice
of compatible Morse function $f$ and modification of the vector field
(compare Subsection~\ref{HolDisk:subsec:SpinCStructures} of~\cite{HolDisk}).

\subsection{Whitney disks and admissibility}
\label{subsec:Admissibility}

Fix a generic $\ell$-pointed Heegaard diagram for $Y$,
$(\Sigma,\alphas,\betas,\ws)$. Given $\x,\y\in\Ta\cap\Tb$, we can
consider Whitney disks from $\x$ to $\y$ relative to $\Ta$ and $\Tb$.
Each such Whitney disk gives rise to a two-chain on $\Sigma$; specifically,
$$\Sigma-\alpha_1-...-\alpha_{g+\ell-1}-\beta_1-...-\beta_{g+\ell-1}$$
consists of a collection of regions $\{\Omega_i\}_{i=1}^m$. Fix a 
reference point $p_i\in \Omega_i$, and let $\cald(\phi)$ denote the 
two-chain
$$\sum_{i=1}^m n_{p_i}(\phi)[\Omega_i],$$
where here $n_{p}(\phi)$ denotes the
algebraic intersection number of $\phi$ with the subvariety
$\{p\}\times\Sym^{g+\ell-2}(\Sigma)$. The element $\cald(\phi)$ specifies
the relative homology class induced from the Whitney disk $\phi$.

Let $\pi_2(\x,\y)$ denote the space of homology classes of Whitney
disks.
Let $n_\w\colon \pi_2(\x,\y)\longrightarrow \Z^\ell$ be the map
which sends $\phi$ to $(n_{w_1}(\phi),...,n_{w_\ell}(\phi))$

Let $\pi_2(\alpha)$ denote the space of homology classes of disks with
boundary in $\Ta$. Clearly, 
we have isomorphisms
\begin{eqnarray}
\label{eq:nwIso}
n_\w \colon \pi_2(\alpha)\stackrel{\cong}\longrightarrow \Z^\ell
&{\text{and}}&
n_\w \colon \pi_2(\beta)\stackrel{\cong}\longrightarrow \Z^\ell.
\end{eqnarray}
There is an exact sequence
\begin{equation}
\label{eq:DecomposePD}
\begin{CD}
0@>>>\Z@>>> \pi_2(\alpha)\oplus \pi_2(\beta)
@>>> \pi_2(\x,\y)@>>>H^1(Y;\Z)@>>> 0;
\end{CD}
\end{equation}

Moreover, we have an exact sequence
\begin{equation}
\label{eq:GroupPerDoms}
\begin{CD}
  0@>>> \PerDoms@>>>\pi_2(\x,\y)@>{n_\w}>> \Z^\ell@>>> 0.$$
\end{CD}
\end{equation}

\begin{defn}
  The group $\PerDoms$ is called the group of {\em periodic domains}.
\end{defn}

In the case where $H^1(Y;\Z)=0$, we have that
$\PerDoms\cong \Z^{\ell-1}$.

As in Subsection~\ref{subsec:Admissibility} of~\cite{HolDisk}, we need
further restrictions on the Heegaard diagram to obtain a reasonable
chain complex (whose homology is the Heegaard Floer homology).
defined.

\begin{defn}
  A generic, balanced $\ell$-pointed Heegaard diagram is called {\em
    weakly admissible} if for any non-trivial homology class
  $\phi\in\pi_2(\x,\x)$ with $n_\w(\phi)=0$ (meaning that
  $n_{w_i}(\phi)=0$ for each $w_i\in\w$),  the domain $\cald(\phi)$
  has both positive and negative local multiplicities.
\end{defn}

\begin{prop}
  \label{prop:Admissibility}
  Let $Y$ be a three-manifold with $H^1(Y;\Z)=0$.  Any balanced
  $\ell$-pointed Heegaard diagram for $Y$ is isotopic to a weakly
  admissible balanced $\ell$-pointed Heegaard diagram.  
\end{prop}

\begin{proof}
  We embed in $\Sigma$ a tree $\Gamma$ whose $\ell$ vertices are the points
  $w_i$. We claim that admissibility can be achieved by isotoping some
  of the $\beta$-curves in a regular neighborhood of $\Gamma$.
  Specifically, if $\gamma$ is an arc connecting $w_1$ to $w_2$ in
  $\Gamma$, then perform an isotopy of the $\beta$ circles in a
  regular neighborhood of $\gamma$ in such a manner that there is a
  pair of arcs $\delta_1$ and $\delta_2$ so that
  $\delta_1\cup\delta_2$ is isotopic to $\gamma$ as an arc from $w_1$
  to $w_2$, but $\delta_1$ is disjoint from the $\beta$ circles while
  $\delta_2$ is disjoint from the $\alpha$ circles. Moreover, we find
  another pair of arcs $\delta_1'$ and $\delta_2'$ so that
  $\delta_1'\cup\delta_2'$ is isotopic to $\gamma$, only now
  $\delta_1'$ is disjoint from the $\alpha$ circles while $\delta_2'$
  is disjoint from the $\beta$ circles
  (see Figure~\ref{fig:Winding} for an illustration).  Isotoping the $\beta$
  circles in a regular neighborhood of all the edges in $\Gamma$ as
  above, we obtain a Heegaard diagram that we claim is weakly admissible.

\begin{figure}
\mbox{\vbox{\epsfbox{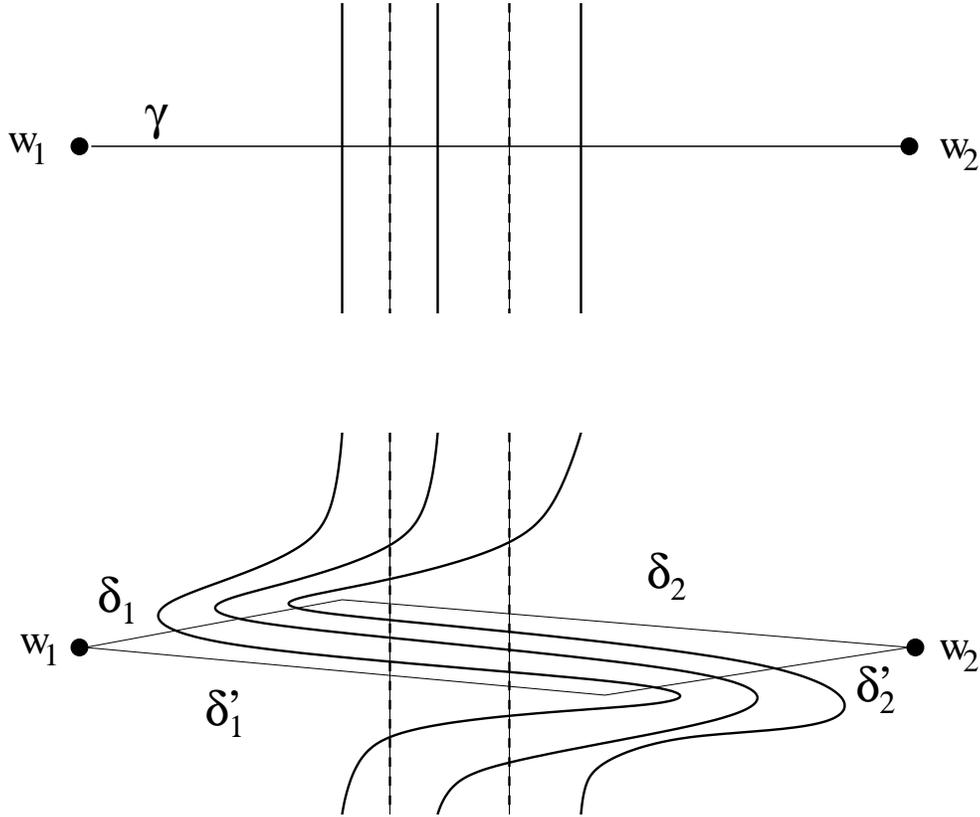}}}
\caption{\label{fig:Winding}
  {\bf{Winding to achieve admissibility}} This is an illustration of
  the procedure described in Proposition~\ref{prop:Admissibility}.
  Here, the dotted lines denote the possible $\alpha$ curves, and the
  solid lines represent $\beta$ curves.}
\end{figure}

  According to Equation~\eqref{eq:DecomposePD}, any $P\in\PerDom$ can
  be decomposed as $P=A+B$, with 
  $A\in\pi_2(\alpha)$ and $B\in\pi_2(\beta)$. The condition
  that $P$ is a periodic domain ensures that $n_\w(A)+n_\w(B)=0$.
  According to Equation~\eqref{eq:nwIso}, $P$ is uniquely determined
  by $n_\w(A)$, modulo addition of $\Sigma$.
  
  Suppose now that $P$ has the property that
  the oriented intersection number of $\partial A$ with $\gamma$
  is non-zero.  Then, at the intermediate endpoint of $\delta_1$, we
  see that $A+B$ has local multiplicity given by $\partial A \cap
  \gamma$, while at the intermediate endpoint of $\delta_1'$, $A+B$
  has local multiplicity given by $\partial B \cap \gamma=-\partial
  A\cap\gamma$. Thus, if for some edge $\gamma$ in $\Gamma$, $\partial
  A \cap \gamma\neq 0$, then $P=A+B$ has both positive and negative
  coefficients.  However, if $P=A+B$, where $A\in\pi_2(\alpha)$ and
  $B\in\pi_2(\beta)$, and $\partial A$ has algebraic intersection
  number equal to zero with each edge in $\Gamma$, then, after
  subtracting off some number ($n_{w_1}(A)$) of copies of $\Sigma$, we
  can write $P=A'+B'$, where $A'\in\pi_2(\alpha)$, $B'\in\pi_2(\beta)$
  and $n_{w_i}(A)=n_{w_i}(B)=0$. According to
  Equation~\eqref{eq:nwIso}, then $A'=0$, and hence $P=0$.
\end{proof}

\subsection{Heegaard diagrams and links}

Fix an $\ell$-pointed Heegaard diagram $(\Sigma,\alphas,\betas,\ws)$
for a three-manifold $Y$, and choose also an additional $\ell$-tuple
of basepoints $\zs=\{z_1,...,z_\ell\}$, with the property that for each
$i=1,...,\ell$, both $w_i$ and $z_i$ are contained in the same
component
$$A_i\subset \Sigma-\alpha_1-...-\alpha_{g+\ell-1}$$
and 
$$B_i\subset \Sigma-\beta_1-...-\beta_{g+\ell-1}.$$
This data gives rise to an oriented, 
$\ell$-component link $L$ in $Y=U_{\alpha}\cup_{\Sigma} U_\beta$.

\begin{defn}
  The diagram $(\Sigma,\alphas,\betas,\ws,\zs)$ as above is said to be
  a {\em $2\ell$-pointed Heegaard diagram} for the oriented link
  $\orL$ in $Y$. 
\end{defn}

Conversely, given an oriented, $\ell$-component link, one can find a
self-indexing Morse function $f\colon Y\longrightarrow \R$ with $\ell$
index zero and three critical points, and $g+\ell-1$ index one and two
critical points, with the additional property that there are two
$\ell$-tuples of
flowlines $\gamma_\w$ and $\gamma_\z$
connecting all the index three and index zero critical points,
so that our oriented link $L$ can be realized as the difference
$\gamma_\z-\gamma_\w$.
Such a Morse function gives rise to a $2\ell$-pointed Heegaard diagram
for $\orL$ in $Y$.

\begin{defn}
  A $2\ell$-pointed Heegaard diagram for a link $\orL\subset Y$ is
  {\em weakly admissible} if the underlying $\ell$-pointed Heegaard
  diagram for $Y$ (gotten by disregarding the $\zs$) is weakly
  admissible.
\end{defn}

\begin{prop}
  If $\orL$ is an oriented $\ell$-component link in $Y$, then
  there is a corresponding (weakly admissible) $2\ell$-pointed Heegaard
  diagram.  Any two (weakly admissible) $2\ell$-pointed Heegaard
  diagrams for the same oriented link $\orL\subset Y$ can be connected by
  a sequence of moves of the following types:
  \begin{itemize}
  \item isotopies and handleslides of the $\alphas$ supported in the
    complement of $\ws$ and $\zs$ 
  \item isotopies and handleslides of the $\betas$ supported in the
    complement of $\ws$ and $\zs$
  \item index one/two stabilizations (and their inverses): forming the
    connected sum of $(\Sigma,\alphas,\betas,\ws)$ with a torus equipped
    with a new pair of curves $\alpha_g$ and $\beta_g$ which meet
    transversally in a single point.
  \end{itemize}
  Moreover, if we start and end with weakly admissible Heegaard
  diagrams, then we can assume that all the intermediate Heegaard
  diagrams are also weakly admissible.
\end{prop}

\begin{proof}
  Without the admissibility hypothesis, the above result follows from
  Morse theory in the usual manner. Admissibility can be achieved as in
  the proof of Proposition~\ref{prop:Admissibility}.
\end{proof}

\subsection{Intersection points and link diagrams}
\label{subsec:RelSpinCLinks}

Given a $2\ell$-pointed  Heegaard diagram for a link, consider the tori
\begin{eqnarray*}
\Ta=\alpha_1\times...\times \alpha_{g+\ell-1}
&{\text{and}}&
\Tb=\beta_1\times...\times \beta_{g+\ell-1}.
\end{eqnarray*}
Given a pair $\x,\y\in\Ta\cap\Tb$ of intersection points, 
we can find paths
\begin{eqnarray*}
a\colon[0,1]\longrightarrow \Ta
&{\text{and}}&
b\colon[0,1]\longrightarrow \Tb
\end{eqnarray*}
with $\partial a = \partial b = \x-\y$. Viewing these paths as 
one-chains in $\Sigma$, supported away from the reference points
$\{w_i\}_{i=1}^\ell$ and $\{z_i\}_{i=1}^\ell$, we obtain a one-cycle
$\epsilon(\x,\y)$ in the complement $Y-L$. This
assignment clearly descends to give a well-defined map
$$\RelDiff_{\w,\z}\colon (\Ta\cap\Tb)\times(\Ta\cap\Tb) \longrightarrow
H_1(Y-L;\Z)$$

\begin{lemma}
  \label{lemma:Filtration1}
  An oriented link in $Y$ gives rise to a map
  $$\Pi\colon H_1(Y-L)\longrightarrow \Z^\ell$$
  (where $\Pi_i(\gamma)$ is the linking number of $\gamma$
  with the component $L_i\subset L$),
  which is an isomorphism in the case where $H^1(Y;\Z)=0$.
  Given an oriented link, and $\x,\y\in\Ta\cap\Tb$, and
  $\phi\in\pi_2(\x,\y)$ is any homology class, then
$$    \Pi(\RelDiff_{\w,\z}(\x,\y))=n_\z(\phi)-n_\w(\phi).
$$
\end{lemma}

\begin{proof}
  The homology class
  $\phi\in\pi_2(\x,\y)$ gives rise to a null-homology of
  ${\underline\epsilon}(\x,\y)$ inside $Y$. This null-homology meets
  the $i^{th}$ component of the
  oriented link $\orL$ with intersection number
  $n_{\z}(\phi)-n_{\w}(\phi)$. The lemma follows at once.
\end{proof}

We can lift ${\underline\epsilon}$ to a map from intersection
points to relative $\SpinC$ structures, generalizing the map
from Subsection~\ref{subsec:IntPtsSpinCStruct}.

Let $(\Sigma,\alphas,\betas,\ws,\zs)$ be a pointed Heegaard diagram
for an oriented link $\orL$. We define the map
\begin{equation}
\label{eq:RelSpinC}
\relspinc_{\ws,\zs}\colon \Ta\cap\Tb\longrightarrow
\RelSpinC(Y,\orL)
\end{equation} as follows. (Note that we write
$\RelSpinC(Y,\orL)$ for the space of relative $\SpinC$ structures on
$Y-\nbd{L}$, cf. Subsection~\ref{subsec:RelSpinC}.)

For this map, we fix a choice as follows. Let $\gamma$ be a gradient
flowline connecting an index zero and index three critical point, and
let $N(\gamma)$ denote a neighborhood of this flowline. One can
construct a nowhere vanishing vector field $\vec{v}$ over $N(\gamma)$,
which has an integral flowline $P$ which enters $N(\gamma)$ from its
boundary, contains $\gamma$ as a subset, and then exits $N(\gamma)$.

Let $f$ be an orientation-preserving involution of $N(\gamma)$. We can
arrange for $-\vec{v}|_{\partial N(\gamma)}$ to agree with
$f^*(\vec{v})|_{\partial(N(\gamma)\cup P)}$. 
Indeed, we can construct $\vec{v}$ in
such a manner that the difference $\vec{v}-f^*(\-\vec{v})$ 
is the Poincar\'e dual of a meridian for $\gamma$, thought of as an
element of $H_1(N(\gamma)-P)$. 

\begin{figure}
\mbox{\vbox{\epsfbox{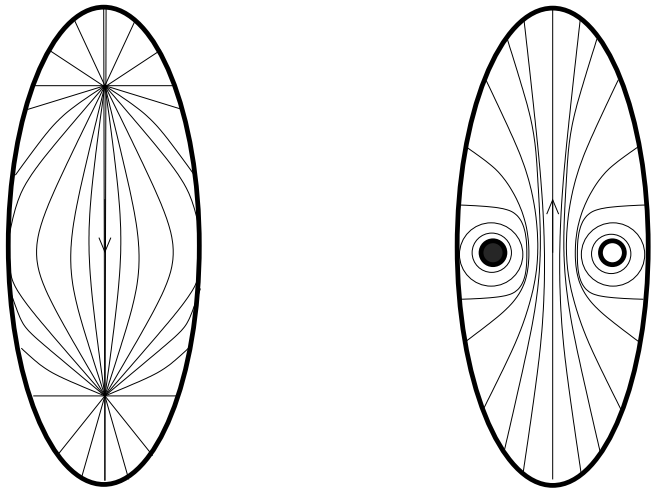}}}
\caption{\label{fig:VectorField}
{\bf{Non-vanishing vector field.}}  On the left, we represent the
gradient vector field in a standard neighborhood of a gradient
flowline $\gamma$ connecting an index zero and three critical point,
by sketching various integral curves. The vector field on the
three-ball is obtained by rotating this vector field
through the axis through the center of the picture.
Replace this vector field with the nowhere vanishing
vector field $\vec{v}$ represented on the right.  Specifically, consider
the vector field represented here (after rotating through the axis in the
center), which vanishes on the oriented meridian represented here by
the pair of circles.  Modify the vector field to point into the page
at the empty circle and out of the page at the dark circle. This gives
the nowhere vanishing vector field $\vec{v}$, with a closed orbit
which coincides with an oriented meridian for $\gamma\subset P$.}
\end{figure}

This vector field is illustrated in Figure~\ref{fig:VectorField}.

Armed with this vector field, we define the map promised in
Equation~\eqref{eq:RelSpinC}. Fix a Morse function $f$ compatible with
the Heegaard diagram $(\Sigma,\alphas,\betas,\ws,\zs)$. Given
$\x\in\Ta\cap\Tb$, consider the flowlines $\gamma_\x$, $\gamma_\w$,
and $\gamma_\z$. We replace the gradient vector field in a neighborhood
of $\gamma_\x$ so as not to vanish there. Similarly, we replace 
the gradient vector field in a neighborhood of $\gamma_\w$ using $\vec{v}$
so that it does not vanish there. In fact, arranging for $P$ to
consist of arcs on $\gamma_\z\cup\gamma_\w$, we obtain in this manner
a vector field on $Y$ which contains $\orL$ as a closed orbit. It is
easy to see that this is equivalent to a vector field on $Y-\nbd{\orL}$
which is a standard non-vanishing vector field on the boundary tori.

\begin{lemma}
  \label{lemma:Filtration} We have that
$$    \spincrel_{\ws,\zs}(\x)-\spincrel_{\ws,\zs}(\y)=
    \PD[{\underline\epsilon}(\x,\y)].$$
where here $\PD$ denotes the Poincar\'e duality map
$$\PD\colon H_1(Y-L)\longrightarrow H^2(Y,L).$$
Indeed, given $\phi\in\pi_2(\x,\y)$, we have that
\begin{equation}
  \label{eq:RelativeDifference}
  \spincrel_{\ws,\zs}(\x)-\spincrel_{\ws,\zs}(\y)
  =\sum_{i=1}^\ell (n_{z_i}(\phi)-n_{w_i}(\phi))\PD[\mu_i],
\end{equation}
where here $\mu_i$ is the meridian for the $i^{th}$ component of $L$,
with its induced orientation from the orientation of $L$.
\end{lemma}

\begin{proof}
The vector fields $\spincrel_{\ws,\zs}(\x)$ and
$\spincrel_{\ws,\zs}(\y)$ differ in a neighborhood of
$\gamma_\x-\gamma_\y$. It is now a local calculation to see that
$\spincrel_{\ws,\zs}(\x)-\spincrel_{\ws,\zs}(\y)=\PD[\gamma_\x-\gamma_\y]$
(compare~\cite[Lemma~\ref{HolDisk:lemma:VarySpinC}]{HolDisk} for the
corresponding statement for $\SpinC$ structures). It is easy to see
that $\gamma_x-\gamma_y$ is homologous to
${\underline\epsilon}(\x,\y)$.

The second remark follows immediately from Lemma~\ref{lemma:Filtration1}.
\end{proof}

The following fact will be useful in studying the symmetry properties
of link Floer homology.

\begin{lemma}
\label{lemma:JAction}
Let $(\Sigma,\alphas,\betas,\ws,\zs)$ be a multiply pointed Heegaard
diagram representing an oriented link $\orL\subset Y$, and 
let $\spincrel_{\ws,\zs}\colon \Ta\cap\Tb\longrightarrow \RelSpinC(Y,L)$
denote the induced map to relative $\SpinC$ structures. Then, $(-\Sigma,\betas,\alphas,\ws,\zs)$
specifies the same link, endowed with the opposite orientation. Let $\spincrel'_{\ws,\zs}\colon \Ta\cap\Tb\longrightarrow \RelSpinC(Y,L)$
denote its induced map. Then, we have that
\begin{equation}
\label{eq:ConjugationFields}
\spincrel_{\ws,\zs}(\x)=J \circ \spincrel'_{\ws,\zs}(\x).
\end{equation}

Also, $(\Sigma,\alphas,\betas,\zs,\ws)$ is another diagram representing the link with the opposite
orientation; let $\spincrel_{\zs,\ws}\colon \Ta\cap\Tb\longrightarrow \RelSpinC(Y,L)$ denote its corresponding
map. Then,
\begin{equation}
\label{eq:OrientationReversal}
\spincrel_{\ws,\zs}(\x)=\spincrel_{\zs,\ws}(\x)+\sum_{i=1}^\ell \PD[\mu_i].
\end{equation}
\end{lemma}

\begin{proof}
  For the first remark, note that if $f$ is a Morse function
  compatible with $(\Sigma,\alphas,\betas,\ws,\zs)$, then $-f$ is a
  Morse function compatible with $(-\Sigma,\betas,\alphas,\ws,\zs)$.
  It is now a straightforward consequence of its definition that if
  $\spincrel_{\ws,\zs}(\x)$ is represented by $\vec{v}$, then
  $-\vec{v}$ represents $\spincrel'_{\ws,\zs}(\x)$. Thus,
  Equation~\eqref{eq:ConjugationFields} follows.
  
  For the second observation, note that the vector fields $\vec{v}$
  and $\vec{w}$ representing $\spincrel_{\ws,\zs}(\x)$ and
  $\spincrel_{\zs,\ws}(\x)$ respectively differ only in a collar
  neighborhood of the boundary. In fact, one can see that in this
  neighborhood, $\vec{w}$ and $\vec{v}$ can be made isotopic away from
  a neighborhood of the meridian, where they point in opposite
  directions. In this way, Equation~\eqref{eq:OrientationReversal}
  follows.
\end{proof}

\subsection{Filling relative $\SpinC$ structures}
\label{subsec:Filling}

Let $\orL\subset Y$ be an oriented link. Fix a component $K_1$ of the
underlying link. We have a natural ``filling map''
$$\Fill{K_1}\colon \RelSpinC(Y,\orL) \longrightarrow
\RelSpinC(Y,\orL-K_1).$$
This is gotten by simply viewing the relative
$\SpinC$ structure on $Y$ relative to $\orL$ as one over $Y$ relative
to $\orL-K_1$. More specifically, if we think of $\RelSpinC(Y,\orL)$
as generated by vector fields which have a closed orbits consisting of
the components of $L$ (traversed with their given orientations), then
we can view these also as relative $\SpinC$ structures which have a
closed orbits consisting of the components of $L-K_1$ (traversed with
their given orientations).

\begin{lemma}
  \label{lemma:COneFill}
  We have the following
$$c_1(\Fill{K_1}(\relspinc))=c_1(\relspinc)+\PD[K_1].$$
\end{lemma}

\begin{proof}
  To verify this, let $F_2$ be an oriented surface in $Y-\nbd{L-K_1}$
  with boundary on $\partial \nbd{L-K_1}$, representing some relative
  homology class $H_2(Y,L-K_1)$. We can assume that $F_2$ meets $K_1$
  transversally.  We can consider the surface $F_2'$ in $Y-\nbd{L}$
  with boundary on $\nbd{L}$ gotten by deleting $\nbd{K_1}\cap F_2$
  from $F_2$.  The homology class of $F_2'$ represents the natural map
  $H_2(Y,L-K_1)\longrightarrow H_2(Y,L)$.
  
  If $\#(K_1\cap F_2)$ denotes the algebraic intersection
  number of $F_2$ with $K_1$ (endowed with the orientation it inherits
  from $\orL$) we have that
  $$\langle c_1(\Fill{K_1}(\relspinc)), [F_2]\rangle
  =\langle c_1(\relspinc),[F_2']\rangle + \#(K_1\cap F_2),$$
  which is equivalent to the claim in the lemma.
\end{proof}

\section{Definition of Heegaard Floer homology for multi-pointed Heegaard diagrams}
\label{sec:DefHFL}

\subsection{Heegaard diagrams for three-manifolds}

For convenience, we work always with Floer homology with coefficients
in $\Field=\Zmod{2}$. We also consider the case where the ambient
manifold is a rational homology three-sphere, and hence
Proposition~\ref{prop:Admissibility} applies.

Let $(\Sigma,\alphas,\betas,\ws)$ be an $\ell$-pointed balanced
weakly admissible Heegaard
diagram for a rational homology three-sphere
$Y$. We define $\CFm(\Sigma,\alphas,\betas,\ws)$ to be the
free module over the polynomial algebra $\Field[U_1,...,U_\ell]$
generated by the intersection points $\Ta\cap\Tb$ inside
$\Sym^{g+\ell-1}(\Sigma)$. There is a module homomorphism
$$\partial^-\colon \CFm(\Sigma,\alphas,\betas,\ws)\longrightarrow
\CFm(\Sigma,\alphas,\betas,\ws)$$
defined by
\begin{equation}
\label{eq:DefD}
\partial^- \x = \sum_{\y\in\Ta\cap\Tb}
\sum_{\{\phi\in\pi_2(\x,\y)\big| \Mas(\phi)=1\}}
\#\UnparModFlow(\phi)\cm 
U_1^{n_{w_1}(\phi)}\cm ...\cm U_\ell^{n_{w_\ell}(\phi)}
\otimes \y.
\end{equation}
Here, as usual, $\ModFlow(\phi)$ denotes the moduli space of
pseudo-holomorphic representatives of the given homology class of
Whitney disks, and $\UnparModFlow(\phi)$ denotes the quotient of this
moduli space by the action of $\R$.  Also, $\Mas(\phi)$ denotes the
expected dimension of $\ModFlow(\phi)$ (i.e. the Maslov index of
$\phi$).  (Indeed, we will find it useful later to use Lipshitz's
cylindrical formulation~\cite{LipshitzCyl} instead; we return to this point
in Section~\ref{sec:Analysis}.)

We will establish the following in
Subsection~\ref{subsec:BoundaryDegenerations}:

\begin{prop}
  \label{prop:MaslovIndexFormula}
  For an $\ell$-pointed, balanced, Heegaard diagram for a rational
  homology three-sphere $Y$, for any $\phi,\phi'\in\pi_2(\x,\y)$, we
  have that
\begin{equation}
  \label{eq:MasIndexFormula}
  \Mas(\phi)-\Mas(\phi')=2\sum_{i=1}^\ell (n_{w_i}(\phi)-n_{w_i}(\phi')).
\end{equation}
\end{prop}

\begin{lemma}
  For an $\ell$-pointed balanced weakly-admissible Heegaard diagram
  for a three-manifold $Y$ and any
  $\x\in\Ta\cap\Tb$, the right-hand-side
  of Equation~\eqref{eq:DefD} consists of only finitely many different
  non-zero
  terms.
\end{lemma}

\begin{proof}
  First, we prove that the coefficient of $U_1^{a_1}\cm...\cm
  U_\ell^{a_\ell}$ in $\y$ is finite, for any given ${\mathbf
    a}=a_1\times...\times a_\ell\in\Z^{\ell}$.  But this follows
  readily from the fact that for any two homology classes
  $\phi,\phi'\in\pi_2(\x,\y)$ with $n_\w(\phi)=n_\w(\phi')$,
  $\phi-\phi'$ is a periodic domain, which hence must have both
  positive and negative coefficients. It follows from this that there
  are at most finitely many homology classes $\phi$ with
  $n_\w(\phi)={\mathbf a}$ and $\cald(\phi)\geq 0$, and hence
  which can support holomorphic representatives.
  
  If $\phi$ admits a holomorphic representative, then
  $n_{w_i}(\phi)\geq 0$ for all $i$; from
  Proposition~\ref{prop:MaslovIndexFormula}, it follows at once that
  for $\phi\in\pi_2(\x,\y)$ with $\Mas(\phi)=1$, the quantity $\sum_i
  n_{w_i}(\phi)$ depends only on $\x$ and $\y$. Thus, given $\x$ and
  $\y$ it follows that there are only finitely many different
  possibilities for $U_1^{a_1}\cdot...\cdot
  U_\ell^{a_\ell}$ with non-zero coefficient
  on the right-hand-side of Equation~\eqref{eq:DefD}.
\end{proof}

\begin{lemma}
\label{lemma:DSquaredZero}
The map $\partial^-$ is a differential on $\CFm(\Sigma,\alphas,\betas,\ws)$.
\end{lemma}

The above lemma is proved in
Section~\ref{sec:IdentifyHeegaardFloers}, where we also
prove that the homology groups are identified with the usual Heegaard
Floer homology of~\cite{HolDisk}:

\begin{theorem}
\label{thm:InvarianceHFm}
The chain homotopy type of the complex 
$\CFm(\Sigma,\alphas,\betas,\ws)$ over the polynomial ring
$\Field[U_1]$  is 
a three-manifold invariant; indeed, its homology coincides
with the Heegaard Floer homology group $\HFm(Y)$.
\end{theorem}

A more precise version of the above theorem can be stated, which respects
the splitting of $\HFm(Y)$ according to its various $\SpinC$ structures.
We do not belabour this point now, as the primary example we have in 
mind here is the case
where $Y=S^3$ (which has a unique $\SpinC$ structure).

There is a simpler variant of the above construction, where we set
$U_1=0$, formally -- i.e. we consider the free module over the polynomial
algebra $\Field[U_2,...,U_\ell]$ generated by $\Ta\cap\Tb$, 
endowed with the differential:

$$
{\partial} \x = \sum_{\y\in\Ta\cap\Tb}
\sum_{\{\phi\in\pi_2(\x,\y)\big| n_{w_1}(\phi)=0, \Mas(\phi)=1\}}
\#\UnparModFlow(\phi)\cm 
U_2^{n_{w_2}(\phi)}\cm ...\cm U_\ell^{n_{w_\ell}(\phi)}
\otimes \y
$$

It is easy to see that the arguments from
Theorem~\ref{thm:InvarianceHFm} also show that the homology of this
chain complex calculates $\HFa(Y)$.

The simplest variant counts holomorphic disks which are disjoint from
all the $w_i$, to obtain a complex $\CFa(\Sigma,\alphas,\betas,\ws)$
i.e. specializing $\CFm(\Sigma,\alphas,\betas,\ws)$
to the case where $U_i=0$ for all $i$. Explicitly, this is the chain 
complex of $\Field$-vector spaces spanned by the intersection points
of $\Ta\cap\Tb$, and equipped with the differential
$$
{\partial} \x = \sum_{\y\in\Ta\cap\Tb}
\sum_{\{\phi\in\pi_2(\x,\y)\big| n_{\ws}(\phi)=0, \Mas(\phi)=1\}}
\#\UnparModFlow(\phi)\cm 
 \y.
$$

For this complex, we have the following:

\begin{theorem}
\label{thm:InvarianceHFaa}
The complex
$\CFa(\Sigma,\alphas,\betas,\ws)$ calculates
$\HFa(Y\#(\#^{\ell-1}(S^2\times S^1)))$.
\end{theorem}

\subsection{Links and multi-filtrations of $\HFa$}
\label{subsec:LinkFiltration}

The multi-pointed Heegaard diagram for an oriented link $\orL$ endows
the chain  complex $\CFm(\Sigma,\alphas,\betas,\ws)$ with a 
relative $\Z^\ell$-filtration, as follows.

Let $(\Sigma,\alphas,\betas,\ws,\zs)$ be the Heegaard diagram for an
oriented link in a three-manifold with $H^1(Y;\Z)=0$.  In 
Subsection~\ref{subsec:RelSpinCLinks}, we defined a map
$$\spincrel_{\ws,\zs}\colon \Ta\cap\Tb\longrightarrow
\RelSpinC(Y,L).$$

We extend the above-defined function to generators of the
chain complex $\CFm(\Sigma,\alphas,\betas,\ws)$ of the
form $U_1^{a_1}\cdot...\cdot U_\ell^{a_\ell}\cdot \x$, using the formula
\begin{equation}
\label{eq:ExtendRelSpinC}
\spincrel_{\w,\z}(U_1^{a_1}\cdot...\cdot U_\ell^{a_\ell}\cdot \x) =\spincrel_{\w,\z}(\x)-a_1\cdot
\PD[\mu_1]-...-a_\ell\cdot \PD[\mu_\ell],
\end{equation}
where here $\{\mu_i\}_{i=1}^\ell$ are meridians for the link
(compatible with its given orientation).

\begin{lemma}
  Let $\orL\subset Y$ be an oriented link in a three-manifold $Y$ with
  $H_1(Y;\Z)=0$, and consider the corresponding identification
  $\Z^\ell\cong H^2(Y,L)$, given by $(a_1,...,a_\ell)\mapsto
  \sum_{i=1}^\ell a_i\cm \PD[\mu_i]$, making $\RelSpinC(Y,L)$ into an
  affine space for $\Z^\ell$.  Then, the function $\spincrel$ (as
  defined in Subsection~\ref{subsec:RelSpinC}, and extended in
  Equation~\eqref{eq:ExtendRelSpinC} above) induces a $\RelSpinC(Y,L)$-filtration on the chain complex
  $\CFm(\Sigma,\alphas,\betas,\ws,\zs)$ (endowed with the differential
  from Equation~\eqref{eq:DefD}).
\end{lemma}

\begin{proof}
  Suppose $U_1^{a_1}\cm...\cm U_\ell^{a_\ell}\cm \y$ appears in
  $\partial \x$ with non-zero coefficient. We must prove that then
  $\spincrel(\x)\geq \spincrel(U_1^{a_1}\cm...\cm U_\ell^{a_\ell}\cm \y)$. But
  in this case, there is a $\phi\in\pi_2(\x,\y)$ with a holomorphic
  representative, and hence $n_\w(\phi), n_\z(\phi)\geq 0$. On the
  other hand, $a_i=n_{w_i}(\phi)$, and hence
  $$\spincrel_{\w,\z}(\x)-\spincrel_{\w,\z}(U_1^{a_1}\cm...\cm U_\ell^{a_\ell}\cm \y)
  =\spincrel_{\w,\z}(\x)-\spincrel_{\w,\z}(\y) + n_{\w}(\phi)
  =n_{\z}(\phi), $$
  in view of Equation~\eqref{eq:RelativeDifference}.
\end{proof}

The following will be verified in Section~\ref{sec:LinkInvariants}:

\begin{theorem}
  \label{thm:InvarianceHFLFilt}
  The $\RelSpinC(Y,L)$-filtered chain homotopy type of the chain complex
  $\CFm(\alphas,\betas,\ws,\zs)$ of $\Z[U_1,...,U_\ell]$-modules 
  is an invariant of the underlying
  oriented link.
\end{theorem}

\begin{defn}
  \label{def:CFLm}
  The $\RelSpinC(Y,L)$-filtered chain homotopy type of the chain complex associated
  to a link will be denoted
  $\CFLm(Y,\orL)$. The homology of the associated graded object is denoted
  $$\HFLm(Y,\orL)=\bigoplus_{\relspinc\in\RelSpinC(Y,L)}\HFLm(Y,\orL,\relspinc).$$
\end{defn}

Explicitly, the associated graded object is generated as a $\Field[U_1,...,U_\ell]$-module by 
intersection points $\Ta\cap\Tb$, endowed with a differential differential akin to Equation~\eqref{eq:DefD},
only now we sum over those $\phi\in\pi_2(\x,\y)$ with $\Mas(\phi)=1$ and also $n_{\zs}(\phi)=0$. Indeed,
the summand in grading $\relspinc$ is generated by symbols $U^{a_1}\cm...\cm U^{a_\ell}\otimes \x$,
where $a_1,...,a_\ell$ are non-negative integers and $\x\in\Ta\cap\Tb$, satisfying the constraint that
$$ \relspinc_{\ws,\zs}(\x)-a_1\PD[\mu_1]-...-a_\ell\PD[\mu_\ell]=\relspinc. $$ The homology of this complex
is the group

$\HFLm(Y,\orL,\relspinc)$.

We can also set $U_1=...=U_\ell=0$, to obtain a filtration of 
a chain complex which, according to Theorem~\ref{thm:InvarianceHFaa},
calculates $\HFa$ of $Y\#(\#^{\ell-1}(S^2\times S^1))$. More concretely,
we consider the chain complex $\CFa(\Sigma,\alphas,\betas,\ws)$ generated
over $\Field$ by intersection points $\Ta\cap\Tb$, endowed with the module homomorphism
\begin{equation}
\label{eq:DefDa}
\partial \x = \sum_{\y\in\Ta\cap\Tb}
\sum_{\{\phi\in\pi_2(\x,\y)\big| \Mas(\phi)=1, n_{\ws}(\phi)=0\}}
\#\UnparModFlow(\phi)\cm  \y.
\end{equation}

As in the case of $\CFLm$, the function $\relspinc_{\ws,\zs}$ endows
$\CFLa(Y,\orL)$ with a filtration.

\begin{defn}
  \label{def:CFLa}
  The $\RelSpinC(Y,L)$-filtered chain homotopy type of the chain
  complex associated to a link will be denoted $\CFLa(Y,\orL)$. The
  homology of the associated graded complex is the link invariant
  $\HFLa(Y,\orL)$.
\end{defn}

More precisely, for $\relspinc\in\SpinC(Y,\orL)$,  
$\HFLa(Y,\orL, \relspinc)$ is the homology of the chain complex generated by
$\x\in\Ta\cap\Tb$ with $\spincrel_{\ws,\zs}(\x)=\relspinc$, endowed with the differential
\begin{equation}
\label{eq:DefDaa}
\partial \x = \sum_{\y\in\Ta\cap\Tb}
\sum_{\{\phi\in\pi_2(\x,\y)\big| \Mas(\phi)=1, n_{\ws}(\phi)=n_{\zs}(\phi)=0\}}
\#\UnparModFlow(\phi)\cm  \y.
\end{equation}

In the introduction, no mention was made of relative $\SpinC$
structures; rather, link Floer homology was described as a group graded by
elements of $\iH$.  For the case of links in $S^3$, the equivalence of
these two points of view is given as follows. Given an element $h=\sum
a_i [\mu_i]\in \iH$, there is a unique relative $\SpinC$ structure
$\relspinc$ with the property that
$$c_1(\relspinc)- \sum_{i=1}^{\ell} \PD[\mu_i] = 2 \PD[h].$$
The group $\HFLa(L,h)$ from the introduction, then, is the 
group $\HFLa(L,\relspinc)$. This convention is quite
natural, as we shall see in Section~\ref{subsec:Notation} below.

In practice, it can be taxing to calculate relative $\SpinC$
structures. It is much simpler, rather, to define the link filtration
on relative terms, declaring  that 
\begin{eqnarray}
\label{eq:RelOrdering}
\x>\y&{\text{if and only if}}&
n_{\zs}(\phi)\geq n_{\ws}(\phi),
\end{eqnarray} where here $\phi$ is any element in
$\pi_2(\x,\y)$, for example in the case where $H_1(Y;\Z)=0$.  (The
equivalence of this with our earlier point of view is a direct consequence of 
Equation~\eqref{eq:RelativeDifference}.)  This determines the filtration
only up to an overall shift (by an element of $H^2(Y,L;\Z)$), but this indeterminacy can be removed
using the symmetry properties of the link invariant, cf.
Section~\ref{sec:Symmetry} below.

\section{Analytic input}
\label{sec:Analysis}

We will be concerned in this section with gluing results for
pseudo-holomorphic disks.  Consider an $\ell$-pointed Heegaard diagram
$(\Sigma,\alphas,\betas,\{w_1,...,w_\ell\})$, and let $d=g+\ell-1$.
Given $\x,\y\in\Ta\cap\Tb$ and a homology class of Whitney disk
$\phi\in\pi_2(\x,\y)$, we can form the moduli space $\ModFlow(\phi)$.

Recall that these moduli spaces have Gromov compactifications,
cf.~\cite{Gromov}, \cite{McDuffSalamon}, \cite{FloerLag}, \cite{FOOO}.
For a given moduli space of pseudo-holomorphic Whitney disks,
these Gromov compactifications include possibly 
moduli spaces of pseudo-holomorphic
Whitney disks connecting other intersection points,  moduli
spaces of pseudo-holomorphic spheres, and finally also
moduli spaces of further degenerate disks called
{\em boundary degenerations}. More formally,
given a point $\x\in\Ta$, we let $\pi_2^{\beta}(\x)$ denote the space
of homology classes of maps $$\left\{ u\colon [0,\infty)\times
  \R\longrightarrow \Sym^g(\Sigma)\Bigg|
\begin{array}{l}
u(\{0\}\times \R)\subset \Tb \\
\lim_{z\goesto \infty}u(z)=\x 
\end{array}
\right\}.$$
Such a map is called a $\beta$-boundary degeneration.
We let $\pi_2^{\alpha}(\x)$ denote the set of $\alpha$-boundary
degenerations, defined analogously.

Loosely speaking, Gromov's compactness theorem states that a sequence
of pseudo-holomorphic curves representing $\phi$, has a subsequence
which converges locally to a ``broken flow-line'', consisting of
collection of pseudo-holomorphic flow-lines $\{\phi_i\}$, a collection
of $\alpha$- and $\beta$-boundary degenerations $\{\psi_j\}$, and
finally a collection of pseudo-holomorphic spheres $\{S_k\}$ with
$$\sum\cald(\phi_i)+\sum\cald(\psi_j)+\sum\cald(S_k)=\cald(\phi).$$

We will also find it necessary at several future points to study ends
of moduli spaces as the Heegaard surface is degenerated. We turn to
this more formally as follows.

\subsection{Gluing moduli spaces}
\label{subsec:FiberedProducts}

Let $(\Sigma_1,\alphas_1,\betas_1,z_1)$ and
$(\Sigma_2,\alphas_2,\betas_2,z_2)$ be a pair of Heegaard diagrams,
where here $\alphas_i$ and $\betas_i$ are $d_i$-tuples of attaching
circles in $\Sigma_i$. We can form their connected sum at $z_1\in
\Sigma_1$ and $z_2\in \Sigma_2$ to obtain a new surface
$\Sigma=\Sigma_1\#\Sigma_2$, endowed with sets of attaching circles
$\alphas_1\cup\alphas_2$ and $\betas_1\cup\betas_2$. We will assume
that $d_i\geq g_i$, and write $d=d_1+d_2$ and $g=g_1+g_2$.

We will need here descriptions of the moduli spaces of flowlines in
the connected sum diagram, in terms of moduli spaces for the two
pieces.

Fix points $\x_i,\y_i\in{\mathbb T}_{\alpha_i}\cap{\mathbb
T}_{\beta_i}\subset \Sym^{d_i}(\Sigma_i)$, which in turn give rise to
points $\x_1\times \x_2$, $\y_1\times\y_2$ in
$\Sym^{d_1+d_2}(\Sigma_1\#\Sigma_2)$.  Fix a pair of homology
classes of Whitney disks
$\phi_i\in\pi_2(\x_i,\y_i)$ with
$n_{z_1}(\phi_1)=n_{z_2}(\phi_2)$. These can be combined naturally to
form a homology class
$\phi_1\#\phi_2\in\pi_2(\x_1\times\x_2,\y_1\times\y_2)$.
Specifically, the local multiplicities of $\phi_1\#\phi_2$ at each
domain $\cald\subset \Sigma_1\#\Sigma_2$ is the local multiplicity of
$\phi_i$ at the corresponding $\cald\subset \Sigma_1$ or $\Sigma_2$.

Conversely, each homology class
$\phi\in\pi_2(\x_1\times\y_1,\x_2\times\y_2)$ can be uniquely
decomposed as $\phi=\phi_1\#\phi_2$ for some pair of
$\phi_i\in\pi_2(\x_i,\y_i)$ with $n_{z_1}(\phi_1)=n_{z_2}(\phi_2)$.

Moreover, given complex structures $j_1$ and
$j_2$ on $\Sigma_1$ and $\Sigma_2$, we can form a complex structure 
$J(T)$ with neck-length $T$ as follows: find conformal 
disks $D_1$ and $D_2$ about
$z_1$ and $z_2$, and form the connected sum
$$\Sigma(T)=(\Sigma_1-D_1)\# \left([-T-1,T+1]\times S^1\right) \# (\Sigma_2-D_2),$$
under identifications $\partial D_1\cong \{-T\}\times S^1$
and $\{T\}\times S^1 \cong \partial D_2$. As $T\goesto\infty$, the
conformal structure on $\Sigma(T)$ converges to
the nodal curve $\Sigma_1\vee\Sigma_2$.

\begin{theorem}
\label{thm:GenericFiberedProduct}
Fix diagrams $(\Sigma_i,\alphas_i,\betas_i,z_i)$ for $i=1,2$ as above,
where here $\alphas_i$ and $\betas_i$ are $d_i$-tuples of attaching
circles.  Given a homology class $\phi=\phi_1\#\phi_2$ for the
connected sum of the two diagrams, we have that
$$\Mas(\phi)=\Mas(\phi_1)+\Mas(\phi_2)-2k,$$ 
where $k=n_{z_1}(\phi)=n_{z_2}(\phi)$. Suppose that
$\ModFlow(\phi)\neq \emptyset$
for a sequence of almost-complex structures $J(T_i)$ 
with $T_i\goesto \infty$. Then, the moduli spaces of broken pseudo-holomorphic
flowlines representing $\phi_1$ and $\phi_2$ (i.e. the Gromov
compactifications of these two moduli spaces) are non-empty. Finally,
suppose that $\Mas(\phi_1)=1$, $\Mas(\phi_2)=2k$, and also that $d_2>g_2$;
and consider the
maps
\begin{eqnarray*}
\rho_1\colon \ModFlow(\phi_1) \longrightarrow \Sym^{k}(\CDisk) 
&{\text{and}}&
\rho_2\colon \ModFlow(\phi_2) \longrightarrow \Sym^{k}(\CDisk),
\end{eqnarray*}
where here
$$\rho_i(u)=u^{-1}(\{z_i\}\times \Sym^{d_i-1}(\Sigma_i)).$$
If the fibered product of $\ModFlow(\phi_1)$ and $\ModFlow(\phi_2)$
$$\ModFlow(\phi_1)\times_{\Sym^k(\CDisk)} \ModFlow(\phi_2)
=\{u_1\times u_2 \in\ModFlow(\phi_1)\times \ModFlow(\phi_2)\big|
\rho_1(u_1)=\rho_2(u_2)\}$$
is a smooth manifold, then there is an identification of this moduli
space with the moduli space $\ModFlow(\phi)$, for sufficiently 
long connected sum length.
\end{theorem}

There are three assertions in the above theorem: one concerns the
Maslov index, one the existence of weak limits, and the third is a
gluing result. They are arranged in order of difficulty; the third
requires the most work.  One could approach this problem from the
point of view of degenerating $\Sym^{d}(\Sigma_1\#\Sigma_2)$ as the
connected sum degenerates into $\Sigma_1\vee\Sigma_2$, following the
approach to stabilization invariance from~\cite{HolDisk}.  In this
case, the limiting symplectic space is
$\Sym^{d}(\Sigma_1\vee\Sigma_2)$, which has a fairly complicated
singular set (consisting of those $d$-tuples where one element is the
singular point $p\in \Sigma_1\vee\Sigma_2$). In particular, the
singular set is, in itself, not a symplectic manifold but rather a
singular space whose singularities consist of those $d$-tuples where
at least two elements are the point $p$.  Under the present
circumstances, holomorphic disks we wish to resolve -- whose boundary
lies from in the the top stratum
$\Sym^{d_1}(\Sigma_1)\times\Sym^{d_2}(\Sigma_2)\subset
\Sym^{d}(\Sigma_1\vee\Sigma_2)$ -- meet the singular stata (consisting of
those tuples where at least one point is the singular point $z\in
\Sigma_1\vee\Sigma_2$) in a complex codimension one subset, i.e. where
at least two coordinates agree with this singular point.

A much simpler approach can be given using Lipshitz's cylindrical
reformulation of Heegaard Floer homology, cf.~\cite{LipshitzCyl}.
With this reformulation, then, the gluing problem takes place 
in a four-manifold, along a singular set which is a manifold,
(placing it on roughly an equal footing with the proof of stabilization
invariance for cylindrical reformulation, cf.~\cite{LipshitzCyl}). 
This kind of degeneration has been extensively studied in the literature,
cf.~\cite{IonelParker}, \cite{LiRuan}, \cite{EGH}, \cite{BEHWZ},
and of course~\cite{LipshitzCyl}. We
turn to this approach, first recalling the basic set-up of Lipshitz's picture.

\subsection{Lipshitz's cylindrical formulation of Heegaard Floer homology}

The starting point of Lipshitz's formulation is that a holomorphic
disk $u\colon\CDisk\longrightarrow \Sym^d(\Sigma)$ corresponds to a
holomorphic curve in $\CDisk\times\Sigma$, which, for any
$p\in\CDisk$, meets the fiber $\Sigma\times \{p\}$ in the 
set of $d$ points $u(p)$.  Thus, one could reformulate the chain
complex defining Heegaard Floer homology as counting certain
pseudo-holomorphic surfaces in $\CDisk\times \Sigma$. 

This can be made more precisely as follows.
Consider 
the four-manifold 
$$W=\Sigma\times[0,1]\times\R,$$ equipped
with two projection maps
\begin{eqnarray*}
  \pi_{\Sigma}\colon W \longrightarrow \Sigma
  &{\text{and}}&
  \pi_{\CDisk}\colon W\longrightarrow [0,1]\times\R.
\end{eqnarray*}
(As usual, we think of the unit disk in the complex plane $\CDisk$
as the conformal compactification of the infinite strip $[0,1]\times \R$,
obtained by adding points $\pm i$ at infinity.)
Endow $W$ with an almost-complex structure $J$ tamed by a natural split
symplectic form on $W$, which is translation invariant in the
$\R$-factor, and for which the projection $\pi_{\CDisk}$ is
a pseudo-holomorphic map.
For example, the product complex structure -- which will be called
a {\em split} complex structure -- satisfies these conditions,
but it is often useful to perturb this. However, to ensure
positivity, it is convenient to choose points $z\in \Sigma$,
and require that $J$ is split in a neighborhood of $z$ 
(in $\Sigma$) times $[0,1]\times \R$. Such an almost-complex structure
$J$ is called {\em split near $z$}.

Consider next a Riemann surface $S$ with boundary,
$d$ ``positive'' punctures $\{p_1,...,p_d\}$ and $d$ ``negative''
punctures $\{q_1,...,q_d\}$ on its boudnary. 

Lipshitz considers pseudo-holomorphic 
maps $$\uL\colon S\longrightarrow W,$$
satisfying the following conditions:
\begin{itemize}
    \item $\uL$ is smooth.
    \item $\uL(\partial S)\subset \left(\alphas\times\{1\}\times \R\right)\cup
      \left(\betas\times\{0\}\times\R\right)$.
    \item No component in the image of $\uL(S)$ is contained in a fiber of 
      $\pi_{\CDisk}$.
    \item For each $i$, $\uL^{-1}(\alpha_i\times\{1\}\times \R)$
      and $\uL^{-1}(\beta_i\times\{0\}\times \R)$
      consist of exactly one component of $\partial S-\{p_1,...,p_d,q_1,...,q_d\}$
    \item The energy of $\uL$ is finite.
    \item $\uL$ is an embedding.
    \item Any sequence of points in $S$ converging to $q_i$ resp. $p_i$
	is mapped under $\pi_{\CDisk}$ to a sequence of points
	whose second coordinate converges to $\{-\infty\}$ resp $\{+\infty\}$.
\end{itemize}
We call holomorphic curves of this type {\em cylindrical flow-lines}.
Thinking of the complex disk $\CDisk$ as a compactification of
$[0,1]\times \R$, a map $\uL$ as above
can be extended continuously to a map of the closure of $S$ into
${\overline W}=\Sigma\times\CDisk$.  We say that $\uL$ connects $\x$
to $\y$ if the image of this extension meets $\Sigma\times \{-i\}$
in the points $\x\times\{-i\}$ and it meets $\Sigma\times\{i\}$ in the
points $\y\times \{i\}$.

Projecting such a map $\uL$ onto $\Sigma$, we obtain a relative
two-chain in $\Sigma$ relative to $\alphas\cup\betas$, whose local
multiplicity at some point $z\in\Sigma$ is given by the intersection
number $$\nL_z(\uL)=\#\Big(\uL\cap
\left(\{z\}\times[0,1]\times\R\right)\Big).$$ Conversely, given
$\phi\in\pi_2(\x,\y)$, let $\ModFlowL(\phi)$ denote the moduli space
of cylindrical flow-lines $\uL$ which induce the same two-chain as
$\phi$. 

It is sometimes also useful to consider the analogue of boundary degenerations
in this cylindrical context. 

\begin{defn}
Consider a Riemann surface $S$ with
boundary and $d$ punctu      
res $\{p_1,...,p_d\}$ on its boundary. Consider now
pseudo-holomorphic maps ${\widetilde u}\colon S \longrightarrow \Sigma
\times [0,\infty)\times \R$ which are finite energy, smooth
embeddings, sending the boundary of $\Sigma$ into $\betas$, containing
no component in the fiber of the projection to $[0,\infty)\times\R$,
so that each component of ${\widetilde
  u}^{-1}(\beta_i\times\{0\}\times\R)$ consists of exactly one
component of $\partial S-\{p_1,...,p_d\}$. Such a map is called
a {\em cylindrical boundary degeneration}. For such a map, the point
at infinity is mapped into a fixed $\x\in\Tb$. These maps can be
organized into moduli spaces of $\ModDegL(\psi^\beta)$ indexed by homology
classes $\psi\in H_2(\Sigma,\betas)$. A corresponding definition can also
be made with $\alphas$ playing the role of $\betas$.
\end{defn}

In~\cite{LipshitzCyl}, Lipshitz sets up a theory analogous to Heegaard
Floer homology (in the case where $d=g$), counting elements of
$\ModFlowL(\phiL)$. In setting up this theory, he establishes the
necessary transversality properties: for a generic choice of
almost-complex structure $J$ over $W$ as above, all the moduli spaces
$\ModFlowL(\phiL)$ with $\Mas(\phiL)<0$ are empty; non-empty moduli
spaces $\ModFlowL(\phiL)$ with $\Mas(\phiL)=0$ consist of constant
flowlines, and moduli spaces with $\Mas(\phi)=1$ are smooth
one-manifolds. He also shows that these moduli spaces have the
necessary Gromov compactifications analogous to those for the moduli
spaces of holomorphic disks in a more traditional Lagrangian set-up.
Indeed, in~\cite{LipshitzCyl} Appendix A, Lipshitz establishes
identifications $\ModFlowL(\phi)\cong \ModFlow(\phi)$ for suitably
generic choices of almost-complex structures, in cases where
$\Mas(\phi)=1$. In particular, if we consider the map defined as in
Equation~\eqref{eq:DefD}, only using moduli spaces
$\ModFlowL(\phi)/\R$ in place of $\ModFlow(\phi)/\R$, then the two
maps actually agree.  Although Lipshitz considers
the case where $d=g$, the logic here applies immediately to prove
the corresponding result in the case where $d>g$, as well.

With this cylindrical formulation, now the denegeration considered in
Theorem~\ref{thm:GenericFiberedProduct} becomes more transparent.
Specifically, the degeneration of $\Sigma$ to $\Sigma_1\vee\Sigma_2$
corresponds to a generation of $W$ into
$W_1=\Sigma_1\times[0,1]\times\R$ and
$W_2=\Sigma_2\times[0,1]\times\R$, two symplectic manifolds which meet
along the hypersurface $[0,1]\times\R$, joined along
$\{p_1\}\times[0,1]\times\R\subset W_1$ and
$\{p_2\}\times[0,1]\times\R\subset W_2$. (Compare
also~\cite{IonelParker}, \cite{LiRuan}, \cite{EGH}, \cite{BEHWZ}.
Finally, observe that this is precisely the set-up which Lipshitz uses
in his proof of stabilization invariance for the cylindrical theory,
see \cite{LipshitzCyl}).  

More precisely, we start with almost-complex structures $J_1$ and
$J_2$ on $W_1$ and $W_2$ and neighborhoods $D_1$ and $D_2$ of $z_1$
and $z_2$.  We assume that $J_1$ and $J_2$ which are split on
$D_1\times[0,1]\times\R$ and $D_2\times [0,1]\times\R$ respectively.
From this, we construct a complex structure $J(T)$ on
$$W(T)=\Sigma(T)\times[0,1]\times\R,$$
which agrees with $J_1$ near $(\Sigma_1-D_1)\times[0,1]\times \R$,
and $J_2$ near $(\Sigma_2-D_2)\times[0,1]\times\R$,  and which 
is split over $S^1\times [-T-1,T+1]\times [0,1]\times\R$.

\begin{defn}
        \label{def:PreGlued}
        A {\em pre-glued flowline representing the homology class
          $\phi=\phi_1\#\phi_2\in\pi_2(\x_1\times\y_1 ,\x_2\times
          \y_2)$} is a pair of cylindrical flow-lines ${\widetilde
          u}_1\in\ModFlowL(\phiL_1)$ and ${\widetilde
          u}_2\in\ModFlow(\phiL_2)$ satisfying the matching condition
        $$(\pi_{\CDisk}\circ\uL_1)
        \left(\pi_{\Sigma_1}\circ\uL_1)^{-1}(\{z_1\})\right) =
        (\pi_{\CDisk}\circ\uL_2)\left(\pi_{\Sigma_2}\circ
          \uL_2)^{-1}(\{z_2\})\right).$$  Similarly, a {\em pre-glued
          $\alpha$-boundary degeneration} is a pair of
        $\alpha$-boundary degenerations ${\widetilde
          v}_1\in\ModFlow(\phi_1)$ and ${\widetilde
          v}_2\in\ModFlow(\phi_2)$ satisfying the analogous matching
        condition $(\pi_{\CDisk}\circ {\widetilde
          u})\left(\pi_{\Sigma_1}\circ {\widetilde
            v}_1)^{-1}(\{z_1\})\right) = (\pi_{\CDisk}\circ
        {\widetilde u}\left(\pi_{\Sigma_2}\circ {\widetilde
            v}_2)^{-1}(\{z_2\})\right)$. A similar definition can also
        be made for $\beta$-boundary degenerations.
\end{defn}

The curves $\alpha_1\cup...\cup \alpha_{g+\ell-1}$ divide $\Sigma$ into
$\ell$ regions $A_1,...,A_\ell$. 
Choose reference points $w_i$, one in each $A_i$.

It will be useful to have the following:

\begin{lemma}
  \label{lemma:MaslovBoundaryDegenerations}
  Given ${\widetilde \psi}\in\pi_2^\beta(\x)$, we have that 
  $\Mas({\widetilde\psi})=2\sum_{i=1}^\ell n_{w_i}({\widetilde \psi})$.
\end{lemma}

\begin{proof}
  This follows readily from the excision principle for the 
  the linearized $\dbar$ operator, to reduce to the case of a disk.
  (See the proof of Theorem~\ref{thm:CountBoundaryDeg} for a more
  detailed discussion of a related problem.)
\end{proof}

We interrupt now our path to Theorem~\ref{thm:GenericFiberedProduct},
paying off an earlier debt, supplying the following quick consequence
of the above lemma:

\vskip.2cm
\noindent{\bf{Proof of Proposition~\ref{prop:MaslovIndexFormula}}}
According to Equation~\eqref{eq:GroupPerDoms}, the homology classes of
$\phi$ and $\phi'$ differ by the juxtaposition of a boundary
degeneration with $\cald(\psi)=\sum_i (n_{w_i}(\phi)-n_{w_i}(\phi'))$,
whose index, according to
Lemma~\ref{lemma:MaslovBoundaryDegenerations}, is given by
$2\sum_i (n_{w_i}(\phi)-n_{w_i}(\phi'))$.  The result now follows
from the additivity of the index under juxtaposition.
\qed \vskip.2cm

We give another consequence of Lemma~\ref{lemma:MaslovBoundaryDegenerations}.
By arranging for the almost-complex structure on $W$ to be split near
neighborhoods of all the $w_i$, one can arrange for the usual
positivity principle to hold: if a moduli space $\ModDeg(\psi^\beta)$
is non-empty, then all the $n_{w_i}(\psi^\beta)\geq 0$.  It follows
from this, together with
Lemma~\ref{lemma:MaslovBoundaryDegenerations}, that if $\psi$ is a
homology class of boundary degenerations which contains a non-constant
pseudo-holomorphic representative, then $\Mas(\psi)\geq 2$. This 
principle is used in the following:

\vskip.2cm
\noindent{\bf{Proof of the cylindrical analogue of
Theorem~\ref{thm:GenericFiberedProduct}}.}
We turn to Theorem~\ref{thm:GenericFiberedProduct}, using 
moduli spaces of cylindrical flowlines in place of pseudo-holomorphic
Whitney disks. In this context, the restriction maps are to be replaced
by maps
$$\rhoL_i\colon \ModFlowL(\phiL_i)\longrightarrow \Sym^k(\CDisk)$$
given by
$$\rhoL_i(\uL)=(\pi_{\CDisk}\circ \uL)\left((\pi_{\Sigma}\circ \uL)^{-1}(\{z_i\})\right).$$

The formula for the Maslov index follows readily from the excision
principle for the linearized $\DBar$ operator, using the cylindrical
formulation. Consider a sequence of pseudo-holomorphic curves
$\{v_t\}_{t\in\Z}\in\ModFlowL_{J(t)}(\phi)$, where here the subscript
$J(t)$ denotes the almost-complex structure induced on $W(t)$ with
neck-length $t$ as described earlier. Using Gromov's compactness
theorem, after passing to a subsequence, $v_{t}$ converges locally to
a pseudo-holomorphic curve in the symplectic manifold $$W(\infty)
\cong \Big(W_1-\{z_1\}\times[0,1]\times \R\Big) \coprod
\Big(W_2-\{z_2\}\times[0,1]\times \R\Big),$$
which in turn can be
completed to a pseudo-holomorphic curve in $W_1$ and one in $W_2$.
More precisely, we obtain a broken flow-line whose components consist
of pre-glued flowlines and boundary degenerations, finally also nodal
curves supported entirely inside fibers $\Sigma_1\vee\Sigma_2$.  The
representatives in the moduli spaces for $\phi_1$ and $\phi_2$
respectively are gotten by ignoring the matching conditions.

In the case where $\Mas(\phi_1)=1$, the limiting process generically
gives rise rise to an unbroken, preglued flowline, according to the
following dimension counts. Specifically, taking a Gromov
compactification, we obtain a pseudo-holomorphic representative
$\uL_1$ of $\phi_1$, and also a possibly 
broken flow-line representing $\phi_2$.

Assuming that this broken flow-line contains no components which are
closed curves, there is some component of it $u_2$ with the property
that $u_1$ and $u_2$ represents a pre-glued flowline in the sense of
Definition~\ref{def:PreGlued}.  We claim that $u_2$ in fact represents
$\phi_2$.  If $\uL_2$ did not represent $\phi_2$, then it represents a
homology class $\phi'$ which is a component in the Gromov
compactification of $\phi'$.  If this compactification contains
additional boundary degenerations, the Maslov index of $\phi'$ is at
least $2$ smaller than the Maslov index of $\phi_2$ (in view of
Lemma~\ref{lemma:MaslovBoundaryDegenerations}). Moreover, if the
compactification contains other flows, those will serve only to
further decrease the Maslov index of $\phi'$ relative to that of
$\phi$. In sum, in the case where $\phi'$ does not agree with
$\phi_2$, its Maslov index $\Mas(\phi')<2k$. But given $\Delta\in
\Sym^k(\CDisk)$, the moduli space
$$\rho_2^{-1}(\Delta)=\{\uL\in\ModFlowL(\phiL')\big| 
(\pi_\CDisk\circ \uL_2)\left(\pi_\Sigma\circ \uL_2)^{-1}(\{z_2\})\right)=\Delta\} $$
has expected dimension $\Mas(\phi')-2k<0$. Thus, for a generic choice
of $\Delta$ (gotten by $\rho_1(\phi_1)$), this space is empty. 

It is easy to rule out also the case that the Gromov limiting
representing representing $\phi_2$ cannot have any closed
components. To this end, suppose it has some components which
represent the homology class $m[\Sigma_2]$ for some $m>2$. After
deleting these components, we are left with a homology class $\phi'$
with $\Mas(\phi_2')=2(d_2-g_2+1)$
(c.f. Lemma~\ref{lemma:MaslovBoundaryDegenerations}). Some component
$u_2$ of this Gromov compactification has $\rho_2(u_2)=\Delta'$, where
$\Delta'$ is obtained from $\rho_1(\phi_1)$ by deleting $m$ points.
But the moduli space of such points has expected dimension given by
\begin{equation}
\label{eq:eDimModSp}
\Mas(\phi_2')-2(k-m) 
= 2k-2(k-m)-2(d_2+1-g_2)m
=-2m(d_2-g_2)<0,
\end{equation}
and hence it is empty (here, of course, is where
we used our assumption that  $d_2>g_2$). 

Thus, we have established that a sequence of holomorphic
representatives for $\phi$ has a Gromov limit (as we stretch the neck)
to a pre-glued flowline representing $\phi_1$ and $\phi_2$.
Conversely, given a pre-glued flowline, one can obtain a
pseudo-holomorphic curve in $\ModFlowL(\phi)$ by gluing
(cf.~\cite{LipshitzCyl} for further details on this gluing problem,
and~\cite{BEHWZ} for a discussion of gluing in a very general
context).
\qed

\subsection{Counting boundary degenerations}
\label{subsec:BoundaryDegenerations}

We let $\ModDegL(\psi)$ denote the moduli
space of pseudo-holomorphic boundary degenerations
in the homology class of  $\psi$.  Note that
${\mathbb P}{\mathrm SL}_2(\R)$ acts on the this moduli space, and we
let $\UnparModDegL(\psi)$ denote the quotient by this action.
The principles used in proving Theorem~\ref{thm:GenericFiberedProduct}
can also be used to count boundary degenerations.

\begin{theorem}
\label{thm:CountBoundaryDeg}
Consider $\Sigma$, a surface of genus $g$, equipped with a
set of attaching circles $\alpha_1,...,\alpha_{g+\ell-1}$ for a handlebody.
If $\cald(\psi)\geq 0$ and $\Mas(\psi)=2$, then
$\cald(\psi)=A_i$ for some $i$; and indeed in this case
\begin{equation}
  \label{eq:CountBoundaryDegenerations}
  \#\UnparModDegL(\psi)\equiv 
  \left\{\begin{array}{ll}
      0 \pmod{2} & {\text{if $\ell=1$}} \\
      1\pmod{2} & {\text{if $\ell>1$.}}
    \end{array}\right.
\end{equation}
\end{theorem}

\begin{proof}
  It follows readily from
  Lemma~\ref{lemma:MaslovBoundaryDegenerations} that if $\psi$
  is a non-zero homology class of boundary degenerations with 
  $\cald(\psi)\geq 0$, then $\Mas(\psi)\geq 2$, with equality iff
  $\cald(\psi)=A_i$. In the case where $\cald(\phi)=A_i$, it remains
  to verify Equation~\eqref{eq:CountBoundaryDegenerations}.  The case
  where $\ell=1$ has already been established
  in~\cite[Proposition~\ref{HolDisk:thm:GromovInvariant}]{HolDisk} for
  the usual Heegaard-Floer moduli spaces and~\cite{LipshitzCyl} for
  the cylindrical version.
  Consider some region $A_i$, which we re-name simply $A'$. Recall
  that this is a Riemann surface with boundary, equipped with $m+p$
  curves $\alpha_1,...,\alpha_{m+p}$, of which the first $m$ comprise
  its boundary, and the rest are pairwise disjoint, embedded circles in
  the interior. We have fixed also
  $\x\in\alpha_1\times...\times \alpha_{m+p}$.  We reduce to the case where
  $p=0$, by de-stabilizing.  Next, we reduce to the case where the
  number of boundary circles is one. To this end, we can write $A'=A\#
  D$, where $D$ here is a disk with boundary $\alpha_{m}$, and $A$ is
  a planar surface-with-boundary $\alpha_1,...,\alpha_{m-1}$.  Denote
  the connected sum point in $A$ by $z_1$ and the one in $D$ by $z_2$.
  Degenerating the connected sum tube, we obtain a fibered product
  description
  $\ModDegL(A')=\ModDegL(A)\times_{\Sym^1(\CDisk)}\ModDegL(D)$, where the
  fibered product is taken over the maps
  \begin{eqnarray*}
    \rho_1\colon \ModDegL(A)\longrightarrow \Sym^1(\CDisk)
    &{\text{and}}&
    \rho_2\colon \ModDegL(D)\longrightarrow \Sym^1(\CDisk)
  \end{eqnarray*}
  defined as before.
  It is easy to see, though that $\ModDegL(D)$ is smooth, and the 
  map $\rho_2$ is in fact a diffeomorphism. Thus, it follows
  that $\ModDegL(A')\cong \ModDegL(A)$. In this manner, we have reduced
  to the case where $p=0$ and $m=1$, which is, once again, the case where
  $A=D$. Now, since the map $\rho_2$ is ${\mathbb P}{\mathrm SL}_2(\R)$-equivariant, we see
  that $\#\UnparModDeg(D)$ is precisely the degree of the map $\rho_2$,
  which is one.
\end{proof}

\subsection{Notational remark}

The notation of Theorem~\ref{thm:GenericFiberedProduct} suggests using
moduli spaces of pseudo-holomorphic Whitney disks, while its proof
uses cylindrical flow-lines.  One could close the gap by either
appealing to an identification between the two kinds of moduli spaces,
cf.~\cite{LipshitzCyl}, or simply adopting the cylindrical point of
view instead; either approach has the same final outcome. In view of
this remark, we henceforth drop the notational distinction between
cylindrical and more traditional moduli spaces.

\section{Heegaard Floer homology for multi-pointed Heegaard diagrams
revisited}
\label{sec:IdentifyHeegaardFloers}

Having set up the analytical preliminaries,
we now prove the invariance properties promised in
Section~\ref{sec:DefHFL}.

\vskip.2cm

\noindent{\bf{Proof of Lemma~\ref{lemma:DSquaredZero}.}}
In the usual proof that $\partial^2=0$ from Floer homology (cf.
Theorem~\ref{HolDisk:thm:DSquaredZero} for a proof in the original
case, and~\cite{LipshitzCyl} for the cylindrical formulation), this is
observed by counting ends of two-dimensional moduli spaces.

Specifically, fix intersection points $\x$ and $\q$, and a vector
${\mathbf a}=(a_1,...,a_\ell)\in\Z^\ell$. We consider the ends of
$$\coprod_{\{\phi\in\pi_2(\x,\q)\big|\Mas(\phi)=2,
n_\w(\phi)={\mathbf{a}}\}}
\UnparModFlow(\phi).$$
In the case where $\x\neq\q$, these ends cannot contain any boundary
degenerations, since according to
Lemma~\ref{lemma:MaslovBoundaryDegenerations}, these all carry Maslov
index at least $2$ (and hence if they appear in the boundary, the
remaining configuration has Maslov index $\leq 0$, and hence, if it is
non-empty, it must consist of the constant flowline alone). 
Thus, the ends in this case are modeled on
$$\coprod_{\y\in\Ta\cap\Tb}\coprod_{\{\phi_1\in\pi_2(\x,\y),
\phi_2\in\pi_2(\y,\q)\big| \phi_1*\phi_2=\phi\}}
\UnparModFlow(\phi_1)\times \UnparModFlow(\phi_2).$$
And the total count of these ends are given by
\begin{equation}
  \label{eq:CountEnds}
\sum_{\y\in\Ta\cap\Tb} \sum_{\{\phi_1\in\pi_2(\x,\y),
  \phi_2\in\pi_2(\y,\q)\big| n_{\w}(\phi_1)+n_{\w}(\phi_2)={\mathbf a}\}}
(\#\UnparModFlow(\phi_1))\cm (\#\UnparModFlow(\phi_2)),
\end{equation}
which on the one hand must be even, on the other hand, it is
easily seen to be the $U_1^{a_1}\cm...\cm U_\ell^{a_\ell}\cm \q$-component of
$(\partial^{-})^2(\x)$.  In the case where $\q=\x$, there are
additional terms, which count boundary degenerations meeting constant
flowlines, whose total signed count is
$$\sum_{\{\psi\in \pi_2^\alpha(\x)\big| n_{\w}(\psi)={\mathbf a},
  \Mas(\psi)=2\}} \#\UnparModDeg^{\alpha}(\psi) + \sum_{\{\psi\in
  \pi_2^\beta(\x)\big| n_{\w}(\psi)={\mathbf a}, \Mas(\psi)=2\}}
\#\UnparModDeg^{\beta}(\psi).
$$
According to Theorem~\ref{thm:CountBoundaryDeg}, this quantity
vanishes.  More precisely, $\#\UnparModDeg^{\alpha}(\psi)\equiv
0\pmod{2}$ except in the case where $\ell>1$ and
$\cald(\psi)$ is one of the
components $A_i$ of $\Sigma-\alpha_1-...-\alpha_g$, so that
$\mathbf a$ has the form that 
$a_i=0$ for all but one value of $i$, where the component $a_i=1$.
In this case $\#\UnparModDeg^{\alpha}(\psi)\equiv 1\pmod{2}$, but there is also
a unique cancelling $\psi'$ with $\cald(\psi')=B_i$ and
$\#\UnparModDeg^{\alpha}(\psi')\equiv 1\pmod{2}$.
Thus, we are left once again 
with a sum as in Equation~\eqref{eq:CountEnds} which
can be interpreted as the $U_1^{a_1}\cm...\cm U_\ell^{a_\ell}\cm \q$-component of
$(\partial^{-})^2(\x)$. 
\qed

\subsection{Simple stabilization invariance}
\label{subsec:SimpleStabilizations}

The aim of this subsection is to prove that the homology of the chain
complex $\CFm(\Sigma,\alphas,\betas,\ws)$ is invariant under a particularly
simple index
zero/three stabilizations (c.f. Proposition~\ref{prop:HeegaardMoves}).

Specifically, recall that if we have a balanced Heegaard diagram
$(\Sigma,\alphas,\betas,\ws)$ with
$\alphas=\{\alpha_1,...\alpha_{g+\ell-1}\}$,
$\betas=\{\beta_1,...,\beta_{g+\ell-1}\}$, $\ws=\{w_1,...,w_\ell\}$, we can
construct a new balanced Heegaard diagram by introducing a new pair of
separating curves $\alpha_{g+\ell}$ and $\beta_{g+\ell}$ and a new
basepoint $w_{\ell+1}$, so that $\alpha_{g+\ell}$ is isotopic to
$\beta_{g+\ell}$ in  $\Sigma-w_1-...-w_{\ell+1}$.
Let $\alphas'=\{\alpha_1',...,\alpha_{g+\ell-1}',\alpha_{g+\ell}\}$,
$\betas'=\{\beta_1',...,\beta_{g+\ell-1}',\beta_{g+\ell}\}$, where for
$i=1,...,g+\ell-1$, $\alpha_i'$ resp. $\beta_i'$ is obtained from
$\alpha_i$ resp.  $\beta_i$ by a small isotopic translate.
As in Section~\ref{sec:HeegaardDiagrams}, we say that
this new diagram $(\Sigma,\alphas',\betas',\ws\cup \{w_{\ell+1}\})$
is  obtained from $(\Sigma,\alphas,\betas,\ws)$ by an index
zero/three stabilization. 

We call the stabilization {\em simple} if $\alpha_{g+\ell}$
bounds a disk in $\Sigma$ whose closure is disjoint from the other $\alpha_i$.

\begin{figure}
\mbox{\vbox{\epsfbox{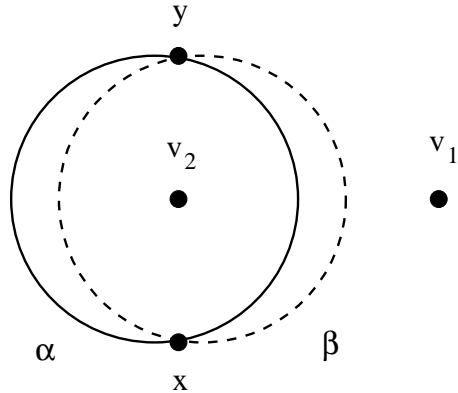}}}
\caption{\label{fig:GenusZeroStabilization}
{\bf{Model doubly-pointed Heegaard diagram for $S^3$.}}  We have
illustrated here the doubly-pointed Heegaard diagram for $S^3$
$(S,\alpha,\beta,\{v_1,v_2\})$ considered above. The picture takes
place on the two-sphere $S$.}
\end{figure}

Let $S$ be the two-sphere, and $\alpha\subset S$ be an embedded curve
which divides $S$ into two regions, each of which contains a basepoint
$v_1$ or $v_2$. Let $\beta$ be a small Hamiltonian isotopic translate
of $\alpha$ (in the complement of $v_1$ and $v_2$), meeting $\alpha$
in two points $x$ and $y$, cf.
Figure~\ref{fig:GenusZeroStabilization}.  Of course,
$(S,\alpha,\beta,\{v_1,v_2\})$ represents a balanced diagram for
$S^3$.

\begin{lemma}
  \label{lemma:SimpleDiagram}
  If $\phi\in\pi_2(x,x)$ or $\pi_2(y,y)$ with $\cald(\phi)\geq 0$,
  then
  \begin{equation}
    \label{eq:MasInequality}
    \Mas(\phi)-2n_{v_i}(\phi)\geq 0.
  \end{equation}
  Moreover, the chain complex $\CFm(S,\alpha,\beta,\{v_1,v_2\})$ is
  given by
  $$
  \begin{CD}
    \Field[U_1,U_2]\cm x @>{U_1-U_2}>>\Field[U_1,U_2]\cm y.
  \end{CD}
  $$
\end{lemma}

\begin{proof}
  Inequality~\eqref{eq:MasInequality} follows easily from
  Lemma~\ref{lemma:MaslovBoundaryDegenerations}.

  For the claim about the chain complex, it is easy to see that boundary
  maps from $y$ to $x$ come in cancelling pairs, whereas there are two
  holomorphic disks from $x$ to $y$ with Maslov index one,
  and they are bigons, one of which contains $v_1$ the other $v_2$.
\end{proof}

\begin{defn}
Let $(\Sigma_1,\xis,\etas,\{u_1,...,u_{\ell}\})$ and $(\Sigma_2,\gammas,\deltas,
\{v_1,...,v_{m}\})$
be a pair of balanced multi-pointed Heegaard diagrams, and choose basepoints
$u_1\in\us$, $v_1\in \vs$. We can form their {\em connected sum}
$(\Sigma,\alphas,\betas,\ws)$, a balanced $\ell+m-1$-pointed Heegaard diagram
whose underlying surface $\Sigma$ is obtained from $\Sigma_1$ and $\Sigma_2$
by forming the connected sum at the points $u_1$ and $v_1$. 
We let $\alphas=\xis\cup\gammas$, $\betas=\etas\cup\deltas$, only now
thought of as curves in $\Sigma$, and $\ws=\{w_1,u_2,...,u_{\ell},
v_2,...,v_m\}$, where $w_1$ is some reference point on the connected sum neck.
We denote  this connected sum by
$(\Sigma_1,\xis,\etas,\{u_1,...,u_{\ell}\})\#(\Sigma_2,\gammas,\deltas,
\{v_1,...,v_{m}\})$. 
\end{defn}

In particular, given an arbitrary multi-pointed the Heegaard diagram
$(\Sigma,\alphas,\betas,\ws)$, a simple index zero/three stabilization
can be thought of as the connected sum with
$(S,\alpha,\beta,\{v_1,v_2\})$ as considered in
Lemma~\ref{lemma:SimpleDiagram}.

Fix a Heegaard multi-diagram $(\Sigma,\alphas,\betas,\ws)$, choose
$\x,\y\in\Ta\cap\Tb$, and $w_1\in \ws$. Correspondingly, let
$\rho^{w_1}\colon \ModFlow(\phi)\longrightarrow \Sym^k(\CDisk)$, where
$k=n_{w_1}(\phi)$, denote the map which assigns to the Whitney disk
$u$ the divisor $u^{-1}(w_1\times\Sym^{g+\ell-2}(\Sigma))$ or, in the
cylindrical formulation, $\pi_{\CDisk}(\pi_{\Sigma}\circ
u)^{-1}(\{w_1\})$. Given $\phi\in\pi_2(\x,\y)$ and $0<t<1$, we let
$M(\phi,t)$ denote the moduli space of pseudo-holomorphic maps $u$
representing $\phi$ with the additional constraint that $(t,0)\in
\rho^{w_1}(u)$.

\begin{lemma}
  \label{lemma:PushOutPoint} Let $\phi\in\pi_2(\x,\y)$ be a homotopy
  class of Whitney disks for a balanced Heegaard diagram
  $(\Sigma,\{\alpha_1,...,\alpha_{g+\ell-1}\},
  \{\beta_1,...,\beta_{g+\ell-1}\},\{w_1,...,w_\ell\})$.  If
  $\Mas(\phi)=2$, then $M(\phi,t)$ is generically a zero-dimension
  moduli space. Moreover, there is a number $\epsilon>0$ with the
  property that for all $t\leq \epsilon$, the only possible non-empty
  moduli spaces $M(\phi,t)$ with $\Mas(\phi)=2$ consist of moduli
  spaces for $\phi\in\pi_2(\x,\y)$ with $\x=\y$, and indeed $\phi$ is
  obtained by splicing a $\Mas(\psi)=2$ boundary degeneration to the
  constant flowline. Indeed, for this moduli space, $$\# M(\phi,t)=
  \left\{\begin{array}{ll} 0 \pmod{2} & {\text{if $\ell=1$}} \\
  1\pmod{2} & {\text{if $\ell>1$.}}  \end{array}\right.$$
\end{lemma}

\begin{proof}
  The dimension statement is clear. Note that the admissibility
  hypothesis ensures that there are at most finitely many
  $\phi\in\pi_2(\x,\y)$ consisting of moduli spaces for with
  $\Mas(\phi)=2$ and $\cald(\phi)\geq 0$, and hence only finitely many
  homotopy classes for which $M(\phi,t)$ could possibly be non-empty
  for some $t$.  Consider one such homotopy class, and suppose that
  $M(\phi,t)$ is non-empty for a sequence of $t\goesto 0$. Let 
  $u_t\in M(\phi,t)$ be a corresponding
  sequence of pseudo-holomorphic curves.  Taking their Gromov limit,
  we obtain a broken flow-line $u$ representing the homology class
  $\phi$.  Since $\rho^{w_1}(u_t)$ contains points arbitrarily close
  to the line $\{0\}\times \R$
  but $w_1$ does not lie in any of the $\betas$,
  we can conclude that the Gromov limit must contain a component which is a
  non-trivial boundary-degeneration $\psi$. According to 
  Lemma~\ref{lemma:MaslovBoundaryDegenerations}, we can conclude that
  $\Mas(\psi)\geq 2$.  Thus, the remaining configuration $\phi-\psi$
  has non-positive Maslov index, and it also has a pseudo-holomorphic
  representative. This forces it to be a constant flowline.
  
  We have thus established that there is a real number $\epsilon>0$
  such that if $M(\phi,t)$ is non-empty for any $0<t\leq \epsilon$,
  then $\phi$ is obtained by splicing a $\Mas(\psi)=2$ boundary
  degeneration to a constant flowline.  The result now follows by
  gluing $\psi$ to the constant flow-line, and applying
  Theorem~\ref{thm:CountBoundaryDeg}.
\end{proof}

We will also need a result about a suitable generalization of $M(\phi,t)$,
but only for the Heegaard diagram
$(S,\alpha,\beta,\{v_1,v_2\})$
introduced above.

Given a divisor $\Delta\in\Sym^k(\CDisk)$, let 
$$M(\phi,\Delta)=\{u\in\ModFlow(\phi)\big| \rho^{v_1}(\phi)=\Delta\}.$$

\begin{lemma}
  \label{lemma:PushOutDivisor}
  Consider the Heegaard diagram $(S,\alpha,\beta,\{v_1,v_2\})$ as
  above, with the two intersection points $x, y\in \alpha\cap\beta$.
  Fix also a generic $\Delta\in \Sym^k(\CDisk)$ for some positive integer $k$. 
  Then,
  $$\sum_{\{\phi\in\pi_2(a,a)|\Mas(\phi)=2k, n_{v_2}(\phi)=0\}} 
  \#M(\phi,\Delta)\equiv 1 \pmod{2},$$
  for $a=x$ or $y$.
\end{lemma}

\begin{proof}
  Consider the case where $a=x$.
  Let $S(x,\Delta)=\sum_{\{\phi\in\pi_2(x,x)\big|\Mas(\phi)=2k, n_{v_2}(\phi)=0\}}\#M(\phi,\Delta)$.
  We claim first that $S(x,\Delta)$ depends on
  $\Delta$ only through its total weight $k$. Specifically, if
  $\Delta_0, \Delta_1$ are two generic points in $\Sym^k(\CDisk)$,
  then $S(x,\Delta_0)=S(x,\Delta_1)$. This follows from the following
  Gromov compactness argument.  Let $\{\Delta_t\}_{t\in [0,1]}$ be a
  path in $\Sym^k(\CDisk)$ connecting $\Delta_0$ and
  $\Delta_1$. Consider the one-dimensional moduli space
  $$\bigcup_{\{t\in[0,1],\phi\in\pi_2(x,x)|\Mas(\phi)=2k\}}
  M(\phi,\Delta_t).$$ This has four types of ends.  The first type
  appear in the expression $$\bigcup_{\left\{ \begin{tiny}
  \begin{array}{c} \phi_1\in\pi_2(x,y) \\ \phi_2\in\pi_2(y,x) \\ t\in
  [0,1] \end{array}\end{tiny} \big| \begin{tiny} \begin{array}{c}
  \Mas(\phi_1)=1 \\ n_{v_2}(\phi_1)=0=n_{v_2}(\phi_2) \\ \Mas(\phi_2)=2k-1 \end{array}
  \end{tiny}\right\}} \UnparModFlow(\phi_1)\times
  M(\phi_2,\Delta_t).$$ The total number of such ends is zero, since
  $x$ is a cycle (in the chain complex for $(S,\alpha,\beta,\{v_1,v_2\})$).  The second type of ends appear in
  $$\bigcup_{\left\{ \begin{tiny} \begin{array}{c} \phi_1\in\pi_2(x,y)
  \\ \phi_2\in\pi_2(y,x) \\ t\in [0,1] \end{array}\end{tiny} \big|
  \begin{tiny} \begin{array}{c} \Mas(\phi_1)=2k-1, \\ \Mas(\phi_2)=1,
  \\ n_{v_2}(\phi_1)=0=n_{v_2}(\phi_2) \end{array} \end{tiny}\right\}}
  M(\phi_1,\Delta_t)\times \UnparModFlow(\phi_2).$$ 
  The total number
  of such ends is zero, since $y$ is a cycle.  Note that $n_{v_1}(\phi_1)=0$
  or $n_{v_2}(\phi_2)=0$ for any such degeneration, since otherwise
  our divisor family $\Delta_t$ remains in a compact portion of the interior
  of the disk for all $t\in[0,1]$.

  The third and fourth
  types of ends appear in $$\bigcup_{\{\phi\in\pi_2(x,x) \big|
  \Mas(\phi)=2k, n_{v_2}(\phi)=0\}}M(\phi,\Delta_0)$$ and the related expression
  $$\bigcup_{\{\phi\in\pi_2(x,x) \big|
  \Mas(\phi)=2k,n_{v_2}(\phi)=0\}}M(\phi,\Delta_1).$$ Thus, taken together, the
  total number of ends is given by
  $S(x,\Delta_0)-S(x,\Delta_1)\pmod{2}$, which must therefore vanish.
    
    Now let $\Delta_T$ consists of $k$ points each of which have
    horizontal component $\leq \epsilon$ (chosen as in
    Lemma~\ref{lemma:PushOutPoint}), and with vertical spacing of at
    least $T$ between them. Taking a limit as $T\goesto\infty$, we see
    that the ends of the parameterized moduli space
    $$\bigcup_{\{T\in
      [1,\infty)\}}\bigcup_{\{\phi\in\pi_2(x,x)\big|\Mas(\phi)=2k\}}
    M(\phi,\Delta_T)$$
    consist of $M(\phi,\Delta_1)$ (the end where $T=1$), and a product
    $\prod_{i=1}^k M(\phi,t_i)$, in the notation of
    Lemma~\ref{lemma:PushOutPoint}. Combining this observation with
    the result of that lemma, we see that $S(x,\Delta_1)\equiv 1\pmod{2}$.

    The case where $a=y$ follows similarly.
\end{proof}

\begin{prop}
  \label{prop:StabilizationInvarianceOne}
  Suppose that $(\Sigma,\alphas',\betas',\ws')$ is obtained from
  $(\Sigma,\alphas,\betas,\ws)$ by a simple index zero/three
  stabilization.  Then, there is an identification of
  $\Field[U_1,...,U_\ell]$-modules
  $$\HFm(\Sigma,\alphas,\betas,\ws) \cong
  \HFm(\Sigma,\alphas',\betas',\ws').$$
  Indeed, multiplication by
  $U_{\ell+1}$ is identified with multiplication by $U_1$.
\end{prop}

\begin{proof}
  We show that for suitable choices of almost-complex structures,
  $\CFm(\Sigma,\alphas',\betas',\ws')$ is identified with the mapping
  cone of a map
  \begin{equation}
    \label{eq:StabilizationMap}
    \begin{CD}
      \CFm(\Sigma,\alphas,\betas,\ws)[U_{\ell+1}]
      @>{U_{\ell+1}-U_1}>>
      \CFm(\Sigma,\alphas,\betas,\ws)[U_{\ell+1}], 
    \end{CD} 
  \end{equation}
  where here
  $\CFm(\Sigma,\alphas,\betas,\ws)[U_{\ell+1}]$ denotes the chain
  complex
  \[\CFm(\Sigma,\alphas,\betas,\ws)\otimes_{\Field[U_1,...,U_\ell]}
  \Field[U_1,...,U_{\ell+1}].\]
  
  We degenerate the Heegaard surface, to realize it as a connected sum
  $$(\Sigma,\alphas,\betas,\ws)\#(S,\alpha,\beta,v_1,v_2)$$ (with
  sufficiently large connected sum neck) where now $w_1$ is identified
  with $v_1$; or more precisely, $w_1$ and $v_1$ are the corresponding
  connected sum points, and we choose a new distinguished point $w_1'$
  to lie in this connected sum region. Note that here $w_{\ell+1}$ 
  is given by $v_2$. We will use a complex structure
  on the connected sum surface with a very long connected sum length
  (and in fact, we will move the connected sum point $v_1$ close to 
  the circle $\beta$ in $S$, as explained below).
  
  Clearly, $\Ta'\cap\Tb'=(\Ta\cap\Tb)\times \{x,y\}$; i.e.
  intersection points come in two types (containing $x$ and $y$), and
  differentials can be of four types. Write $C_x$ resp $C_y\subset
  \CFm(\Sigma',\alphas',\betas',\ws')$ as the submodule generated by
  intersection points containing $x$ resp $y\in \alpha\cap\beta$; i.e.
  we have a module splitting
  $C'=\CFm(\Sigma',\alphas',\betas',\ws',\zs')\cong C_x\oplus C_y$.
  
  We begin by considering the $C_x$-component of the differential of a
  generator in $C_x$. Apply Theorem~\ref{thm:GenericFiberedProduct} to some
  homotopy class $\phi\in\pi_2(\x\times \{x\},\y\times \{x\})$, in the
  case where $\ell>1$ .  Writing $\phi=\phi_1\#\phi_2$, we have that
  $$\Mas(\phi)=\Mas(\phi_1)+\Mas(\phi_2)-2n_{v_1}(\phi).$$
  If
  $\#\ModFlow(\phi)\neq 0$, then $\cald(\phi_2)\geq 0$ (since $\phi_2$
  has a pseudo-holomorphic representative). As an easy consequence of
  Lemma~\ref{lemma:MaslovBoundaryDegenerations}, it follows that
  $\Mas(\phi_2)-2n_{v_1}(\phi)\geq 0$, with equality iff
  $n_{v_2}(\phi_2)=0$. In the case where, $n_{v_2}(\phi_2)\neq 0$ we
  force $\Mas(\phi_1)\leq -1$, and hence it is generically empty.
  Thus, we are left with the case of
  $\Mas(\phi_2)-2n_{v_1}(\phi)=0$, and $n_{v_2}(\phi)=0$.  In this
  case, Theorem~\ref{thm:GenericFiberedProduct} shows that the
  $\y\times x$-component of $\partial^-(\x\times x)$ is given by
  $$\sum_{\{\phi_1\in\pi_2(\x,\y)\big| \Mas(\phi_1)=1\}} \sum_{u_1\in
  \UnparModFlow(\phi_1)}
  \sum_{\phi_2\in\pi_2(x,x)}\#\{u_2\in\ModFlow(\phi_2) \big|
  \rho_1(u_1)=\rho_2(u_2)\}.$$ In this expression, we consider $u_1\in
  \UnparModFlow(\phi_1)$ as an actual map (rather than only one modulo
  re-parameterization) by taking the representatives with the property
  that the projection of $\rho_1(u_1)$ onto the $\R$ factor contains
  $0$, and no positive real number.  According to
  Lemma~\ref{lemma:PushOutDivisor}, then, for each $u_1\in
  \UnparModFlow(\phi_1)$,
  $$\sum_{\phi_2\in\pi_2(x,x)}\#\{u_2\in\ModFlow(\phi_2) \big|
  \rho_1(u_1)=\rho_2(u_2)\}\equiv 1\pmod{2}.$$ Thus, the $\y\times x$
  component of $\partial^-(\x\times x)$ is identified with the $\y$
  component of $\partial^-(\x)$ for the original Heegaard diagram.

  Thus we have identified the $C_x$ component of the differential of
  an element in $C_x\subset C'$ with the differential coming from
  the obvious identification of $C_x$ with the chain complex for the original diagram
  $\CFm(\Sigma,\alphas,\betas,\ws)$.
  In the same manner, the differential within $C_y$ is identified with
  the differential of the original diagram.
  
  We consider now the $C_x$-component of the differential of a
  generator in $C_y$. Again, we expresss $\ModFlow(\phi)$ as a fibered
  product over $\phi_1\in\pi_2(\x,\y)$ with $\phi_2\in\pi_2(y,x)$.
  Now, an application of
  Lemma~\ref{lemma:MaslovBoundaryDegenerations} shows that if
  $\cald(\phi_2)\geq 0$, then $\Mas(\phi_2)-2n_{v_1}(\phi_2)\geq
  1$, with equality iff $n_{v_2}(\phi_2)=0$. In
  the case of equality, we have that $\Mas(\phi_1)=0$, and hence it
  must be constant. This forces $\phi_2$ to be one of the two flows
  from $y$ to $x$ with $\Mas(\phi_2)=1$ and $\cald(\phi_2)\geq 0$.
  Each of these homotopy classes admits a unique holomorphic
  representative, and hence the differential cancels: the $C_x$
  component of the boundary of something in $C_y$ is trivial.
  
  Finally, we consider the $C_y$ component of the differential of a
  generator in $C_x$. Splitting a homotopy class $\phi$ in the fibered
  sum description, we once again have
  $$\Mas(\phi_1)+\Mas(\phi_2)-2n_{v_1}(\phi_1)=1,$$
  and hence the
  condition that $\Mas(\phi_1)\geq 0$ (which is needed for its
  corresponding moduli space to be non-empty) translates into the
  condition that $\Mas(\phi_2)-2n_{v_1}(\phi_2)\leq 1$.  Moreover,
  for $\phi_2\in\pi_2(x,y)$ with $\cald(\phi_2)\geq 0$, we have that
  $\Mas(\phi_2)-2n_{v_1}(\phi_2)\geq -1$. Parity considerations exclude
  the possibility that the quantity equals $0$.
  
  In the case where $\Mas(\phi_2)-2n_{v_1}(\phi_2)=1$, we conclude
  that $\Mas(\phi_1)=0$, and hence it must be constant, forcing
  $n_{v_1}(\phi_2)=0$. Thus, $\phi_2\in \pi_2(x,y)$ must be the unique
  homotopy class with $\Mas(\phi_2)=1$, $\cald(\phi_2)\geq 0$, 
  $n_{v_1}(\phi_2)=0$, and $n_{v_2}(\phi_2)=1$.
  Thus, the corresponding component of 
  $\partial (\x\times x)$ is $U_{\ell+1}\cm \x\times y$
  (note that in the connected sum diagram, the reference point $v_2$
  corresponds to the new variable $U_{\ell+1}$).
  
  In the remaining cases, 
  $\Mas(\phi_1)=2$ and
  $\Mas(\phi_2)-2n_{v_1}(\phi_2)=-1$.  Since $\Mas(\phi_2)\geq 0$, we
  conclude that $n_{v_1}(\phi_2)>0$, while the condition that
  $\cald(\phi_2)\geq 0$ and $\Mas(\phi_2)-2n_{v_1}(\phi_2)=-1$ readily
  forces $n_{v_2}(\phi_2)=0$.
  
  Suppose that $n_{v_1}(\phi_2)=1=\Mas(\phi_2)$. In this case,
  $\phi_2$ is constrained to be a homotopy class with a unique
  holomorphic representative up to translation.  Let $u_2$ be a
  holomorphic representative of $\phi_2$, and
  $\Delta_2=u_2^{-1}(\{v_1\})$. Note that this consists of a single
  point up to translation in $\CDisk$; after suitable translation, we
  arrange for $\Delta_2=(t,0)\in\CDisk$.  We wish now to apply
  Theorem~\ref{thm:GenericFiberedProduct}, with now $\Sigma$ playing
  the role of $\Sigma_2$ in the statement of that theorem. For this,
  we need to assume that $\ell$, the number of marked points $\ws$, is
  greater than one (so that we are taking more than the $g^{th}$
  symmetric product of $\Sigma$, as required in the hypothesis of the
  last part of Theorem~\ref{thm:GenericFiberedProduct}). We return to
  the case where $\ell=1$ at the end of this proof.
  
  \vskip.2cm {\noindent{\bf{Completion of the proof of the proposition
        when $\ell>1$}.}}  According to
  Theorem~\ref{thm:GenericFiberedProduct}, count of points in all
  Maslov index one moduli spaces $\phi_1\#\phi_2$, where
  $\phi_2\in\pi_2(x,y)$ has $n_{v_1}(\phi_2)=\Mas(\phi_2)=1$ is given
  by the map $$\delta_1\colon
  C_x\longrightarrow
  C_y$$
  defined by
  $$\delta_1(\x)=\sum_{\y\in\Ta\cap\Tb}\sum_{\{\phi_1\in\pi_2(\x,\y)\big|
  \Mas(\phi_1)=2, n_{w_1}(\phi_1)=1\}} \Big(\#M(\phi_1,t)\Big)\cm
  \left(\prod_{i=1}^{\ell} U_1^{n_{w_i}(\phi_1)}\right)\cm \y,$$ in
  the notation of Lemma~\ref{lemma:PushOutPoint}.  By choosing $v_1$
  sufficiently close to $\beta$, we can arrange for $t$ to be
  arbitrarily close to $0$. According to
  Lemma~\ref{lemma:PushOutPoint} for suitable choice of $t$, this
  count is given by \begin{equation} \label{eq:MovePointPart}
  \delta_1(\x)=U_1\cm \x \end{equation}

  There are in principle other terms which count homotopy classes with
  $n_{v_1}(\phi_2)>1$ (and corresponding to factorizations of $\phi$
  as $\phi_1$ and $\phi_2$ with $k=n_{w_1}(\phi_1)>1$,
  $\Mas(\phi_1)=2$).  Our claim is that for a sufficiently large
  parameter $\tau$, these homotopy classes have trivial contribution.
  Here, $\tau$ parameterizes the choice of connected sum point
  $v_1(\tau)$ in $S$, with the limit
  $\lim_{\tau\goesto\infty}v_1(\tau)$ given by a point $v_1^{\infty}$
  on the curve
  $\beta$. Note that we have already taken $\tau$ large to apply
  Lemma~\ref{lemma:PushOutPoint} in establishing
  Equation~\eqref{eq:MovePointPart}.
  
  Suppose now that for a sequence of $\tau$ going to infinity, the
  moduli spaces $\ModFlow_{\tau}(\phi)$ are non-empty for all choice
  of connected sum neck length. Then, for all sufficiently large
  $\tau$, the fibered product $\ModFlow(\phi_1)\times_{\Sym^k(\CDisk)}
  \ModFlow(\phi_2)$ is non-empty, where here $k=n_{v_1}(\phi_2)$.
  
  Thus, we obtain a sequence of pseudo-holomorphic representatives
  $u_1^{\tau}\times u_2^{\tau}$ of the fibered product.  Clearly,
  there are Gromov limits ${\overline u}_1^\infty$ and ${\overline
  u}_2^\infty$ as $\tau\goesto\infty$ of the curves $u_1^{\tau}$ and
  $u_2^\tau$. By dimension counts, since $\Mas(u_1^{\tau})=2$, there
  are only three possible types of limit for ${\overline u}_1^\infty$:
  either it is a strong limit to a pseudo-holomorphic disk, or it is a
  weak limit to a singly-broken flowline, or it contains a boundary
  degeneration, in which case the remaining component must be a
  constant flowline.  In this latter case, $k=1$, and hence it has
  been covered earlier.

  Suppose that the limit ${\overline u}_1^\infty=u_1^{\infty}$ is not
  a broken flowline, and let $u_2^{\infty}$ denote the matching
  component of ${\overline u}_2^{\infty}$ (i.e. ${\overline
  u}_2^\infty$ could {\em a priori} be a broken flow-line, but it has
  some component $u_2^{\infty}$ with the property that
  $\rho_1(u_1^{\infty})=\rho_2(u_2^{\infty})$ .  However, since
  $v_1(\tau)$ limits to $\beta$, we see that $\rho_2(u_2^{\infty})$
  contains some points on the $\beta$-boundary. To achieve this, we
  must have a sequence $u_1^\tau$ with arbitrarily large $\tau$, with
  $\rho_1(u_1^{\tau})$ containing points arbitrarily close to the line
  $\{0\}\times \R$. According to Lemma~\ref{lemma:PushOutPoint}, this
  forces $k=1$, a case considered earlier.
  
  In the remaining case, the Gromov limit ${\overline u}_1^\infty$ is
  given as a broken flow-line ${\overline u}_1^\infty=a\circ
  b$. Again, by simple dimension counts, we see that
  $\Mas(a)=\Mas(b)=1$. There is also a corresponding 
  ${\overline u}_2^\infty=a'\circ b'$. Since $\Mas(\phi_2)$ is odd, we
  can conclude that so is either $\Mas(a')$ or $\Mas(b')$. Suppose it
  is $\Mas(a')$ which is odd. It is easy to see that for
  $(a')^{-1}(v_1^{\infty})$ contains points on the boundary
  $\{0\}\times \R$. Moreover, the same reasoning as before, with $a$
  and $a'$ playing the roles of $u_1^\infty$ and $u_2^\infty$
  respectively, shows that $a$ is supported in a homotopy class which
  admits holomorphic representatives $a^\tau$ with $\rho_1(a^\tau)$
  containing points arbitrarily close to $\{0\}\times \R$.  But this
  is impossible, as $w_1$ is disjoint from $\betas$, and $a^{\tau}$
  has to be one of finitely many holomorphic disks up to
  translation. The case where $\Mas(b')$ is odd follows {\em mutas
  mutandis}.

  Putting the above facts together, we obtain the desired
  identification of $\CFm(\Sigma,\alphas',\betas',\ws')$ with the
  mapping cone of Equation~\eqref{eq:StabilizationMap}, giving an
  expression $$\HFm(\Sigma,\alphas',\betas',\ws')
  \cong\frac{\HFm(\Sigma,\alphas,\betas,\ws)[U_{\ell+1}]}{U_{\ell+1}-U_1},$$
  at least in the case where $\ell>1$.
  
  \vskip.2cm {\noindent{\bf{Proof of the proposition when $\ell=1$}.}}
  In this case, Theorem~\ref{thm:GenericFiberedProduct} cannot be
  applied directly to identify the $C_y$ component of $C_x$. 
  
  Specifically, consider a homology class of $\phi=\phi_1\#\phi_2$,
  with $\Mas(\phi)=1$, $\phi_2\in\pi_2(x,y)$, and
  $n_{v_1}(\phi_2)=1=\Mas(\phi_2)$.  Take a Gromov limits of elements
  of $\ModFlow(\phi)$, as the connected sum is degenerated.  Again,
  this limits to the unique (up to translation) holomorphic
  representative $u_2$ of $\phi_2$, but now we cannot assume that the
  Gromov limit ${\overline u}_1$ from the other side is a flow-line.
  More specifically, the dimension counts (cf.
  Equation~\eqref{eq:eDimModSp}) which ruled out the possibility of a
  closed component in the Gromov compactification no longer apply.
  Indeed, by pushing the connected sum point $v_1$ sufficiently close
  to $\beta$, we can use Lemma~\ref{lemma:PushOutPoint} to rule out
  the case where the component of $u_1$ in the Gromov compactification
  which matches with $u_2$ is an actual cylindrical flow-line.
  Rather, it must be a closed holomorphic curve representing the
  homology class $[\Sigma_1]$ (with multiplicity one).  The argument
  from stabilization invariance as in~\cite{LipshitzCyl} applies now to
  show that $\UnparModFlow(\phi)=1$.  More precisely, the Gromov limit
  $u_1$ equals a constant flow-line meeting a copy of $\Sigma_1$, with
  $\rho_2(u_2)=p$.  In this case, gluing can be used to show that
  $\#\ModFlow(\phi)=C\cm \#\ModFlow(\phi_2)$, where here $C$ is the
  count of representatives $\CDisk\# \Sigma_1$ with a marked point $q$
  to $\Sigma_1\times [0,1]\times \R$ in the homology class
  $[\Sigma_1]$ which maps $u(p)$ to $(w_1,p)$. The fact that $C=1$, in
  the case where $g_1=1$ is calculated in the proof of stabilization
  invariance of cylindrical Heegaard Floer homology (cf. Appendix B
  of~\cite{LipshitzCyl}). One can show that $C=1$ for arbitrary $g$
  follows from this case by, for example, by realizing $\Sigma_1$ as a
  connected sum of $g_1$ copies of a genus one surface, and
  degenerating along all the necks.  Now, $\ModFlow(\phi)=1$ follows
  from the elementary calculation that $\#\ModFlow(\phi_2)=1$. 
  
  Thus, we have computed that the counts of all Maslov index one
  moduli spaces of the form $\phi_1\#\phi_2$, where $\phi_2(x,y)$ is
  the flowline with $n_{v_1}(\phi_2)=\Mas(\phi_2)=1$ is given by the
  same map $\delta_1$ as in Equation~\eqref{eq:MovePointPart}.  We can
  rule out the case where $\Mas(\phi_1)=2$ and $k>1$ much as before.
  Specifically, by pushing the connected sum point close to $\beta$,
  Lemma~\ref{lemma:PushOutPoint} shows that the corresponding
  Gromov limit $u_1$ must contain at least one closed curve component.
  Removing this component, we are left with a broken flow-line with
  Maslov index zero, and hence, a constant flowline, hence showing
  that we needed to be in the case where $k=1$. The argument is now completed
  as before.
\end{proof}

\subsection{Model calculations}

The aim of the present subsection is to continue with the methods from
the previous subsection to perform some model calculations which will
be useful for establishing handleslide invariance.

\begin{defn}
  Let $C$ be a chain complex over the polynomial algebra
  $\Field[U_1,...U_\ell]$. We say that it is of 
  $\Field[U]$-type if for $i\neq j$, 
  there are chain homotopies $U_i\simeq U_j$
  (thought of as endomorphisms of $C$).
\end{defn}

Of course, if a chain complex $C$ is of $\Field[U]$-type, then its
homology is a module over the polynomial algebra $\Field[U]$, where
$U$ acts by multiplication by any of the $U_i$.

Let $(\Sigma,\alphas,\betas,\ws)$ be a Heegaard diagram, and fix two
preferred basepoints $w_1, w_2\in \ws$. Let $\Sigma'$ be a surface
obtained by attaching a one-handle in the neighborhood of $w_1$ and
$w_2$, and fix a pair of circles $\alpha_0$ and $\beta_0$ which are
supported inside the handle, each a small isotopic translate of one
another, separating $w_1$ and $w_2$. Equivalently, we
form a double-connected sum of $\Sigma$ along $w_1$ and $w_2$
with the sphere $S$ as in Figure~\ref{fig:GenusZeroStabilization}.

\begin{prop}
  \label{prop:StabilizationInvarianceTwo}
  Let $(\Sigma,\alphas,\betas,\ws)$ be a Heegaard diagram, and let
  $(\Sigma',\alphas',\betas',\ws')$ be the Heegaard diagram obtained
  by attaching a one-handle in the above sense, so that if the
  original describes a three-manifold $Y$, then the second diagram
  describes $Y'=Y\# (S^2\times S^1)$.  Suppose that
  $\HFm(\Sigma,\alphas,\betas,\ws)$ of $\Field[U]$-type, then the same
  is true of $\HFm(\Sigma',\alphas',\betas',\ws')$; and indeed
  $\CFm(\Sigma',\alphas',\betas',\ws')\simeq \CFm(\Sigma,\alphas,\betas,\ws)\oplus
        \CFm(\Sigma,\alphas,\betas,\ws)$.
\end{prop}

\begin{proof}
  We analyze as in Proposition~\ref{prop:StabilizationInvarianceOne},
  only this time stretching near both $w_1$ and $w_2$.
  In this case, homotopy classes $\phi$ break as fibered products of
  homotopy classes $\phi_1$ for $(\Sigma,\alphas,\betas,\ws)$, and 
  $\phi_2$ for $(S,\alpha,\beta,v_1,v_2)$. Now we have
  $$\Mas(\phi)=\Mas(\phi_1)+\Mas(\phi_2)-2n_{v_1}(\phi_2)-2n_{v_2}(\phi_2).$$
  
  The same dimension counts as in the proof of
  Proposition~\ref{prop:StabilizationInvarianceOne} express
  $\CFm(\Sigma',\alphas',\betas',\ws')$ as a mapping cone of a map
  \begin{equation} \label{eq:StabilizationMapTwo} \begin{CD}
  \CFm(\Sigma,\alphas,\betas,\ws) @>{U_1-U_2+\delta}>>
  \CFm(\Sigma,\alphas,\betas,\ws).  \end{CD} \end{equation} As in the
  proof of that proposition, the first chain complex is identified
  with $C_x$ and the second with $C_y$.  The map $\delta$ counts those
  disks which have $n_{v_1}(\phi_2)+n_{v_2}(\phi_2)\geq 2$. By moving
  the connected sum points as in
  Proposition~\ref{prop:StabilizationInvarianceOne}, we can arrange
  that these terms contribute trivially.

  Thus, since $U_1\simeq U_2$, 
  the chain map from Equation~\eqref{eq:StabilizationMapTwo} is
  null-homotopic, and it follows at once that 
  $\CFm(\Sigma',\alphas,\betas',\ws')$ is chain homotopic
  to the direct sum of two copies of $\CFm(\Sigma,\alphas,\betas,\ws)$.
\end{proof}

Let $\Sigma$ be an oriented two-manifold, and $\alphas$ be a
collection of attaching circles, and let $\ws$ be a collection of
basepoints, one in each component of
$\Sigma-\alpha_1-...-\alpha_{g+\ell-1}$.  Let $\alphas'$ be a set of
attaching circles obtained as small exact Hamiltonian translates of
$\alpha_i$.  Specifically $\alpha_i'$ is obtained as an exact
Hamiltonian translate of $\alpha_i$, so that
$\alpha_i\cap\alpha_j'=\emptyset$ unless $i=j$, in which case the
intersection consists of two points of transverse intersection.
Moreover, the Hamiltonian isotopy never crosses any of the $w_i\in
\ws$. Clearly, $(\Sigma,\alphas,\alphas',\ws)$ represents
$\#^{g}(S^2\times S^1)$, where $g$ is the genus of $\Sigma$.

\begin{prop}
  \label{prop:SmallPerturbation}
  $\CFm(\Sigma,\alphas,\alphas',\ws)$ is of $\Field[U]$-type, and its
  homology is isomorphic to $\Field[U]\otimes \Wedge^* V$, where $V$
  is a $g$-dimensional vector space.
\end{prop}

\begin{proof}
  We can reduce to the case where $g=0$ by repeatedly applying
  Proposition~\ref{prop:StabilizationInvarianceTwo}.
  In the case where $g=0$, the proposition is proved after repeated
  applications of
  Proposition~\ref{prop:StabilizationInvarianceOne}. 
\end{proof}

As an example, suppose that $(\Sigma,\alphas,\betas,\ws)$ is an
admissible multi-pointed Heegaard diagram, and suppose that $\betas'$
is obtained from $\betas$ by a small perturbation as in
Proposition~\ref{prop:SmallPerturbation}. According to
Proposition~\ref{prop:SmallPerturbation}, there is an element
$[\Theta_{\beta\beta'}]\in\HFm(\Sigma,\betas,\betas',\ws)$ of maximal degree.

We can define a corresponding map
$$\Phi_{\alpha\beta\beta'}\colon \CFm(\Sigma,\alphas,\betas,\ws)
\longrightarrow \CFm(\Sigma,\alphas,\betas',\ws)$$
by 
$$\Phi_{\alpha\beta\beta'}(\xi)
=\sum_{\y\in\Ta\cap\Tb}\sum_{\{\psi\in\pi_2(\x,\Theta_{\beta\beta'},\y)\}}
\#\ModFlow(\psi) \left(\prod_{i=1}^\ell U_i^{n_{w_i}(\psi)} \right)\cm \y.$$
This map is a chain map (compare~\cite[Section~\ref{HolDisk:sec:HolTriangles}]{HolDisk}). (In this notation, we implicitly assume that $[\Theta_{\beta\beta'}]$
is represented by a single intersection point $\Theta_{\beta\beta'}$; more
generally, our map is gotten by summing triangle maps over the various
intersection points whose sum represents the homology class.)

\begin{prop}
  \label{prop:SmallPerturbationInvariance}
  If the curves in $\betas'$ are sufficiently close to those in
  $\betas$, chosen so that each $\beta_i'$ meets $\beta_i$ in precisely two
  intersection points, then the
  map $\Phi_{\alpha\beta\beta'}$ defined above induces an isomorphism
  in homology.
\end{prop}

\begin{proof}
  This follows as
  in~\cite[Proposition~\ref{HolDisk:prop:Isomorphism}]{HolDisk}. The
  point is that for each $\x\in\Ta\cap\Tb$, there is a corresponding
  nearest point $\x'\in\Ta\cap{\mathbb T}_{\beta'}$, and also a
  corresponding small triangle
  $\psi\in\pi_2(\x,\Theta_{\beta\beta'},\y)$.  Thus, the map obtained
  by counting only these smallest triangles induces an isomorphism of
  chain groups.  Using the energy filtration, it follows that
  $\Phi_{\alpha\beta\beta'}$ is an isomorphism of chain complexes.
\end{proof}

\subsection{Handleslide invariance}

We can now adapt the proof of handleslide invariance of Heegaard Floer
homology as in
~\cite[Section~\ref{HolDisk:sec:HandleSlides}]{HolDisk} to
establish handleslide invariance of $\HFm$ in the present context.

Specifically, start with a Heegaard diagram
$(\Sigma,\alphas,\betas,\ws)$.  Let $\gammas$ be obtained from
$\betas$ by a single handleslide. 

\begin{prop}
  \label{prop:HandleSlides}
  There is an identification
  $\HFm(\Sigma,\alphas,\betas,\ws)\cong\HFm(\Sigma,\alphas,\gammas,\ws)$.
\end{prop}

\begin{proof}
  As a first step, we claim that $\HFm(\Sigma,\betas,\gammas,\ws)\cong
  \Field[U]\otimes \Wedge^* V$, where $V$ is a $g$-dimensional
  $\Field$-vector space. To see this, de-stabilize using
  Propositions~\ref{prop:StabilizationInvarianceOne} and
  \ref{prop:StabilizationInvarianceTwo} to the case where there are
  exactly two $\beta$ curves and two $\gamma$ curves. There are three
  cases, according to whether $g=g(\Sigma)$ is $0$, $1$, or $2$. The
  case where $g=2$ is established in~\cite{HolDisk}; the other two
  cases are easily established by the same calculation.
  
  With this said, there is a canonical top-dimensional generator
  $\Theta_{\beta\gamma}$ of $\HFm(\Sigma,\betas,\gammas,\ws)$. The
  handleslide map
  $$
  \Psi_{\alpha\beta\gamma}\colon
  \HFm(\Sigma,\alphas,\betas,\ws)\longrightarrow
  \HFm(\Sigma,\alphas,\gammas,\ws)
  $$
  is defined by counting
  pseudo-holomorphic triangles: $\Psi_{\alpha\beta\gamma}=
  f_{\alpha\beta\gamma}(\cdot\otimes\Theta_{\beta\gamma})$; i.e. here

  $$f_{\alpha\beta\gamma}(\x\otimes\Theta_{\beta\gamma})=\sum_{\y\in\Ta\cap\Tc}\sum_{\{\psi\in\pi_2(\x,\Theta_{\beta\gamma},\y)\big|\Mas(\psi)=0\}}
  \#\ModFlow(\psi) \cm \left(\prod_{i=1}^{\ell}
    U_i^{n_{w_i}(\psi)}\right)\cm \y.$$
  
  Now, according to associativity of the triangle maps, we see that
  the composite
  $$\Psi_{\alpha\gamma\beta'}\circ\Psi_{\alpha\beta\gamma}
  =f_{\alpha\beta\beta'}(\cdot\otimes f_{\beta\gamma\beta'}(
  \Theta_{\beta\gamma}\otimes\Theta_{\gamma\beta'})).$$
  Next, we
  verify that $f_{\beta\gamma\beta'}(
  \Theta_{\beta\gamma}\otimes\Theta_{\gamma\beta'})=\Theta_{\beta\beta'}$
  is the canonical top-dimensional generator of
  $\HFm(\Sigma,\betas,\betas',\ws)$, which is calculated in
  Proposition~\ref{prop:SmallPerturbation}; i.e.  we obtain the map
  $\Phi_{\alpha\beta\beta'}$ studied in
  Proposition~\ref{prop:SmallPerturbationInvariance}.
  The map is an isomorphism now according to that proposition.
\end{proof}

\subsection{Invariance}
\vskip.2cm We prove the invariance of $\HFm$ as introduced in
Section~\ref{sec:DefHFL}.

\noindent{\bf Proof of Theorem~\ref{thm:InvarianceHFm}.}
We show that $\HFm$ is invariant under the four Heegaard moves
from Proposition~\ref{prop:HeegaardMoves}, or, more specifically, 
their admissible versions as in Proposition~\ref{prop:Admissibility}.

Isotopy invariance follows exactly as in~\cite[Section~\ref{HolDisk:sec:Isotopies}]{HolDisk}.

Handleslide invariance was established in Proposition~\ref{prop:HandleSlides}.

Invariance under index one and two stabilizations follows from exactly
as in~\cite[Section~\ref{HolDisk:sec:Stabilization}]{HolDisk}.

Finally, it remains to consider invariance index zero/three
stabilizations.  In this case, we introduce a new pair of (isotopic)
curves $\alpha_{g+\ell}$ and $\beta_{g+\ell}$, each of which
separates $\Sigma$ in two. In the case where the stabilization is
simple in the sense of Subsection~\ref{subsec:SimpleStabilizations},
i.e. $\alpha_{g+\ell}$ bounds a disk containing $w_{\ell+1}$ and
none of the other $\alpha_i$ (or $w_i$), invariance is established in
Proposition~\ref{prop:StabilizationInvarianceOne}.  It is easy to see that
general index
zero/three stabilizations can be achieved by simple index zero/three
stabilizations, followed by a sequence of handleslides.  \qed

The case of $\CFa$ can be handled more quickly:

\vskip.2cm
\noindent{\bf Proof of Theorem~\ref{thm:InvarianceHFaa}.}
Let $(\Sigma,\alphas,\betas,\ws)$ be the multi-pointed Heegaard
diagram, with $\ws=\{w_1,...,w_\ell\}$. Let $P$ be a planar suface
with $\ell$ boundary components. Form the surface
$\Sigma'=(\Sigma-\nbd{\ws})\cup P$.    Consider
$(\Sigma',\alphas,\betas,w)$, where here $w$ is chosen in the region
$P\subset \Sigma'$. It is easy to see that
$(\Sigma',\alphas,\betas,w)$ is a pointed Heegaard diagram for
$Y\#(\#^{\ell-1}(S^2\times S^1))$. Admissibility of the original
diagram corresponds precisely to admissibility for this new diagram.
Moreover, the intersection points of $\Ta\cap\Tb$ in
$\Sym^{g+\ell-1}(\Sigma)$ are identical with those in
$\Sym^{g+\ell-1}(\Sigma')$. Finally, holomorphic curves in
$\Sym^{g+\ell-1}(\Sigma)$ which are disjoint from all the $w_i$ are
identical with those in $\Sym^{g+\ell-1}(\Sigma')$ which are disjoint
from $w$. The result now follows.
\qed

\section{Invariance of link invariants}
\label{sec:LinkInvariants}

Invariance of the filtered chain homotopy type of the link filtration
$\CFLm(S^3,\orL)$ is an easy consequence of the methods from
Section~\ref{sec:IdentifyHeegaardFloers}. 

\vskip.2cm
\noindent{\bf Proof of Theorem~\ref{thm:InvarianceHFLFilt}.}
We must verify that $\CFLm(\Sigma,\alphas,\betas,\ws,\zs)$ is invariant
under isotopies and handleslides supported in the complement of $\ws$
and $\zs$, and also index one/two stabilizations.

If $(\Sigma',\alphas',\betas',\ws',\zs')$ is obtained from an index
one/two stabilization, then the methods from
Section~\ref{HolDisk:sec:Stabilization} of~\cite{HolDisk} actually
give an isomorphism of chain complexes
$\CFLm(\Sigma',\alphas',\betas',\ws',\zs')\cong
\CFLm(\Sigma,\alphas,\betas,\ws,\zs)$.
Performing the stabilization away from $\ws$ and $\zs$, we see that 
the assignment to relative $\SpinC$ structures is unaffected.

Handleslides amongst the $\betas$ are understood as follows. As before,
if $\gammas$ is obtained from $\betas$ by a handleslide, then there is a 
chain map (c.f. Proposition~\ref{prop:HandleSlides})
$$\Psi_{\alpha\beta\gamma}\colon
\CFLm(\Sigma,\alphas,\betas,\ws)\longrightarrow
\CFLm(\Sigma,\alphas,\gammas,\ws),$$
defined by counting holomorphic
triangles. Since all triply-periodic domains are disjoint from $\ws$
and $\zs$, it follows easily that if $\x\leq \y$, then
$\Psi_{\alpha,\beta,\gamma}(\x)\leq \Psi_{\alpha,\beta,\gamma}(\y)$.
We thus consider the induced map on the associated graded object, 
i.e. the map ${\widehat \Psi}_{\alpha\beta\gamma}$ induced by counting pseudo-holomorphic
triangles which are disjoint from $\ws$ and $\zs$ alike. Note that 
$$\spincrel_{\ws,\zs}(\x)=\spincrel'_{\ws,\zs}({\widehat \Psi}_{\alpha\beta\gamma}(\x)),$$
where here $\spincrel'_{\ws,\zs}$ is the map associated to the Heegaard diagram
$(\Sigma,\alphas,\gammas,\ws,\zs)$, since the triangles used in 
the definition of ${\widehat\Psi}$ this count
are disjoint from $\ws$ and $\zs$. Now, 
adapting the
argument from Proposition~\ref{prop:HandleSlides}, we have that the map
$${\widehat\Psi}_{\alpha\gamma\beta'}\circ{\widehat\Psi}_{\alpha\beta\gamma}
={\widehat f}_{\alpha\beta\beta'}(\cdot\otimes \Theta_{\beta\beta'})$$
is an isomorphism. Specifically, here
$${\widehat f}_{\alpha\beta\beta'}(\x\otimes \Theta_{\beta\beta'})
=\sum_{\y\in\Ta\cap\Tb}\sum_{\{\psi\in\pi_2(\x,\Theta_{\beta\beta'},\y)\big|\Mas(\psi)=0,
  n_\w(\psi)=n_\z(\psi)=0\}}\#\ModFlow(\psi) \y.$$
The fact that this
is an isomorphism follows readily as in
Proposition~\ref{prop:SmallPerturbationInvariance}, with the
observation that the small triangles considered there do not contain
any of the the basepoints $\{w_i,z_i\}_{i=1}^\ell$. Since
${\widehat\Psi}$ induces an isomorphism of $\RelSpinC(Y,L)$-graded
complexes, it follows formally that $\Psi_{\alpha\beta\gamma}$ induces an
isomorphism of $\RelSpinC(Y,L)$-filtered complexes.  \qed

\subsection{Link homology as an associated graded object}
We consider now the link filtration.

\vskip.2cm
\noindent{\bf{Proof of Theorem~\ref{thm:FilterHFLa}.}}
From its construction, we see that $\HFLa$ is the homology of the
graded object associated to a $\Z^\ell$-filtration of
$\CFa(\Sigma,\alphas,\betas,\ws)$, which, according to
Theorem~\ref{thm:InvarianceHFaa}, calculates
$\HFa(\#^{\ell-1}(S^2\times S^1))$.  
But $\HFa(\#^{\ell-1}(S^2\times S^1))$
is identified with the exterior algebra of $H_1(\#^{\ell-1}(S^2\times S^1);\Field)$ (cf. for example~\cite[Lemma~\ref{HolDisk:lemma:Tori1}]{HolDisk}).
Invariance of the
spectral sequence is a consequence of
Theorem~\ref{thm:InvarianceHFLFilt}.  \qed

\vskip.2cm
\noindent{\bf{Proof of Theorem~\ref{thm:HFLaHFLm}.}}
Fix $\relspinc\in\RelSpinC(Y,L)$.  Recall that the chain complex
computing $\HFLm(Y,\orL,\spincrel)$ consists of symbols
$U^{a_1}\cm...\cm U^{a_\ell}\otimes \x$, 
where $a_1,...,a_\ell$ are non-negative
integers and $\x\in\Ta\cap\Tb$, satisfying the constraint that
$$
\relspinc_{\ws,\zs}(\x)-a_1\PD[\mu_1]-...-a_\ell\PD[\mu_\ell]=\relspinc.
$$
This set is endowed with the differential as in
Definition~\ref{eq:DefD}, except that now we sum over only those
$\phi\in\pi_2(\x,\y)$ with $n_{\z}(\phi)=0$. The further filtration on
$\CFLa(Y,\orL,\relspinc)$ is given by the map $U^{a_1}\cm...\cm
U^{a_\ell}\otimes \x\mapsto (a_1,...,a_\ell)$.  Its associated graded
object is clearly the one stated in the theorem.  \qed

\vskip.2cm
\noindent{\bf{Proof of Theorem~\ref{thm:FilterHFLm}.}}
The spectral sequence arizes from the filtration of 
$\CFLm(S^3,\orL)$. 
The fact that the homology of this total complex is $\Z[U]$
follows from the fact that the underlying chain complex 
associated to $(\Sigma,\alphas,\betas,\ws)$ has homology $\HFm(S^3)$
(according to Theorem~\ref{thm:InvarianceHFm}),
which in turn is known to be $\Z[U]$ (cf.~\cite{HolDisk}).
\qed

\subsection{Forgetful functors}
\label{subsec:Forgetfuls} 

Let $\orL$ be an $\ell$-component oriented link, which we assume to be
in $S^3$ for simplicity. Distinguish a component $K_1\subset L$, and
consider the filling map
$$\Fill{K_1}\colon\RelSpinC(Y,L)\longrightarrow \RelSpinC(Y,L-K_1)$$
studied in Subsection~\ref{subsec:Filling}.
Under this map, a $\RelSpinC(Y,L)$-filtered chain complex
can be viewed as a $\RelSpinC(Y,L-K_1)$-filtered chain complex.

Let $M$ be the rank two graded vector space with one generator in
grading $0$ and another in grading $-1$. 

\begin{prop}
\label{prop:Forgetfuls}
Let $\orL$ be an oriented, $\ell$-component link in $S^3$, and
distinguish the first component $K_1$. Consider the filtration
$\CFLa(\orL)$, viewed as a $\RelSpinC(Y,L-K_1)$-filtered chain
complex, via the filling map using the distinguished component
$K_1$ of $L$.  The filtered chain homotopy type of this complex is
identified with $\CFLa(\orL-K_1)\otimes M$.
\end{prop}

In more elementary terms, a $\Z^{\ell}$-filtered chain complex
can be viewed as a $\Z^{\ell-1}$-filtered one, by forgetting the first
term in the relative filtration. The above proposition says
that for the relatively $\Z^{\ell}$-filtered chain complex $\CFLa(Y,\orL)$,
if we forget the first term in the relative filtration, we obtain
the relatively $\Z^{\ell-1}$-filtered complex
$\CFLa(\orL-K_1)\otimes M$. 

In particular, if we consider the higher differentials inducing a
chain complex on $\HFLa(\orL)$ (from Theorem~\ref{thm:FilterHFLa}) and
take their homology in the first component (in the sense of
Subsection~\ref{subsec:HomologicalProjection}), we obtain
$\HFLa(\orL')\otimes M$ (up to some overall translation of the grading
by relative $\SpinC$ structures).

\vskip.2cm
\noindent{\bf{Proof of Proposition~\ref{prop:Forgetfuls}.}}
We begin with the relative $\Z^{\ell}$-filtered statement.

In view of Propositions~\ref{prop:StabilizationInvarianceOne}
and~\ref{prop:StabilizationInvarianceTwo}, and keeping track of
filtrations, we see that the $\Z^{\ell-1}$-filtered chain complex
obtained from $\CFLm(S^3,\orL)$ obtained by forgetting about the first
factor is identified with the mapping cone of $$
\begin{CD}
\CFLm(S^3,\orL')[U_1] @>{U_1-U_i}>> \CFLm(S^3,\orL')[U_1],
\end{CD}
$$
where $\orL'=\orL-K_1$ and $i$ is some integer $>1$.  This identification is gotten by reducing
to the case of simple stabilizations (in the sense of Subsection~\ref{subsec:SimpleStabilizations}), after handlesliding across the basepoint $z_1$. 

Specializing to $U_j\equiv 0$ for all $j$,
we obtain two copies of $\CFLa(S^3,\orL')$, verifying the simplified
statement of Proposition~\ref{prop:Forgetfuls}, using only 
the relative $\Z^{\ell}$-filtered statement. 

To obtain the version stated here, we must add the following
observations. If $(\Sigma,\alphas,\betas,\ws,\zs)$ is a
multiply-pointed Heegaard diagram for $\orL$, and $\x\in\Ta\cap\Tb$,
then the map obtained by handlesliding across $z_1$ still preserves
$\Fill{K_1}(\relspinc_{\ws,\zs}(\x))$. Moreover, in the identification
of Proposition~\ref{prop:StabilizationInvarianceOne}
between the complex obtained from applying the forgetful map to
$\CFLa(S^3,L)$ and the complex
$\CFLa(S^3,L-K_1)\otimes M$, we have an identification between the
induced relative $\SpinC$ structures.
\qed \vskip.2cm

\section{Symmetry}
\label{sec:Symmetry}

We consider the symmetry properties of $\HFLa(\orL)$. The properties
we give here are formally analogous to corresponding properties of 
knot Floer homology, cf.~\cite{Knots}, \cite{RasmussenThesis}
(compare also~\cite{HolDisk}).

\begin{prop}
  \label{prop:ReverseLink}
        Let $\orL$ be an oriented link. Then,
        we have an identification
        $$\HFLa_{*}(\orL,\relspinc)\cong \HFLa_{*}(-\orL,J{\overline\relspinc})$$
\end{prop}

\begin{proof}
  Suppose $(\Sigma,\alphas,\betas,\ws,\zs)$ is a pointed Heegaard
  diagram representing $\orL$.  Then,
  $(-\Sigma,\betas,\alphas,\ws,\zs)$ is also a pointed Heegaard
  diagram for the link $L$, only now it represents the opposite
  orientation $-\orL$. The two Heegaard diagrams give assignments
  \begin{eqnarray*}
    \spincrel_{\ws,\zs}\colon \Ta\cap\Tb\longrightarrow \RelSpinC(Y,L)
    &{\text{and}}&
    \spincrel_{\ws,\zs}'\colon \Ta\cap\Tb\longrightarrow \RelSpinC(Y,L)
  \end{eqnarray*}
  respectively.  Let $\CFaa$ denote the chain complex whose homology is
  $\HFa$, i.e. this is the associated graded object for the link
  filtration of $\CFa$. Explicitly, its generators are
  $\x\in\Ta\cap\Tb$, and it is endowed with the differential from
  Equation~\eqref{eq:DefDaa}.
  
  Of course, the generators for $\CFa(\Sigma,\alphas,\betas,\ws,\zs)$
  and $\CFa(-\Sigma,\betas,\alphas,\ws,\zs)$ coincide. Indeed, by
  pre-composing with a reflection on the disk, we see that the
  differentials for the two chain complexes coincide. 
  
  However, $\relspinc_{\ws,\zs}(\x)=J\relspinc_{\ws,\zs}'(\x)$, as in
  Lemma~\ref{lemma:JAction}.  The result follows.
\end{proof}

\begin{prop}
  \label{prop:JSymmetry}
  Let $\orL$ be an oriented link in $S^3$.  Let $f\colon
  \RelSpinC(Y,L) \longrightarrow \Z$ be the map defined as follows.
  Writing $$c_1(\spincrel)=\sum_{i=1}^\ell a_i\cm \PD[\mu_i],$$
  we have that $f(\spincrel)=\sum_{i=1}^\ell a_i$.
  Then, we have an identification
  $$\HFLa_{d}(\orL,\relspinc)\cong
  \HFLa_{d-f(\relspinc)+\ell}(\orL,J({\relspinc})+\sum_{i=1}^\ell \PD[\mu_i])$$
\end{prop}

\begin{proof}
  Write $\mu=\sum_{i=1}^\ell \mu_i$.
  As in Lemma~\ref{lemma:JAction}, we see that
  $(\Sigma,\alphas,\betas,\ws,\zs)$ and
  $(\Sigma,\alphas,\betas,\zs,\ws)$ are Heegaard diagrams for $\orL$
  and $-\orL$ respectively. The generators of $\CFaa$ and the
  differentials coincide, but the two maps to relative $\SpinC$
  structures differ, with
  $\spincrel_{\ws,\zs}(\x)=\spincrel_{\zs,\ws}(\x)+\PD[\mu]$. Thus, we see that there is a function $g\colon
  \RelSpinC(Y,L) \longrightarrow \Z$ with the property that

  $$\HFLa_{d}(\orL,\relspinc)\cong\HFLa_{d+g(\relspinc)}\left(-\orL,\relspinc-\mu\right).$$
  Suppose that $f(\relspinc_1)-f(\relspinc_2)=k$, then there is a
  flowline $\phi\in\pi_2(\x,\y)$ with $\relspinc(\x)=\relspinc_1$ and
  $\relspinc(\y)=\relspinc_2$, $n_{\ws}(\phi)=0$, $n_{\zs}(\phi)=k$.
  Now, in $(\Sigma,\alphas,\betas,\ws,\zs)$, it follows that
  $\gr(\x)-\gr(\y)=\Mas(\phi)$, whereas for
  $(\Sigma,\alphas,\betas,\zs,\ws)$, $\gr'(\x)-\gr'(\y)=\Mas(\phi)-2k$
  (where here we use $\gr'$ to denote the grading induced on the second diagram).
  It follows at once that there is some constant $c=c(\orL)$
  (depending on the oriented link) with the property that for all
  $\relspinc\in\RelSpinC(Y,L)$, $g(\relspinc)=-f(\relspinc)+c(\orL)$.
  
  Combining the above remarks with Proposition~\ref{prop:ReverseLink}, we have that
  \begin{equation}
    \label{eq:aSymmetry}
    \HFLa_{d}(\orL,\relspinc)
  =
  \HFLa_{d-f(\relspinc)+c(\orL)}(\orL,J(\relspinc)-\mu)
\end{equation}

  On the other hand, two applications of Equation~\eqref{eq:aSymmetry} shows that
  $$\HFLa_{d}(\orL,\relspinc)=\HFLa_{d-f(\relspinc) +c(\orL)) - f(J
    \relspinc+\mu)+c(\orL)}(\orL,\relspinc),$$
  and hence, since $f(J(\relspinc)+\mu)=-f(\relspinc)+2n$ that
  $c(\orL)=n$.
\end{proof}

There is a natural map 
$$i\colon \SpinC(Y,L) \longrightarrow \SpinC(-Y,L),$$
since the notion of equivalence classes of nowhere vanishing vector fields
does not refer to the orientation of the ambient manifold. Note, however,
that $c_1(i(\relspinc))=-c_1(J(\relspinc))$, since the induced two-plane field
on $\vec{v}^{\perp}$ gets opposite orientations depending on the orientation of $Y$.

\begin{prop}
  Let $\orL$ be an oriented link, and let $r(\orL)$ denote its mirror. 
  Then, 
  $$\HFLa_d(\orL,\relspinc)=\HFLa^{-d+\ell-1}(r(\orL),i(\relspinc))$$
\end{prop}

\begin{proof}
  If $(\Sigma,\alphas,\betas,\ws,\zs)$ represents $\orL$, we can think of 
  $(-\Sigma,\alphas,\betas,\ws,\zs)$ as representing its mirror. It is easy to
  see that this latter chain complex is the dual complex for the former, i.e.
  $$\CFa_{d}(-\Sigma,\alphas,\betas,\ws,\zs)\cong
  \Hom(\CFa_{-d}(\Sigma,\alphas,\betas,\ws,\zs)),$$
  so that the
  assignment from the intersection points to relative $\SpinC$
  structures is intertwined with the natural map $i$ referred to
  above. The shift in grading follows from the fact that both total
  complexes for $\CFa$ have homology isomorphic to $\Wedge^*
  \Field^{\ell-1}$, graded so that its top-most term has grading $0$
  (cf. Theorem~\ref{thm:FilterHFLa}).
\end{proof}

\subsection{Notational remarks}
\label{subsec:Notation}

Relative $\SpinC$ structures are not very concrete objects, and
hence it is sometimes awkward to describe link Floer homology as filtered
by them. We could alternatively use relative two-dimensional homology,
and indeed the above results  suggest
the following convention.

For $\kappa\in H^2(Y,L)$, we could define $\HFLa(\orL,\kappa)$ to be
$\HFLa(\orL,\relspinc)$, where here $\relspinc$ is the relative
$\SpinC$ structure determined by the formula
$$c_1(\relspinc)+\sum_{i=1}^\ell\PD[\mu_i]=\kappa.$$
Given $\kappa\in H^2(Y,L)$, write
$$\kappa=\sum_{i=1}^\ell a_i\cm\PD[\mu_i].$$
We let $|\kappa|$ be the quantity
$$|\kappa|=\sum_{i=1}^\ell a_i.$$

With these conventions, then,
Proposition~\ref{prop:JSymmetry} reads
$$\HFLa_{d}(\orL,\kappa)=\HFLa_{d-|\kappa|}(\orL,-\kappa).$$

Also, according to Proposition~\ref{prop:Forgetfuls},
$$\HFLa(\orL-K_1,\kappa+\PD[K_1])$$
is gotten by introducing extra
differentials on $\bigoplus_{s\in\Z}\HFLa(\orL,\kappa+2s\cm
\PD[\mu_1])$.

It is simpler yet to consider $\HFLa(\orL)$ as graded by elements
of $\iH(L)$; i.e. given $h\in \iH(L)$, write
$\HFLa(\orL,h)=\HFLa(\orL,\relspinc)$, where 
$\relspinc$ is chosen so that
\begin{equation}
  \label{eq:DefHFLaH}
  c_1(\relspinc)+ \sum_{i=1}^{\ell}\PD[\mu_i]= 2\cm \PD[h].
\end{equation}
Now, Proposition~\ref{prop:JSymmetry} reads
\begin{equation}
  \label{eq:SymmetryH}
  \HFLa_{d}(\orL,h)=\HFLa_{d-|h|}(\orL,-h).
\end{equation}
Note that we could have
had the grading set be inside $H$, rather than $\iH$, by dropping the
factor of $2$ in Equation~\eqref{eq:DefHFLaH}. However, we
have chosen the present formula so that
$\HFLa(K,s\cm\PD[\mu])=\HFKa(K,s)$, in the case where $K$ is a
(one-component) knot.

\section{Euler characteristics}
\label{sec:Euler}

The aim of this section is to calculate the Euler characteristic
of link homology, establishing Equations~\eqref{eq:EulerHFLa} and
~\eqref{eq:EulerHFLm}.

\subsection{The case of $\HFLa$}

Let $X$ be a connected, $m$-dimensional CW complex, and let
$H=H_1(X;\Z)/\Tors$. The {\em Reidemeister-Franz torsion} of $X$ is
defined as follows (see for example~\cite{MilnorTorsion}, \cite{Turaev}). Let $\Kfrac$
denote the field of fractions of the group-ring $\Z[H]$.  Let
${\widetilde X}$ denote the covering space determined by the natural
homomorphism of $\pi_1(X,x)\longrightarrow H$. The action by $H$ gives
the CW complex of ${\widetilde X}$, 
$C_*({\widetilde X};\Z)$ a complex
of free $\Z[H]$-modules. In particular, its
differentials
$${\widetilde \partial}_i\colon C_i({\widetilde X};\Z) \longrightarrow
C_{i-1}({\widetilde X};\Z)$$ are all $\Z[U]$-equivariant.  A {\em
fundamental family of cells for ${\widetilde X}$} is a a collection of
cells for ${\widetilde X}$, each of which projects to exactly one cell
of $X$.  The Reidemeister-Franz torsion is defined to vanish if
$H_*(C_*({\widetilde X})\otimes_{\Z[H]}\Kfrac)\neq 0$.  Otherwise, we proceed as follows.
Suppose that $c$ and $c'$ are two ordered bases for some $\Kfrac$ vector space.
Then, their top exterior products differ by an element which we denote by 
$[c/c']$ (i.e. $\det(c)=[c/c']\cm \det(c')$).
Let $b_i$ be a sequence of vectors in $C_{i}({\widetilde X};\Z)$
whose image under ${\widetilde \delta}_i$ forms a basis for this image.
By our assumption, ${\widetilde \partial}_i(b_{i+1}) b_i$ is a basis for
$C_i({\widetilde X};\Z)$; let $c_i$ be a basis for $C_i({\widetilde X};\Z)$
coming from a fundamental family of cells for ${\widetilde X}$.
Then, the Reidemeister-Franz torsion 
$\Torsion(X)\in \frac{\Kfrac}{\pm H}$ is defined to be
$$\Torsion(X)=\prod_{i=0}^m [\partial_i (b_{i+1}) b_i/c_i]^{(-1)^{i+1}}.$$ 

In the case where $X$ is the complement of an $\ell$-component link
$L\subset S^3$, $H\cong \Z^{\ell}$, and $\Torsion(X)\in \Z[H]$. Indeed,
in this case the Reidemeister-Franz torsion,
which we now denote by 
$\Torsion(L)$, is related to the (multi-variable) Alexander polynomial by
the formula
$$
\Torsion(L) = \left\{
\begin{array}{ll}
\Delta_L(T)/(T-1) & {\text{if $\ell=1$}} \\
\Delta_L(T_1,...,T_{\ell}) & {\text{if $\ell>1$.}} 
\end{array}
\right.
$$
Although this defines $t(L)$ only up to an overall sign and
translation by elements of $H$, the latter indeterminacy can be 
resolved by taking the representative of $t(L)$ which satisfies the symmetry
$t(L)={\overline t(L)}$, where here $f\in\Z[H]\mapsto {\overline f}$
is the map induced by taking $h$ to $h^{-1}$. (The existence
of such a symmetric representative is standard in the theory
of torsion, cf.~\cite{MilnorTorsion}, \cite{Turaev}.)

\begin{prop}
  \label{prop:EulerHFa}
  Given a link $L$, we have an identification
  \begin{equation}
    \label{eq:EulerHFaSymm}
    \sum_{h\in\iH}
    \chi (\HFLa_*(L,h)) \cm e^{h} = \pm 
  \left(\prod_{i=1}^{\ell}
    (T_i^{\OneHalf}-T_i^{-\OneHalf})
  \right)\cm \Torsion(L).
  \end{equation}
  (Here, we think of $T_i=e^{m_i}$, where $m_i$ are meridians for $L$.)
\end{prop}

\begin{proof}
Let $(\Sigma,\alphas,\betas,\ws,\zs)$ be a $2\ell$-pointed Heegaard diagram
for $S^3$ subordinate to $L\subset S^3$.

By stabilizing if necessary, we can arrange for the following
conditions to hold:
\begin{itemize}
\item  the genus of $\Sigma$ is $g$,
with $\beta_1,...,\beta_\ell$ forming meridians for the various components of
$L$, 
\item there are circles $\alpha_1,...,\alpha_\ell$, with the property that
$\alpha_i\cap \beta_j$ for $i=1,...,\ell$ and $j=1,...,g$ is empty unless $i=j$,
in which case $\alpha_i$ and $\beta_i$ meet in a single point.
\item $w_i$ and $z_i$ can be connected by an arc $\delta_i$ 
which is disjoint from all the $\alpha_j$ with $j=1,...,g$, 
and it is also disjoint from  all the $\beta_j$ with $j\neq i$,
meeting $\beta_i$ transversally in a single point.
\end{itemize}

Let $a_i$ be the one-handle in $S^3$ corresponding to the circle $\alpha_i$,
and $b_j$ be the two-handle corresponding to the attaching circle $\beta_j$.
The zero-handles $A_1,...A_\ell$, the one-handles $a_1,...,a_{g+\ell}$,
and the two-handles $b_{\ell+1},...b_{\ell+g}$  together form a handle decomposition
for $S^3-L$. Let ${\widetilde A}_i$, ${\widetilde a}_j$, and ${\widetilde b}_k$
(with $i=1,...,\ell$, $j=1,...,g+\ell$, $k=\ell+1,...,\ell+g$) denote a collection of
lifts to the cover ${\widetilde M}$.
Note that we can label $\{{\widetilde A}_i\}_{i=1}^{\ell}$ so that
\begin{equation}
\label{eq:ChooseZeroHandles}
{\widetilde \partial} {\widetilde a}_i = (1-T_i){\widetilde A}_i.
\end{equation}
for all $i=1,...,\ell$.

Now,
$$\Torsion({\widetilde M})
= [\{{\widetilde\partial}{\widetilde a}_i\}_{i=1}^\ell
/\{{\widetilde A}_i\}]^{-1}
\cm 
[\{{{\widetilde \partial}\widetilde b_k}\}_{k=\ell+1}^{\ell+g}
 \{ {\widetilde a}_i\}_{i=1}^\ell \}/
 \{{\widetilde a}_j\}_{j=1}^{\ell+g}]
\cm [\{\widetilde b_k\}_{k=\ell+1}^{\ell+g} /\{\widetilde b_k\}_{k=\ell+1}^{\ell+g}]^{-1}.$$
Of course, 
$[\{\widetilde b_k\}_{k=\ell+1}^{\ell+g} /\{\widetilde b_k\}_{k=\ell+1}^{\ell+g}]=1$.
Moreover, it is not difficult to see that
$$[\{{{\widetilde \partial}\widetilde b_k}\}_{k=\ell+1}^{\ell+g}
 \{ {\widetilde a}_i\}_{i=1}^\ell \}/
 \{{\widetilde a}_j\}_{j=1}^{\ell+g}]
=[\{\partial b_k\}_{k=\ell+1}^{\ell+g}/\{a_k\}_{k=\ell+1}^{\ell+g}]
= \#\left({\widetilde \Ta}\cap{\widetilde \Tb}\right).$$
Finally, 
according to Equation~\eqref{eq:ChooseZeroHandles}, 
$$ [\{{\widetilde\partial}{\widetilde a}_i\}_{i=1}^\ell
/\{{\widetilde A}_i\}]
= \prod_{i=1}^{\ell}(1-T_i).$$

Observe next that for any $\relspinc\in\RelSpinC(Y,L)$, we have that
$$\sum_{\relspinc\in\RelSpinC(Y,L)} \chi (\HFLa(L,\relspinc_0+h)) e^h \dot{=}
\#({\widetilde\Ta}\cap{\widetilde\Tb}).$$
This establishes Equation
Equation~\eqref{eq:EulerHFaSymm}, up to an overall translation by
units in the Laurent polynomials in $T_1,...,T_\ell$. The
indeterminacy in translation by the various $T_i$ is resolved by
observing that both sides of the equation are symmetric under the
involution of the ring of Laurent polynomials induced by $T_i\mapsto
T_i^{-1}$ (for all $i$): the left-hand-side is invariant according to
Proposition~\ref{prop:JSymmetry}, while the right-hand-side is symmetric 
by basic properties of torsion, cf.~\cite{MilnorTorsion}.
\end{proof}

\subsection{The case of $\HFLm$}

\begin{prop}
We have an identification
$$\sum_{h\in \iH}\chi (\HFLm(L,h))  \cm  e^h \cdot{=}
\pm \Torsion(L).$$
\end{prop}

\begin{proof}
  Each generator $\x\in\HFLa(L)$ gives rise to infinitely many generators
  $$T_1^{-a_1}\cm...\cm T_\ell^{-a_\ell} 
  \cm \x$$ for all $(a_1,...,a_\ell)\geq 0$, 
	with $\Filt(e^{h}\cm \x)=\Filt(\x)-h$.
  Thus, it follows that
  $$\sum_{h\in\iH}\chi (\HFLm(L,h))\cm
  e^{h} =\left(\sum_{\{(a_1,...,a_\ell)\in \Z^\ell\big| a\geq 0\}} 
  T_1^{-a_1}\cm...\cm T_{\ell}^{-a_\ell}\right)\cm
  \left(\sum_{h\in\iH}\chi (\HFLa(L,h)) \cm
    e^{h}\right).$$
  Since
$$\sum_{\{(a_1,...,a_\ell)\in \Z^\ell\big| a\geq 0\}} 
  T_1^{-a_1}\cm...\cm T_{\ell}^{-a_\ell}=\frac{1}{\prod_{i=1}^{\ell}(1-T_i^{-1})},$$
  the proposition now follows from Proposition~\ref{prop:EulerHFa}.
\end{proof}

\section{Relationship with knot Floer homology}
\label{sec:Relate}

Our aim here is to show that certain aspects of the link Floer
homology considered here can be extracted from the knot Floer homology
of links considered in~\cite{Knots}.  Specifically, we prove
Theorem~\ref{thm:IdentifyWithLinkHomology}, which shows in particular
that the total rank of $\HFLa$ agrees with the total rank of the knot
Floer homology for an oriented link. From this, we quickly deduce the
link Floer homology of alternating links, in
Subsection~\ref{subsec:AltLinks}.

First, we briefly recall the construction from~\cite{Knots}.  Suppose
that $L$ is an $\ell$-component oriented link in $S^3$. Then, we can
construct an oriented knot in $\#^{\ell-1}(S^2\times S^1)$ as follows.
Fix $2\ell-2$ points $\{p_i,q_i\}_{i=1}^{\ell-1}$ in $\orL$ which are
paired off in such a manner that if we formally identify each $p_i$
with $q_i$ in $L$, we obtain a connected graph.  We view $p_i$ and
$q_i$ as the feet of a one-handle to attach to $S^3$. Attaching all
$\ell-1$ of these one-handles, we obtain a three-manifold
$\kappa(S^3,\{p_i,q_i\})$ which is diffeomorphic to
$\#^{\ell-1}(S^2\times S^1)$. Inside each one-handle, we can find a
band along which to perform a connected sum of the component of $\orL$
containing $p_i$ with the component containing $q_i$. We choose the
band so that the induced orientation of its boundary is compatible
with the orientation of the link (moreover the band is chosen to be
always transverse to the foliation of the one-handle by two-spheres).
Our hypotheses on the number and distribution of the distinguished
points ensures that the newly-constructed link gotten by performing
all $\ell-1$ of the connected sums is in fact a single-component knot.
We denote this link by $\kappa(\orL,\{p_i,q_i\})$ inside
$\kappa(S^3,\{p_i,q_i\})\cong \#^{\ell-1}(S^2\times S^1)$. It is not
difficult to see that the diffeomorphism type of the pair
$(\kappa(S^3,\{p_i,q_i\}),\kappa(\orL,\{p_i,q_i\}))$ depends on only
the underlying oriented link $\orL$ (cf.
\cite[Proposition~\ref{Knots:prop:LinksToKnots}]{Knots}); hence we
denote this object simply by the knot $\kappa(\orL)\subset
\#^{\ell-1}(S^2\times S^1)$.

With these preliminaries, the {\em knot Floer homology groups}
$$\HFKa(\orL)=\bigoplus_{s\in\Z}\HFKa(\orL,s)$$
of an
oriented link are defined to be the knot Floer homology
groups of $\kappa(\orL)\subset \#^{\ell-1}(S^2\times S^1)$:
$$\HFKa(\#^{\ell-1}(S^2\times S^1),\kappa(\orL))
=\bigoplus_{s\in\Z}
\HFKa(\#^{\ell-1}(S^2\times S^1),\kappa(\orL),s).$$

In our proof of Theorem~\ref{thm:IdentifyWithLinkHomology},
it will be useful to pass from a Heegaard diagram for
a link $\orL\subset S^3$ to a corresponding Heegaard diagram
for $\kappa(L)\subset \#^{\ell-1}(S^2\times S^1)$.
This is done as follows.
Fix a Heegaard diagram $(\Sigma,\alphas,\betas,\ws,\zs)$ for an
oriented, $\ell$-component link $\orL\subset S^3$. 
Write $\ws=\{w_1,...,w_\ell\}$, $\zs=\{z_1,...,z_\ell\}$.  We attach
$\ell-1$ one-handles to $\Sigma$ to obtain a new surface $\Sigma'$ as
follows. The feet of the $i^{th}$ one-handle is attached along a
neighborhood of $z_i$ and $w_{i+1}$.  Letting $w=w_1$ and $z=z_\ell$,
we get a two-pointed Heegaard diagram $(\Sigma',\alphas,\betas,w,z)$ for
$\kappa(\orL)\subset \#^{\ell-1}(S^2\times S^1)$
(compare the proof of Theorem~\ref{thm:InvarianceHFaa}).
Admissibility for $(\Sigma,\alphas,\betas,\ws,\zs)$
translates into admissibility
for $(\Sigma',\alphas,\betas,w,z)$.


\vskip.2cm
\noindent{\bf{Proof of Theorem~\ref{thm:IdentifyWithLinkHomology}.}}
There is a one-to-one correspondence between generators for
$\CFLa(\orL)$ and $\CFKa(\kappa(\orL))$. 

Collapse the relative $\Z^\ell$ filtration on $\CFLa(S^3,\orL)$ to a
relative $\Z$-filtration by
$$\Filt_o(\x)-\Filt_o(\y)=\sum_{i}(n_{z_i}(\phi)-n_{w_i}(\phi)),$$
where here $\phi\in\pi_2(\x,\y)$ is any Whitney disk in
$\Sym^{g+\ell-1}(\Sigma)$. 
Recall also that 
the $\Z$-filtration on $\CFa(\#^{\ell-1}(S^2\times S^1),\kappa(\orL))$ is
given by $$\Filt'(\x)-\Filt'(\y)=n_{z_\ell}(\phi')-n_{w_1}(\phi')$$
where here $\phi'\in\pi_2(\x,\y)$ is any Whitney disk in
$\Sym^{g+\ell-1}(\Sigma')$. 

Clearly, any Whitney disk $\phi'\in\pi_2(\x,\y)$ 
in $\Sym^{g+\ell-1}(\Sigma')$
gives rise naturally to a Whitney disk $\phi\in\pi_2(\x,\y)$    
in $\Sym^{g+\ell-1}(\Sigma)$, having the property that    
$n_{z_i}(\phi)=n_{w_{i+1}}(\phi)$ for $i=1,...\ell-1$. 
From these conditions, it is immediate to see that
$$\Filt_o(\x)=\Filt'(\x),$$
identifying one of the factors in the
bigrading for both theories.

Moreover, it is also straightforward to see that
$\Mas(\phi)=\Mas(\phi')$, and hence the relative homological gradings
of both theories are identified. A shift in the absolute homological
bigradings by $\frac{\ell-1}{2}$ can be identified in a model
calculation, according to which the top-dimensional generator
$\HFa(\#^{\ell-1}(S^2\times S^1))$ generates $\HFa(S^3)$ in the
corresponding multiply-pointed Heegaard diagram.
We can also remove the indeterminacy in the
relative $\Z$-filtration, using the fact that both theories are
enjoy the symmetry property (Proposition~\ref{prop:JSymmetry}).

Having established that both complexes have the same generators with
the same bigradings, it remains
to show that the differentials are identified; i.e.  for suitable
choices of complex structure on $\Sigma'$, the holomorphic disks
correspond, as well.

By stabilizing and renumbering curves if necessary, we can arrange for
the diagram $(\Sigma,\alphas,\betas,\ws,\zs)$ to have the following
properties:
\begin{itemize}
\item for positive, odd $i$ less than or equal to $\ell$,
$w_i$ an $z_i$ are connected by an arc $\delta_i$ which meets $\alpha_i$
in a single point but which is disjoint from all the other $\alpha_i$ and
$\beta_j$,
\item for positive, even $i$ less than or equal to $\ell$,
$w_i$ an $z_i$ are connected by an arc $\delta_i$ which meets $\beta_i$
in a single point but which is disjoint from all the other $\alpha_i$ and
$\beta_j$,
\end{itemize}
Our aim is to exhibit a complex structure on $\Sigma'$ with the
property that for all $\x,\y\in\Ta\cap\Tb$ and all
$\phi'\in\pi_2(\x,\y)$ with $n_{w_1}(\phi')=n_{z_\ell}(\phi')=0$ and
$\Mas(\phi')=1$ (i.e. these are the disks which contribute in the
differential for $\CFKa(\kappa(\orL))$), if $\ModFlow(\phi')$ is
non-empty, then we must have $n_{z_i}(\phi')=n_{w_i}(\phi')=0$ for all
$i$, as well.

The complex structure is obtained by thinking of $\Sigma'$ as obtained
from $\Sigma$ by adding one-handles 
with a long connected sum length attached connecting $z_i$ to $w_{i+1}$
(for $i=1,...,\ell$). We choose the feet of these handles
$z_i$ and $w_{i+1}$ in a family
$\{z_i(t)\}_{t\in [0,\infty)}$ and $\{w_{i+1}(t)\}_{t\in[0,\infty)}$
starting out at $z_i(0)$ and $w_{i+1}(0)$ our original $z_i$ and $w_{i+1}$
respectively, 
$z_i(t),w_{i+1}(t)\in\Sigma-\alpha_1-...-\alpha_{g+\ell}-\beta_1-...\beta_{g+\ell}$,
for odd $i$ less than $\ell$,
\[
\lim_{t\goesto\infty} z_i(t)=z_i^\infty\in \alpha_i
\]
while for even $i$, 
\[
\lim_{t\goesto\infty} w_i(t)=w_i^\infty \in \beta_i.
\]
This can be achieved by constraining $w_i(t)$ and $z_i(t)$ to be 
subarcs of $\delta_i$.

We have a two-parameter family of complex structures on $\Sigma'$,
$\Sigma'(\sigma,\tau)$, where $\sigma$ denotes the conformal structure
on the annulus in the connected sum region, and $\tau$ parameterizes
the placement of the feet of the one-handles
$w_i(\tau)$ and $z_i(\tau)$.  

Starting with our homotopy class
$\phi'\in\pi_2(\x,\y)$, thought of representing a cylindrical
flow-line in $\Sigma'$, let $\phi\in\pi_2(\x,\y)$ be the
corresponding Whitney disk in $\Sym^{g+\ell-1}(\Sigma)$. 
If for all
$\tau$, we can find sufficiently large $\sigma$ for which
$\ModFlow_{\sigma,\tau}(\phi')$ is non-empty, then we 
can take a Gromov limit (as $\sigma\goesto\infty$) 
for any choice of $\tau$. This gives a possibly broken pseudo-holomorphic
cylindrical flow-line 
${\overline u}_\tau$ representing $\phi$. 

Observe that since $n_{w_1}(\phi')=0$, it is clear that ${\overline u}_\tau$
contains no components which contain all of $\Sigma$. Moreover, it can
contain no boundary degenerations.  Specifically, each boundary degeneration
satisfies $n_{w_i}(\psi)=n_{z_i}(\psi)$ for all $i=1,...,\ell$, while a
boundary degeneration arising in this manner must also satisfy the
additional conditions that $n_{w_1}(\psi)=n_{z_{\ell}}(\psi)=0$ (as
$n_{w_1}(\phi)=n_{z_{\ell}}(\phi)=0$), forcing all
$n_{w_i}(\psi)=n_{z_i}(\psi)$. This forces $\psi=0$.

Thus, ${\overline u}_{\tau}$ represents a juxtaposition of cylindrical
flow-lines in $\Sigma$. Moreover, each component $u_\tau$ in the
Gromov compactification must satisfy the conditions that
$$\rho^{z_i(\tau)}(u_\tau)=\rho^{w_{i+1}(\tau)}(u_\tau)$$
for
$i=1,...\ell-1$, while $\rho^{w_1}(u_\tau)=\rho^{z_{\ell}}(u_\tau)=0$.
Indeed, the same remarks apply taking $\tau\longrightarrow \infty$, as well:
we obtain a possibly broken flowline ${\overline u}_\infty$,
each of whose components onsist of cylindrical flowlines satisfying
\begin{equation}
  \label{eq:PreGlued}
\rho^{z_i^\infty}(u_\infty)=\rho^{w_{i+1}^\infty}(u_\infty).
\end{equation}

Choose $i$ minimal with the property that $n_{z_i}(\phi')\neq 0$.  If
there is no such $i$, of course, then $\phi'$ represents a flowline
with all $n_{z_i}(\phi')=n_{w_i}(\phi')=0$ as desired. Otherwise, we
assume that $n_{w_i}(\phi')=0$ (returning to the other case later).
Then, some component $u_\infty$ of ${\overline u}_\infty$ has
$\rho^{z_i^\infty}(u_\infty)\neq \emptyset$. In fact, if $i$ is odd,
then $\rho^{z_i^\infty}(u_\infty)$ is contained in $\{1\}\times\R$.
However, since $w_{i+1}^\infty$ is not on one of the $\alpha$-circles,
$\rho^{w_{i+1}}(u_\infty)$ is disjoint from $\{1\}\times \R$,
contradicting Equation~\eqref{eq:PreGlued}. A symmetrical argument
applies when $i$ is even.  In the other case, where
$n_{w_i}(\phi')\neq 0$, we can instead find the minimal $j$ for which
$n_{w_j}(\phi')\neq 0$, and then apply similar reasoning  to derive a
contradiction.  \qed

\subsection{Alternating links}
\label{subsec:AltLinks}

A large class of calculations is provided by Theorem~\ref{thm:AltLink},
which follows quickly from the material established thus far, combined
with results from~\cite{AltKnots}.

\vskip.2cm
\noindent{\bf{Proof of Theorem~\ref{thm:AltLink}.}}
In~\cite[Theorem~\ref{AltKnots:thm:LinkHomology}]{AltKnots}, 
an analogue of Theorem~\ref{thm:AltLink} is established for the knot invariant
of $L$ and the (one-variable) Alexander-Conway polynomial; specifically, 
letting $\Delta_L$ be its Alexander-Conway polynomial, and  writing
$$(T^{-1/2}-T^{1/2})^{n-1}\cm \Delta_L(T) = a_0+\sum_{s>0} a_s
(T^s+T^{-s}),$$ we have that 
$$\HFKa(S^3,L,s)\cong \Field^{|a_s|}_{(s+\frac{\sigma}{2})}.$$
Combining this with  Theorem~\ref{thm:IdentifyWithLinkHomology}, the
result follows.
\qed
\vskip.2cm

\section{The K{\"u}nneth formula for connected sums}
\label{sec:Kunneth}

In this section, we study the behaviour of the link invariant under
connected sums, generalizing corresponding results for knot Floer
homology (cf.~\cite{Knots}, \cite{RasmussenThesis}). We give a variant
for the filtration of $\CFm(S^3)$, and use it to conclude a
corresponding result for $\HFLa$ (Theorem~\ref{thm:KunnethHFLa}).

Let $C$ be a $\Z^{m}$ filtered chain complex over the ring
$\Field[U_1,...,U_m]$ and $C'$ be a $\Z^n$-filtered 
chain complex over $\Field[U_1',...,U_n']$, we can form 
their tensor product
$C\otimes_{U_1\sim U_1'} C'$, thought of as a chain complex over
the polynomial ring
$\Field[U_1,...,U_m,U_1',...,U_n']/U_1\sim U_1'$.
This complex is naturally $\Z^{m+ n-1}$-filtered, by the rule
that the part of 
$C\otimes C'$ in filtration level $(a_1,a_2,...,a_m,b_2,...,b_n)$,
$\Filt(C\otimes_{U_1\simeq U_1'} C',(a_1,a_2,...,a_m,b_2,...,b_n))$,
is generated by the sum over all $a+b=a_1$ of the subcomplexes
$$\Filt(C,a,a_2,...,a_m)\otimes \Filt(C',(b,b_2,...,b_{n})).$$

Consider two links $L_1$ and $L_2$, and fix components $p\in L_1$,
$p'\in L_2$. There is a natural map
$$\#\colon \RelSpinC(L_1)\times\RelSpinC(L_2)\longrightarrow
\RelSpinC(L_1\#_{p\sim p'} L_2),$$
defined as follows. Represent
relative $\SpinC$ structures $\relspinc_1$ and $\relspinc_2$ as
nowhere vanishing vector fields $\vec{v}_1$ and $\vec{v}_2$ which are
tangent to $L_1$ and $L_2$ respectively. Deleting a sufficiently small
ball $B_1$ around $p\in L_1\subset S^3$, the vector field $\vec{v_1}$
points normal to the boundary $S^3-B_1$ at exactly two points, where
the corresponding component of $L_1$ meets $\partial B_1$. Find a
corresponding ball $B_2$ around $p'\in L_2$. We can then identify
$v_1|_{\partial B_1}$ with $v_2|_{\partial B_2}$, and hence extend the
vector field over $S^3\#S^3$, to obtain a nowhere vanishing
vector field containing $L_1\#_{p\sim p'}L_2$ as a closed orbit. This 
determines the stated map. It is easy to see that
$$ (\relspinc_1 +\PD[h_1])\#(\relspinc_2 + \PD[h_2]) = \relspinc_1\#\relspinc_2 + 
\iota_1(\PD[h_1]) + \iota_2(\PD[h_2]),$$
where here 
$$\iota_i\colon H_1(S^3-L_i) \longrightarrow  H_1(S^3-L_1\# L_2)$$
are the natural maps.

We can refine the above tensor product as follows. If $C$ is a
$\RelSpinC(L_1)$-filtered complex over $\Field[U_1,...,U_m]$
and $C'$ is a $\RelSpinC(L_2)$-filtered
one over $\Field[U_1',...,U_n']$, we can form the chain complex
complex
$$C\otimes_{U_1\simeq U_1'} C',$$
equipped with the $\RelSpinC(L_1\#L_2)$-filtration,
where
$$\Filt(C\otimes_{U_1\simeq U_1'} C',\relspinc)$$
is generated by the sum of
$\Filt(C,\relspinc_1)\otimes \Filt(C',\relspinc_2)$ over all
pairs $\relspinc_1\in\RelSpinC(Y,L_1)$ and $\relspinc_2\in\RelSpinC(L_2)$
with $\relspinc=\relspinc_1\#\relspinc_2$.

\begin{theorem} 
        \label{thm:KunnethCFLm}
        Let $L_1$ and $L_2$ be oriented links with $n_1$ and $n_2$
        components respectively. Fix also a point $p\in L_1$ and
        $p'\in L_2$, and form the connected sum $L_1\# L_2$ (near the
        two distinguished points).  Then, the chain complex
        $\CFLm(L_1\#_{p\sim p'} L_2)$ 
        is isomorphic (as a $\RelSpinC(L_1\# L_2)$-filtered chain complex
        over $\Field[U_1,...,U_{n},U_1',...,U'_m]/U_1\sim U_1'$)
        to the tensor product
        $\CFLm(L_1)\otimes_{U_1\sim U_1'}\CFLm(L_2)$, where
        here $U_1$ and $U_1'$ are variables corresponding to the
        components of $L_1$ and $L_2$ distinguished by the points $p$
        and $p'$ respectively.
\end{theorem}

\begin{prop}
\label{prop:ConnSumHeegDiag}
Let $(Y_1,L_1)$ and $(Y_2,L_2)$ are links with multi-pointed Heegaard
diagrams $(\Sigma_1,\alphas_1,\betas_1,\{w_i^1\}_{i=1}^{n_1},\{z_i^1\}_{i=1}^{n_1})$
$(\Sigma_2,\alphas_2,\betas_2,\{w_i^2,z_i^2\}_{i=1}^{n_2})$.  The
connected sum $(Y_1\#Y_2,L_1\#L_2)$, where we connect the first
component of $L_1$ with the first component of $L_2$, can be given a
Heegaard diagram 
$$(\Sigma_1\#\Sigma_2,\alphas_1\cup\alphas_2, \betas_1\cup\betas_2,
\{w_i^1\}_{i=1}^{n_1}\cup\{w_i^2\}_{i=2}^{n_2},
\{z_i^1\}_{i=2}^{n_1-1}\cup\{z_i^2\}_{i=1}^{n_2}),$$
whose underlying surface $\Sigma_1\#\Sigma_2$
is formed by the connected sum of
$\Sigma_1$ and $\Sigma_2$ along the point $z_{1}^1\in \Sigma_1$ and
$w_1^2\in\Sigma_2$. In this new diagram, the pair of reference points
$w_1^1$ and $z_1^2$ represent the connected sum 
of the first component of $L_1$ with the first component of $L_2$.
\end{prop}

\begin{proof}
Straightforward.
\end{proof}

\vskip.2cm
\noindent{\bf{Proof of Theorem~\ref{thm:KunnethCFLm}.}}
Let $\alphas_1'$, $\betas_1'$, and $\alphas_2'$ be small perturbations
of the sets of curves $\alphas_1$, $\betas_1$, and $\alphas_2$
respectively. Correspondingly, there are canonical top-dimensional
intersection points $\Theta_1\in {\Tb}_1\cap {\mathbb T}_{\beta_1'}$ and
$\Theta_2\in {\mathbb T}_{\alpha_2}\cap{\mathbb T}_{\alpha_2'}$.

We have a map
$$\CFm(\alphas_1,\betas_1)\longrightarrow
\CFm(\alphas_1\cup\alphas_2,\betas_1\cup\alphas_2')$$
$$\x_1\mapsto \x_1\times \Theta_2.$$ 
Combining Propositions~\ref{prop:StabilizationInvarianceOne} and
\ref{prop:StabilizationInvarianceTwo}, we see that this  is a chain map.
There is a similar map
$$\CFm(\alphas_2',\betas_2) \longrightarrow
\CFm(\betas_1\cup\alphas_2',\betas_1'\cup\betas_2)$$
defined by 
$$\x_2\mapsto \Theta_1\times \x_2$$
Finally, we have a chain map 
$$\CFm(\alphas_1\cup\alphas_2,\betas_1\cup\alphas_2')\otimes
\CFm(\betas_1\cup\alphas_2',\betas_1'\cup\betas_2)\longrightarrow
\CFm(\alphas_1\cup\alphas_2,\betas_1'\cup\betas_2)$$
defined by counting pseudo-holomorphic triangles.

Putting these together, we obtain a composite chain map
$$\Phi\colon \CFm(\alphas_1,\betas_1)\otimes \CFm(\alphas_2',\betas_2)\longrightarrow
\CFm(\alphas_1\cup\alphas_2',\betas_1'\cup\betas_2).$$
There is also a ``closest point map'' inducing an isomorphism of
modules over the polynomial algebra 
($\Field[U_1,...,U_m,U_1',...,U_{n'}']/U_1\sim U_1'$) 
$$\iota\colon 
\CFm(\alphas_1,\betas_1)\otimes
\CFm(\alphas_2',\betas_2)\longrightarrow
\CFm(\alphas_1\cup\alphas_2',\betas_1'\cup\betas_2)$$
with $\iota(\x_1\times\x_2)=\x_1'\times \x_2'$. Here
$\x_1'\in{\Ta}_1\cap{\Tb}_1'$ is the point nearest to
$\x_1\in{\Ta}_1\cap{\Tb}_1$, while $\x_2'\in {\Ta}_{_2}\cap{\Tb}_2$ is the point
closest to $\x_2\in{\Ta}_2'\cap{\Tb}_2$.  
The closest point map can alternatively be thought of as counting
holomorphic triangles with minimal area. Thus, $\Phi$ is a map of the
form $\iota$ plus lower order terms (provided that the curves
$\alphas'$ resp. $\betas'$ are sufficiently to $\alphas$ resp.
$\betas'$). In particular, $\Theta$ induces an isomorphism of chain
complexes over $\Field[U_1,...,U_m,U_1',...,U_{n'}']/U_1\sim U_1'$.

It is easy to see that 
$$\relspinc_{\ws_1,\zs_1}(\x_1)\#\relspinc_{\ws_2,\zs_2}(\x_2)=\relspinc_{\ws_1\cup\ws_2,\zs_1\cup\zs_2}(\x'_1\times\x_2').$$
Thusm the map $\iota$ induces an isomorphism on the associated
$\RelSpinC(L_1\# L_2)$-graded complex.
The other terms in $\Phi$ have lower order than $\iota$, and hence
$\Phi$ induces a $\RelSpinC(L_1\# L_2)$-filtered isomorphism, as claimed.
\qed
\vskip.2cm

\vskip.2cm
\noindent{\bf{Proof of Theorem~\ref{thm:KunnethHFLa}.}}
This theorem follows readily from Theorem~\ref{thm:KunnethCFLm},
and the principle that an isomorphism of $\Z^\ell$-filtered complexes
induces an isomorphism on the homology of its associated graded object.
\qed

\section{Examples}
\label{sec:Examples}

We turn to some sample calculations. Specifically, we will be
concerned here with calculation the filtered chain homotopy type of
$\CFLa(\orL)$ for various oriented links.  In
Subsection~\ref{subsec:AltLinksAgain}, we give a theorem which can be
used to compute this data explicitly for a two-bridge link. In
Subsection~\ref{subsec:NonAlt}, we give some calculations of this data
for the first two non-alternating links.

\subsection{Two-component alternating links}
\label{subsec:AltLinksAgain}
Although Theorem~\ref{thm:AltLink} concerns only the homology of the
associated graded object $\HFLa(L)$, in many cases, one can also use
it to deduce information about the filtered homotopy type of the full
chain complex.

We study here the filtered chain homotopy type $\CFLa$
of alternating links.
For simplicity, we consider here the case of two-component
links.

\begin{theorem}
\label{thm:TwoComponentLink}
Let $\orL$ be an oriented, two-component link
which admits a connected, alternating projection.
The filtered chain homotopy type of $\CFLa(S^3,\orL)$ is determined
by the following data:
\begin{itemize}
\item the multi-variable Alexander polynomial of $\orL$
\item the signature of $\orL$
\item the linking number of $K_1$ and $K_2$
\item the filtered chain homotopy type $\CFa(S^3,K_1)$ and $\CFa(S^3,K_2)$
        of the two components $K_1$ and $K_2$ of $L$.
\end{itemize}
\end{theorem}

In the course of the proof, we show explicitly how $\CFLa(S^3,\orL)$ is
determined by the data. First, we must set up some notions.

We say that a $\Z\oplus\Z$-filtered chain complex is {\em $E_2$-collapsed}
if it has a splitting
$C=\{C_{i,j}\}_{i,j\in\Z}$ so that its differential
$\partial$ is written has the form
$\partial = D^1 + D^2$, where here
\begin{eqnarray*}
D^1|_{C_{i,j}}=D^1_{i,j}\colon C_{i,j} \longrightarrow C_{i-1,j}
&{\text{and}}&
D^2|_{C_{i,j}}=D^2_{i,j}\colon C_{i,j} \longrightarrow C_{i,j-1}.
\end{eqnarray*}
For example, suppose that $C$ is a $\Z\oplus\Z$ filtered chain complex
which also has an internal grading $g$ (which the differential drops by one).
Suppose moreover that the filtration $(i,j)$ and the grading $g$ are related
by $i+j-g=c$ for some constant $c$. Then, $C$ is $E_2$ collapsed. 
Theorem~\ref{thm:AltLink} ensures that the link Floer homology 
of an alternating link is $E_2$-collapsed.

We give some examples of $E_2$-collapsed chain complexes.
Let $B_{(d)}$ denote the chain complex with
$$\left(B_{(d)}\right)_{i,j}=
\left\{
\begin{array}{ll}
\Field_{(d)} & {\text{if $(i,j)=(0,0)$}} \\
\Field_{(d+1)} & {\text{if $(i,j)\in \{(0,1), (1,0)\}$}} \\
\Field_{(d+2)} & {\text{if $(i,j)=(1,1)$}} \\
0 & {\text{otherwise,}}
\end{array}
\right.$$
and the property that $D^1_{1,1}$ and $D^1_{1,0}$ 
and $D^2_{1,1}$ and $D^2_{0,1}$ are field isomorphisms.

Fix an integer $\ell$, and let $V^{\ell}_{(d)}$ be the chain complex with
$$\left(V^{\ell}_{(d)}\right)_{i,j}=
\left\{
\begin{array}{ll}
\Field_{(d)} & {\text{if $(i,j)=(-j,j)$ with $j=0,...,\ell-1$}} \\
\Field_{(d-1)} & {\text{if $(i,j)=(-j-1,j)$ with $j=0,...,\ell-1$}} \\
0  & {\text{otherwise,}}
\end{array}
\right.$$
where $D=D^1+D^2$ is a sum of maps, where $D^1_{-j,j}$ are field
isomorphisms for all $i$, $D^2_{-j,j}$ is an isomorphism for all
$1\leq j \leq \ell-1$.

Similarly, let $H^\ell_{(d)}$ be the chain complex with
$$\left(H^{\ell}_{(d)}\right)_{i,j}=
\left\{
\begin{array}{ll}
\Field_{(d)} & {\text{if $(i,j)=(i,-i)$ with $i=0,...,\ell-1$}} \\
\Field_{(d-1)} & {\text{if $(i,j)=(i,-i-1)$ with $i=0,...,\ell-1$}} \\
0  & {\text{otherwise,}}
\end{array}
\right.$$
where $D=D^1+D^2$ is a sum of maps,
where $D^1_{i,-i}$ 
is an isomorphism for all $i=1,...,\ell-1$,
while $D^2_{i,-i}$ is an isomorphism for all $i$.

Note that $H_*(B)\cong H_*(V^{\ell}_{(d)})\cong H^*(H^\ell_{(d)})=0$.

There are two basic types of $E_2$-collapsed chain complexes with
non-trivial homology Let $X^{\ell}_{(d)}$ be the complex given by
$$\left(X^{\ell}_{(d)}\right)_{i,j}=
\left\{
\begin{array}{ll}
\Field_{(d)} & {\text{if $i+j=\ell$ and $i, j \geq 0$}} \\
\Field_{(d+1)} & {\text{if $i+j=\ell+1$ if $i, j >0$}} \\
0  & {\text{otherwise,}}
\end{array}
\right.$$
where $D=D^1+D^2$ is a sum of maps,
where $D^1_{i,j}$ and $D^2_{i,j}$ are isomorphisms when
$i+j=\ell+1$ and $i, j>0$, and zero otherwise.

Also, let $Y^{\ell}_{(d)}$ be the complex determined by
$$\left(Y^{\ell}_{(d)}\right)_{i,j}=
\left\{
\begin{array}{ll}
\Field_{(d)} & {\text{if $i+j=\ell$ and $i,j\geq 0$}} \\
\Field_{(d-1)} & {\text{if $i+j=\ell-1$ and $i,j\geq 0$}} \\
0  & {\text{otherwise,}}
\end{array}
\right.$$
where $D=D^1+D^2$ is a sum of maps, where $D^1_{i,j}$ is 
an isomorphism when $i>0$, $j\geq 0$, and $i+j=\ell$,
while $D^2_{i,-i}$ is an 
isomorphism when $i\geq 0$, $j>0$, and $i+j=\ell$, and
all other maps are zero.

For these latter two complexes, we see that $H_*(X^\ell_{(d)})\cong
H_*(Y^\ell_{(d)})\cong\Field_{(d)}$.

It is straightforward to see that for each $E_2$-collapsed chain
complex $C$ is filtered chain homotopy equivalent to 
a filtered chain complex which splits as 
a direct sum of copies of chain
complexes of the type $B_{(d)}[i,j]$, $V^{\ell}_{(d)}[i,j]$
$H^{\ell}_{(d)}[i,j]$, $X^{\ell}_{(d)}[i,j]$, and
$Y^{\ell}_{(d)}[i,j]$, allowing $\ell$, $i$, and $j$ to vary.
(Here, as in Section~\ref{sec:Algebra}, 
given $A$ a $\Z\oplus \Z$-filtered chain complex, we let $A[i,j]$ denote
the $\Z\oplus\Z$ filtered chain complex obtained from $A$ by shifting
the filtration, so that $(A[i,j])_{x,y}=A_{x+i,y+j}$.).

\vskip.2cm
\noindent{\bf{Proof of Theorem~\ref{thm:TwoComponentLink}.}}
In view of Theorem~\ref{thm:AltLink}, $\CFa(S^3,\orL)$ is an
$E_2$-collapsed chain complex, and as such decomposes into summands of
the above five types. Thus, the filtered chain homotopy type is
determined by the number of summands of each type, and their various
parameters ($i$, $j$, $d$, and in four cases, $\ell$). Our goal is to show
how the data assumed in the statement of
Theorem~\ref{thm:TwoComponentLink} can be used to extract all of the
needed information.

Let $E^{\ell}_{(d)}$ denote the $\Z$-filtered chain complex
$$\left(E^{\ell}_{(d)}\right)_{i}
=\left\{\begin{array}{ll}
        \Field_{(d)} & {\text{if $i=0$}} \\
        \Field_{(d-1)} & {\text{if $i=-\ell$}} \\
        0 & {\text{otherwise,}}
\end{array}
\right.$$
endowed with a differential which is an isomorphism from 
$\left(E^{\ell}_{(d)}\right)_{0}$ to 
$\left(E^{\ell}_{(d)}\right)_{-\ell}$. It is easy to see that
each $\Z$-filtered chain complex splits as a sum of complexes of
the form $E^{\ell}_{(d)}[i]$, and also homologically non-trivial complexes of 
the form $\Field$ supported in some fixed degree and filtration level.

Observe that taking the homology in the vertical direction (i.e.
taking the homology with respect to $D^2$ to obtain a $\Z$-filtered
chain complex, endowed with the differential induced from $D^1$) has
the property that it annihilates $B[i,j]$, $H^{\ell}[i,j]$.  Moreover,
taking the vertical homology of $V^{\ell}_{(d)}[i,j]$, we obtain 
$E^{\ell}_{(d)}[i]$. The vertical homology of $X^{\ell}[i,j]$
is a copy of $\Field$ in filtration level $i$, while the horizontal homology
gives $\Field$ in filtration level $j-\ell$; the horizontal homology of
$Y^{\ell}[i,j]$ is a copy of $\Field$ in filtration level $j$
while its vertical homology gives $\Field$ in filtration level $i+\ell$.

On the other hand, according to Proposition~\ref{prop:Forgetfuls}, if
we take the vertical homology of $\CFLa(\orL)$, we obtain the knot
filtration on $\CFa(S^3)$ coming from $K_1$, tensored with
$\Field\oplus \Field$ (supported in two consecutive dimensions) and
shifted over (in its filtration) by $n$, 
and taking the horizontal homology, we obtain the
knot filtration of $K_2$ again shifted by $\frac{n}{2}$ and
tensored with $\Field\oplus\Field$.
Specifically, these two projections are $\Z$-filtered chain complexes,
and as such can be decomposed into one-step complexes $\Field$
(supported in some filtration level), and two-step complexes
$E^{\ell}_{(d)}[i]$, consisting of $\Field$ in filtration level $i$
and $i-\ell$ (trivial otherwise), endowed with a differential which
induces an isomorphism from the part in filtration level $i$ to the
one in filtration level $i-\ell$.

It follows from the above remarks that the summands in $\CFLa(\orL)$ of
the form $V$ are in two-to-one correspondence with the summands of
type $E$ in $\CFa(K_1)$: more precisely, if $E^{\ell}_{(d)}[i]$
appears in $\CFa(K_1)$, there are two summands,
$V^{\ell}_{(d)}[i-\frac{n}{2},d-i-\frac{n}{2}]$, 
and also $V^{\ell}_{(d-1)}[i-\frac{n}{2},d-1-i-\frac{n}{2}]$.
Similarly, the summands in $\CFLa(\orL)$ of the form $H$ are in
two-to-one correspondence with the summands of type $E$ in
$\CFa(K_2)$: more precisely, if $E^{\ell}_{(d)}[i]$ appears in
$\CFa(K_2)$, there are two summands, $H^{\ell}_{(d)}[i-\frac{n}{2},d-i-\frac{n}{2}]$, and also
$H^{\ell}_{(d-1)}[i-\frac{n}{2},d-1-i-\frac{n}{2}]$.

Moreover, since $H_*(C)\cong \Field_{(-1)}\oplus \Field_{(0)}$, there
can be at most two summands of type $X^{\ell}[i,j]$ or
$Y^{\ell}[i,j]$.  We claim that there are in fact only two possible
cases: either the two summands are $X^{\ell}_{0}[i,j]$ and
$X^{\ell-1}_{-1}[i,j]$ or $Y^{\ell}_{-1}[i,j]$ and
$Y^{\ell-1}_{0}[i-1,j-1]$. This follows from the constraints that
$H_*(C)\cong \Field\oplus\Field$ are supported in two consecutive
degrees ($0$ and $-1$), together with the constraint that if we take
the horizontal resp. vertical homologies, we are to get a filtered
chain complex whose homology $\Field\oplus\Field$ is supported in two
consecutive gradings and the same filtration level. In fact, taking
the vertical resp. horizontal homology gives the knot Floer homology
of $K_1$ resp. $K_2$, tensored with $\Field\oplus\Field$, supported in
two consecutive dimesions and filtration level $\tau(K_1)$ resp.
$\tau(K_2)$.

More specifically, writing $\tau_1=\tau(K_1)$, $\tau_2=\tau(K_2)$, and
$n$ for the linking number of $K_1$ with $K_2$, 
we have two cases, according to the sign of
$$\ell=\tau_1+\tau_2+n+\frac{\sigma-1}{2}.$$
If $\ell\geq 0$, $\CFLa(\orL)$ has summands
$$Y^{\ell}_{(0)}[\tau_2+\frac{1-\sigma-n}{2},
              \tau_1+\frac{1-\sigma-n}{2}]
\oplus
Y^{\ell-1}_{(-1)}[\tau_2+\frac{3-\sigma-n}{2},
              \tau_1+\frac{3-\sigma-n}{2}].$$
If  $\ell\leq 0$, then $\CFLa(\orL)$ has summands
$$X^{|\ell|}_{(0)}[\tau_1+\frac{n}{2},\tau_2+\frac{n}{2}]
\oplus
X^{|\ell|-1}_{(-1)}[\tau_1+\frac{n}{2},\tau_2+\frac{n}{2}].$$

The remaining summands of type $B$ (and their precise placement) are
determined by the Alexander polynomial of $L$.
\qed
\vskip.2cm

For example, if $\orL$ is a two-component, two-bridge link,
Theorem~\ref{thm:TwoComponentLink} applies immediately to express the
knot filtration in terms of the Alexander polynomial and signatures of
$\orL$. Note in particular, that in this case, the two components
$K_1$ and $K_2$ are individually unknotted; it follows in particular
that there are no summands of type $V$ or $H$. (Any two-bridge link
can be represented by a four-pointed Heegaard diagram with genus zero,
and hence the holomorphic curve counts defining the
differential take place in the Riemann sphere. Thus,
these counts are  purely
combinatorial, compare~\cite{RasmussenTwoBridge} and
also~\cite[Section~\ref{Knots:sec:Examples}]{Knots}.)

Consider the Hopf link $\orH$. There are two orientations for $\orH$,
distinguished by the signature. We denote the two cases by
$\orH^{\pm}$, with the convention that $\sigma(\orH^{\pm})=\pm (-1)$.
Note also that $\Delta_{\orH}(X,Y)=1$, and hence $\HFLa(\orH)$
consists of four generators.  It follows at once from the above
considerations that 
\[\CFLa(\orH^+)\simeq 
Y^0_{(0)}[\OneHalf,\OneHalf]\oplus Y^{1}_{(-1)}[-\OneHalf,-\OneHalf]
\] while
\[\CFLa(\orH^-)\simeq X^{1}_{(0)}[-\OneHalf,-\OneHalf]\oplus X^0_{(-1)}{[-\OneHalf,-\OneHalf]}.\]

\begin{figure}
\mbox{\vbox{\epsfbox{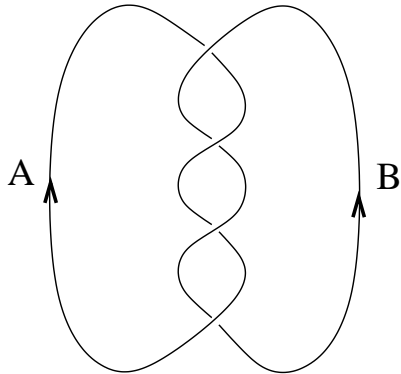}}}
\caption{\label{fig:HopfLinked}
{\bf{The oriented link $\orH_n$ with $n=2$.}}}
\end{figure}

More generally, consider the link $\orH_n$ consisting of two unknotted
circles which link each other both algebraically $n>0$ and
geometrically $n$ times -- i.e. this is the $(2,2n)$ torus link
with the specified orientation (c.f. Figure~\ref{fig:HopfLinked}).

It is easy to see that
$$\Delta_{\orH_{n}}(S,T)= S^{\frac{n-1}{2}} T^{\frac{1-n}{2}}
\sum_{i=0}^{n-1} (S^{-1}T)^i $$
and $\sigma(\orH_{n})=-n$.

It follows at once that $$\CFa(\orH_{n})\simeq
Y^0_{(0)}[\frac{n}{2},\frac{n}{2}]\oplus 
Y^{1}_{(-1)}[\frac{n}{2}-1,\frac{n}{2}-1]
\bigoplus_{i=1}^{n-1}
B_{(-2i)}[-i+\frac{n}{2},-i+\frac{n}{2}].$$

Similarly, if we reverse the orientation of one of the components, we get 
a link denoted $\orH_{-n}$, and we see that
$$\CFa(\orH_{-n})\simeq  X^{n}_{(0)}[-\frac{n}{2},-\frac{n}{2}]\oplus 
X^{n-1}_{(-1)}[-\frac{n}{2},-\frac{n}{2}].$$

\begin{figure}
\mbox{\vbox{\epsfbox{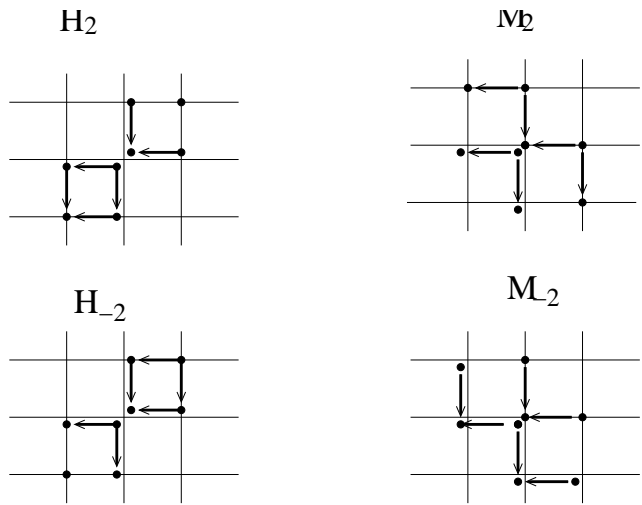}}}
\caption{\label{fig:HopfLinkedAns}
{\bf{Link Floer homologies of $H_n$ with $n=2$.}}   We have illustrated the
link Floer homologies of $H_2$ for the four different possible
orientations. The upper left is $\orH_2$ -- the right-handed, positive
version, as illustrated in Figure~\ref{fig:HopfLinked}.  The upper
right-hand illustrates the link obtained by reversing the orientation
of $B$. The second row illustrates the link Floer homology groups
of the mirrors of the links in the first row.}
\end{figure}

The complexes for $\orH_n$ and $\orH_{-n}$, along with their mirrors,
are illustrated in the case where $n=2$ 
in Figure~\ref{fig:HopfLinkedAns}.

\begin{figure}
\mbox{\vbox{\epsfbox{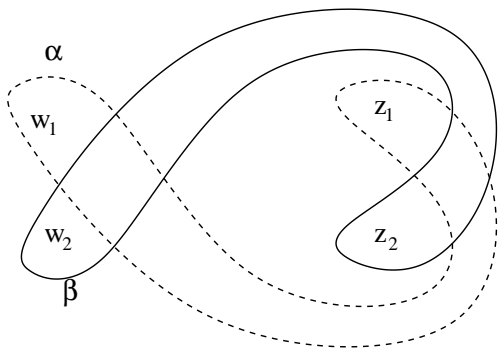}}}
\caption{\label{fig:HopfLinkedHeeg}
{\bf{Heegaard diagram for $\orH_n$ with $n=2$.}}  We have illustrated
the four basepoints and the two curves ($\alpha$ and $\beta$). This
picture takes place on the Riemann sphere.  The eight intersection
points are visible; and indeed each gives rise to a generator in
homology, as illustrated in Figure~\ref{fig:HopfLinkedAns}. The
various orientations are, of course, obtained by permuting the roles
of $w_1$ with $z_1$, while the mirrors can be found by reversing the
roles of $\alpha$ and $\beta$.}
\end{figure}

One can alternatively calculate this directly by by looking at a genus
zero Heegaard diagram, as illustrated in
Figure~\ref{fig:HopfLinkedHeeg}.

\subsection{Two non-alternating examples}
\label{subsec:NonAlt}

\begin{figure}
\mbox{\vbox{\epsfbox{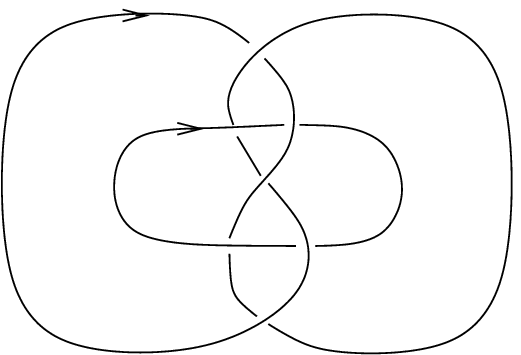}}}
\caption{\label{fig:7n2}
  {\bf{A seven-crossing non-alternating link.}}  This is the link we
  denote $L_1$; it is denoted $7^2_8$ in Rolfsen's notation; $7n_2$
  in Thistlethwaite's.}
\end{figure}

\begin{figure}
\mbox{\vbox{\epsfbox{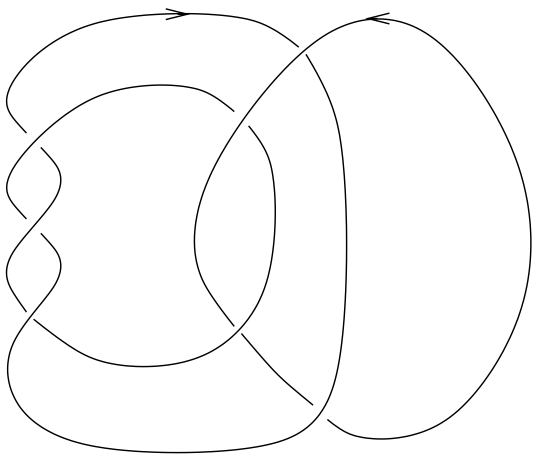}}}
\caption{\label{fig:7n1}
  {\bf{The other seven-crossing non-alternating link.}}  This is the
  link we denote $L_2$; it is denoted $7^2_7$ in Rolfsen's notation;
  $7n_1$ in Thistlethwaite's.}
\end{figure}

We calculate $\CFLa(L)$ for the first two two-component,
non-alternating links $L_1$ and $L_2$.  These are both seven-crossing
links, both obtained as a union of a trefoil and an unknot, denoted
$7^2_8$ and $7^2_7$ respectively in Roflsen's table~\cite{Rolfsen},
and $7n_2$ and $7n_1$ respectively in Thistlethwaite's link
table~\cite{ThistlethwaiteLinkTable}.  They are distinguishable
immediately by the linking number of the two components.  For $7^2_8$,
the two components have linking number zero, while for $7^2_7$ they
have linking number $\pm 2$.  Both links are illustrated in
Figures~\ref{fig:7n2} and~\ref{fig:7n1} respectively.

\subsubsection{The link $L_1$ (a.k.a. $7^2_8$, $7n2$)}

We claim that the link Floer homology groups of the link $L_1$
illustrated in Figure~\ref{fig:7n2}) have the form

\begin{equation}
\label{eq:7n2Ans}
\HFLa(7^2_8,i,j) = 
\left\{
\begin{array}{ll}
\Field_{(0)}^{4}  & {\text{if $(i,j)=(0,0)$}} \\
\Field_{(i+j)}^{2} & {\text{if $i,j\in\Z$, $|i|+|j|=1$}} \\
\Field_{(i+j)} & {\text{if $i,j\in\{\pm 1\}$}} \\
0 & {\text{otherwise.}}
\end{array}
\right.
\end{equation}
In particular, this link, too, is $E_2$-collapsed. Note that is
not necessary to separate two possible orientations: rotation
through a vertical axis gives an identification between two
possible orientations for $L_1$.

\begin{figure}
\mbox{\vbox{\epsfbox{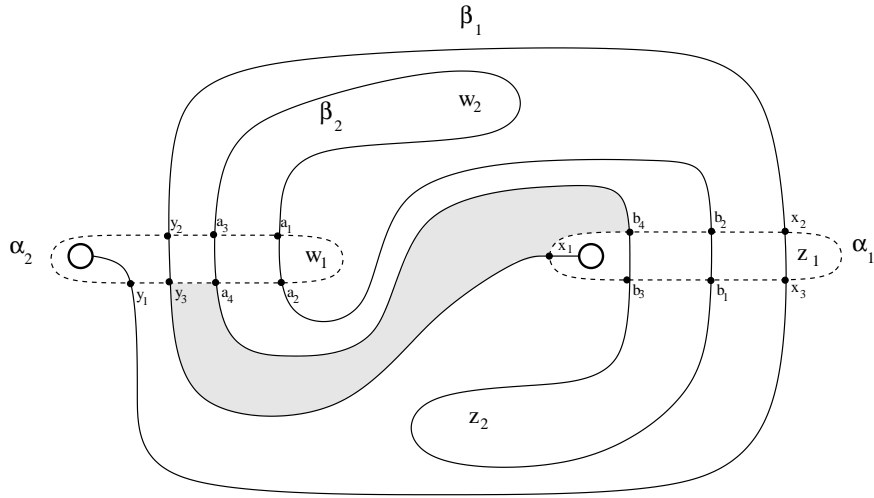}}}
\caption{\label{fig:7n2Heeg}
{\bf{Heegaard diagram for $7^2_8$.}}
Attach a one-handle to the plane at the two circles. This 
gives a Heegaard diagram for $7^2_8$, with the orientation as
given in Figure~\ref{fig:7n2}. The  shaded rectangle
represents a flow-line from $a_4\times x_1$ to $b_4\times y_3$.}
\end{figure}

Equation~\eqref{eq:7n2Ans} can be seen by considering the (admissible)
Heegaard diagram pictured in Figure~\ref{fig:7n2Heeg}, which takes
place in a genus one surface, with attaching circles
$\{\alpha_1,\alpha_2\}$ and $\{\beta_1,\beta_2\}$.

We label
\begin{eqnarray*}
\alpha_1\cap\beta_1=\{x_1,x_2,x_3\} &&
\alpha_1\cap\beta_2=\{b_1,b_2,b_3,b_4\} \\
\alpha_2\cap\beta_1=\{y_1,y_2,y_3\} 
&&\alpha_2\cap\beta_2=\{a_1,a_2,a_3,a_4\}
\end{eqnarray*}
Any two intersection points between $\alpha_i$ and $\beta_j$ can be
connected by a sequence of consecutive embedded Whitney disks. Using
these disks (each of which has $\Mas(\phi)=\pm 1$, depending on its
orientation), it is straightforward to calculate the relative gradings and
filtrations of any two intersection points of the form $\{a_i,x_j\}$
(or any two intersection points of the form $\{b_i,y_j\}$). Finally, 
there is a square at the center of the diagram which represents a Whitney
disk connecting $a_4\times x_1$ to $b_4\times y_3$, which allows one to complete
the calculations all relative gradings and filtration levels of generators.

Indeed, one readily sees that in all but four filtration levels (which
in the normalization of Equation~\eqref{eq:7n2Ans}, are $(i,j)\in
\{-1,-2\}\times \{0,-1\}$) the generators have the same relative
gradings, and hence the ranks and degrees of the homology (up to an
overall translation) are as given in Equation~\eqref{eq:7n2Ans}.

We consider the four remaining filtration levels $(-2,-1)$, $(-2,0)$,
$(-1,-1)$, $(-1,0)$, represented by intersection points
\[
\begin{array}{ll}
\{a_4\times x_1, b_4\times y_3\}, &
\{a_2\times x_1, b_2\times y_3\}, \\
\{a_3\times x_1,b_3\times y_3, b_4\times y_2\}, &
\{a_1\times x_1, a_4\times x_2, b_1\times y_3, b_2\times y_2\}
\end{array}
\]
respectively.
We need to show that the homology groups have ranks $0$, $0$, $1$ and $2$ 
respectively.

To handle the filtration level $(-2,-1)$, one inspects the Heegaard
diagram in Figure~\ref{fig:7n2Heeg}, to find a rectangle connecting
$a_4\times x_1$ to $b_4\times y_3$.  It is easy to see that there are
no other non-negative domains connecting these two disks, hence the
differential annihilates this pair of generators. (For this, it is
useful to observe that all other possible homology classes of disks
which could contribute to the differential are obtained from this
given square by the addition of a periodic domain. But any non-trivial
periodic domain has both positive and negative local multiplicities in
one of the regions adjoining $b_1$, which is disjoint from our given
square. Thus, there are no other non-negative homotopy classes.)

In the same manner, one can find a rectangle to show that the homology
in the filtration level $(-2,0)$ is trivial. Indeed, one can find also
rectangles connecting generators $a_3\times x_1$ to $b_4 \times y_2$
and $a_1\times x_1$ to $b_2\times y_2$, showing that the differentials
in filtrations levels $(-1,-1)$ and $(-1,0)$ are non-trivial. It is
straightforward then to conclude that the groups with their relative
gradings are as given in Equation~\eqref{eq:7n2Ans}. 

To verify the absolute gradings, recall that we need to orient the
knot (and an orientation is implicit in the Heegaard diagram, via the
choices of $w_i$ and $z_j$). If we allow isotopies to cross $z_2$ and
$z_1$, it is easy to shrink $\beta_2$ in Figure~\ref{fig:7n2Heeg}
(canceling out intersection points $b_1$, $b_3$, $b_2$, $b_4$, $a_4$,
$a_2$) and then perform a finger move on $\beta_1$ to cancel point
$x_2$ and $x_3$, to obtain an admissible doubly-pointed Heegaard
diagram for $S^3$ with exactly two generators, $a_1\times x_1$ and
$a_3\times x_1$. It is easy to see that $a_1\times x_1$ represents a
generator for $\HFa(S^3)$, and hence it is supported in degree zero.
This completes the verification of  Equation~\eqref{eq:7n2Ans}.

\begin{figure}
\mbox{\vbox{\epsfbox{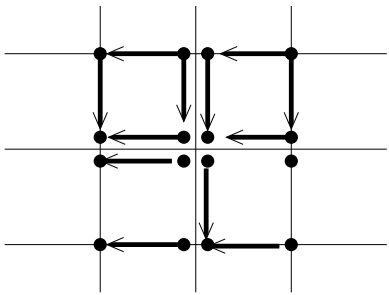}}}
\caption{\label{fig:7n2Ans}
{\bf{Link Floer homology for $L_1$.}}
We have illustrated the link Floer homology for the link $L_1$
as pictured in Figure~\ref{fig:7n2}. The first coordinate
denotes the filtration induced by the trefoil component,
while the second denotes the filtration by the unknot component.}
\end{figure}

As in the alternating case, the higher differentials on $\HFLa(S^3,\orL)$
can be determined by the knot Floer homologies of the two components
(the trefoil and the unknot). 
We illustrate the resulting complex in Figure~\ref{fig:7n2Ans}.

\subsubsection{The link $L_2$ (a.k.a. $7^2_7$, $7n1$)}

The Floer homology of the link $\orL_2$ is not $E_2$-collapsed.  In
fact, we have:
\begin{equation}
\label{eq:7n1Ans}
  \HFLa(\orL_2)=
\left\{
\begin{array}{ll}
  \Field_{(i+j-3)} & {\text{if $(i,j)\in \{0,1\}\times \{1,2\}
    \cup \{0,-1\}\times \{-1,-2\}$}} \\
  \Field_{(-2)}\oplus\Field_{(-3)} & {\text{if $(i,j)=(0,0)$}} \\
  0 & {\text{otherwise.}}
\end{array}
\right.
\end{equation}

\begin{figure}
\mbox{\vbox{\epsfbox{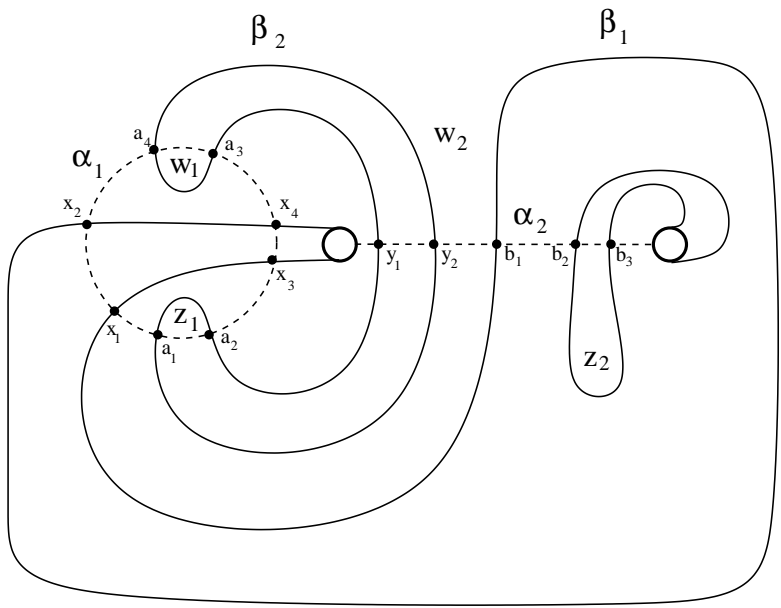}}}
\caption{\label{fig:7n1Heeg}
{\bf{Heegaard diagram for $7^2_7$.}}
Attach a one-handle to the plane at the two circles. This 
gives a Heegaard diagram for $7^2_7$, with the orientation as
given in Figure~\ref{fig:7n1}.}
\end{figure}

To perform the calculation, we draw the Heegaard diagram in
Figure~\ref{fig:7n1Heeg}.
This once again is a genus one Heegaard diagram
with two pairs of attaching circles $\{\alpha_1,\alpha_2\}$ and
$\{\beta_1,\beta_2\}$. Now we label
\begin{eqnarray*}
\alpha_1\cap\beta_1=\{x_1,x_2,x_3,x_4\},
&&\alpha_1\cap\beta_2=\{a_1,a_2,a_3,a_4\} \\
\alpha_2\cap\beta_1=\{b_1,b_2,b_3\},
&&\alpha_2\cap\beta_2=\{y_1,y_2\}
\end{eqnarray*}
Proceeding as before, it is easy to calculate filtration levels and
gradings (up to an overall shift) of generators. There are three
filtration levels, $(-1,1)$, $(-1,0)$, and $(1,0)$, where there are
exactly two generators, one, $(0,1)$, where there are three, and
finally one, $(0,0)$ where there are four. In all the rest, there are
one or zero generators, so in these other filtration levels, 
Equation~\eqref{eq:7n1Ans} is immediately verified.

For $(i,j)=(-1,1)$, there are two generators $\{a_3\times b_3,
x_4\times y_1\}$, and a rectangle can be immediately found to show
that the homology in this filtration level is trivial. The same
remarks apply for $(i,j)=(-1,0)$ and $(1,0)$.  Also, rectangles can be
found connecting $a_4\times b_3$ to $x_2\times y_1$ (and also
$a_4\times b_3$ to $x_4\times y_2$), forcing the differential to be
non-trivial and hence the homology to be one-dimensional.

The remaining case of $(i,j)=(0,0)$ can be either analyzed carefully
in this way, or alternatively, one can argue that taking the
horizontal homologies should give the knot Floer homology of the
trefoil (as in Proposition~\ref{prop:Forgetfuls}) tensored with a
two-dimensional vector space. In particular, the horizontal homology
through the $j=0$ line must be two-dimensional. But since the homology
is trivial for all $i\neq 0$ and $j=0$, it follows that this
horizontal homology is identified with the homology in the filtration
level $(i,j)=(0,0)$. 

Absolute degrees are now computed by considering isotopies which cross
$z_1$ and $z_2$. It is easy to find an isotopy of $\beta_1$ crossing
$z_2$ which cancels $b_1,b_2, x_1,x_2,x_3,_x4$, and then an isotopy of
$\beta_2$ which crosses $z_1$ cancelling $a_1$ and $a_2$. This leaves
a diagram with exactly two generators $a_3\times b_1$ and $a_4\times
b_1$. In fact, $a_4\times b_1$ is the generator which survives in
$\HFa(S^3)$, and hence it must have absolute grading equal to zero.

This completes the verification of Equation~\eqref{eq:7n1Ans}.

\begin{figure}
\mbox{\vbox{\epsfbox{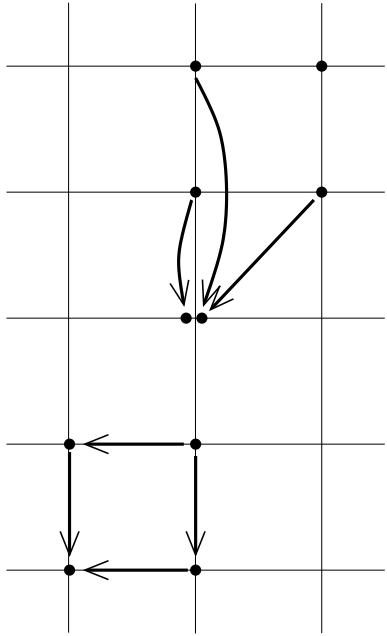}}}
\caption{\label{fig:7n1Ans1}
{\bf{Heegaard diagram for $7^2_7$.}}
Link Floer homology for $L_2$ with the orientation from
Figure~\ref{fig:7n1}. The first coordinate represents the unknot factor,
while the second represents the trefoil.}
\end{figure}

To compute higher differentials, we can no longer argue that the link
complex must be $E_2$-collapsed. However, we still have a number of
constraints: the two homological projections give the homology of the
trefoil and the unknot respectively, and the 
total homology must have rank two (in dimensions $0$ and $1$). These 
constraints suffice to determine the higher differentials uniquely.
We have illustrated these in Figure~\ref{fig:7n1Ans1}.

\begin{figure}
\mbox{\vbox{\epsfbox{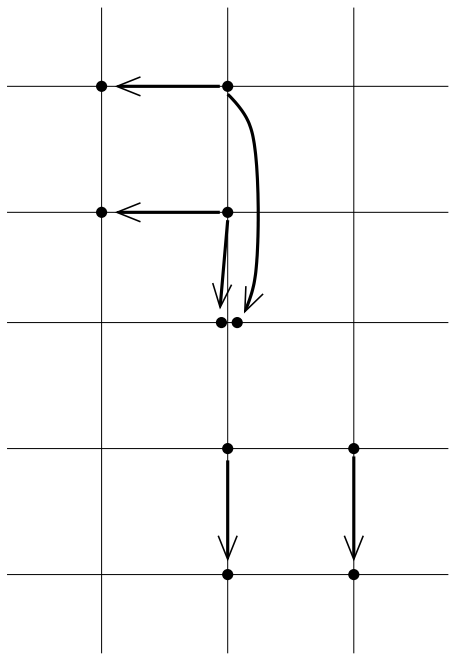}}}
\caption{\label{fig:7n1Ans2}
{\bf{Heegaard diagram for $7^2_7$.}}
Link Floer homology for $7^2_7$, reversing  the orientation
of  with the orientation of the unknot from that depicted in
Figure~\ref{fig:7n1}.}
\end{figure}

Reverse the orientation of one of the two components. Again, we obtain a 
chain complex which is uniquely determined by these constraints. We have
illustrated the unique solution in 
Figure~\ref{fig:7n1Ans2}.

\bibliographystyle{plain}
\bibliography{biblio}

\begin{thebibliography}{10}

\bibitem{BEHWZ}
F.~Bourgeois, Y.~Eliashberg, H.~Hofer, K.~Wysocki, and E.~Zehnder.
\newblock Compactness results in symplectic field theory.
\newblock {\em Geom. Topol.}, 7:799--888 (electronic), 2003.

\bibitem{EGH}
Y.~Eliashberg, A.~Givental, and H.~Hofer.
\newblock Introduction to symplectic field theory.
\newblock {\em Geom. Funct. Anal.}, (Special Volume, Part II):560--673, 2000.
\newblock GAFA 2000 (Tel Aviv, 1999).

\bibitem{FloerLag}
A.~Floer.
\newblock Morse theory for {L}agrangian intersections.
\newblock {\em J. Differential Geometry}, 28:513--547, 1988.

\bibitem{FOOO}
K.~Fukaya, Y-G. Oh, K.~Ono, and H.~Ohta.
\newblock {\em Lagrangian intersection Floer theory---anomaly and obstruction}.
\newblock Kyoto University, 2000.

\bibitem{Gornik}
B.~Gornik.
\newblock Note on {K}hovanov link cohomology.
\newblock math.QA/0402266, 2004.

\bibitem{Gromov}
M.~Gromov.
\newblock Pseudo holomorphic curves in symplectic manifolds.
\newblock {\em Inventiones mathematicae}, 82:307--347, 1985.

\bibitem{IonelParker}
E-N. Ionel and T.~H. Parker.
\newblock Relative {G}romov-{W}itten invariants.
\newblock {\em Ann. of Math. (2)}, 157(1):45--96, 2003.

\bibitem{Khovanov}
M.~Khovanov.
\newblock A categorification of the {J}ones polynomial.
\newblock {\em Duke Math. J.}, 101(3):359--426, 2000.

\bibitem{KhovanovRozansky}
M.~Khovanov and L.~Rozansky.
\newblock Matrix factorizations and link homology.
\newblock math.QA/0401268.

\bibitem{KhovanovRozanskyII}
M.~Khovanov and L.~Rozansky.
\newblock Matrix factorizations and link homology {II}.
\newblock math.QA/0505056.

\bibitem{KMcontact}
P.~B. Kronheimer and T.~S. Mrowka.
\newblock Monopoles and contact structures.
\newblock {\em Invent. Math.}, 130(2):209--255, 1997.

\bibitem{EunSooLee}
E.~S. Lee.
\newblock The support of the {K}hovanov's invariants for alternating knots.
\newblock math.GT/0201105, 2002.

\bibitem{LiRuan}
A-M. Li and Y.~Ruan.
\newblock Symplectic surgery and {G}romov-{W}itten invariants of {C}alabi-{Y}au
  3-folds.
\newblock {\em Invent. Math.}, 145(1):151--218, 2001.

\bibitem{LipshitzCyl}
R.~Lipshitz.
\newblock A cylindrical reformulation of {H}eegaard {F}loer homology.
\newblock math.SG/0502404, 2005.

\bibitem{McDuffSalamon}
D.~McDuff and D.~Salamon.
\newblock {\em {$J$}-holomorphic curves and quantum cohomology}.
\newblock Number~6 in University Lecture Series. American Mathematical Society,
  1994.

\bibitem{MilnorTorsion}
J.~Milnor.
\newblock Whitehead torsion.
\newblock {\em Bull. Amer. Math. Soc.}, 72:358--426, 1966.

\bibitem{DunfieldGukovRasmussen}
J.~Rasmussen N.~M.~Dunfield, S.~Gukov.
\newblock The {S}uperpolynomial for {K}not {H}omologies.
\newblock math.GT/0505662.

\bibitem{AltKnots}
P.~S. Ozsv{\'a}th and Z.~Szab{\'o}.
\newblock Heegaard {F}loer homology and alternating knots.
\newblock {\em Geom. Topol.}, 7:225--254, 2003.

\bibitem{4BallGenus}
P.~S. Ozsv{\'a}th and Z.~Szab{\'o}.
\newblock Knot {F}loer homology and the four-ball genus.
\newblock {\em Geom. Topol.}, 7:615--639, 2003.

\bibitem{Knots}
P.~S. Ozsv{\'a}th and Z.~Szab{\'o}.
\newblock Holomorphic disks and knot invariants.
\newblock {\em Adv. Math.}, 186(1):58--116, 2004.

\bibitem{HolDisk}
P.~S. Ozsv{\'a}th and Z.~Szab{\'o}.
\newblock Holomorphic disks and topological invariants for closed
  three-manifolds.
\newblock {\em Ann. of Math. (2)}, 159(3):1027--1158, 2004.

\bibitem{LinkTwo}
P.~S. Ozsv{\'a}th and Z.~Szab{\'o}.
\newblock Link {F}loer homology and the {T}hurston norm.
\newblock math.GT/0601618, 2006.

\bibitem{RasmussenTwoBridge}
J.~A. Rasmussen.
\newblock Floer homologies of surgeries on two-bridge knots.
\newblock {\em Algebr. Geom. Topol.}, 2:757--789, 2002.

\bibitem{RasmussenThesis}
J.~A. Rasmussen.
\newblock {\em Floer homology and knot complements}.
\newblock PhD thesis, Harvard University, 2003.

\bibitem{RasmussenSlice}
J.~A. Rasmussen.
\newblock Khovanov homology and the slice genus.
\newblock math.GT/0402131, 2004.

\bibitem{Rolfsen}
D.~Rolfsen.
\newblock {\em Knots and links}, volume~7 of {\em Mathematics Lecture Series}.
\newblock Publish or Perish Inc., Houston, TX, 1990.
\newblock Corrected reprint of the 1976 original.

\bibitem{ThistlethwaiteLinkTable}
M.~Thistlethwaite.
\newblock Link table.
\newblock
  {http://katlas.math.toronto.edu/wiki/The\_Thistlethwaite\_Link\_Table}.

\bibitem{Turaev}
V.~Turaev.
\newblock {\em Torsions of 3-manifolds}, volume~4 of {\em Geom. Topol. Monogr.}
\newblock Geom. Topol. Publ., Coventry, 2002.

\end{thebibliography}

\end{document}